\documentclass[12pt,a4paper]{amsart}

\usepackage[american]{babel}
\usepackage{bm,amsmath,amssymb}
\usepackage{amsbsy,amsfonts,array}
\usepackage{verbatim,graphicx}
\usepackage{epstopdf}
\usepackage{subfigure}
\graphicspath{{./immagini/}}
\usepackage{color}
\DeclareGraphicsExtensions{.jpg,.png,.pdf}
\usepackage{dsfont} 
\usepackage{amsthm}
\usepackage{textcomp}                      
\usepackage[mathscr]{euscript}

\usepackage[margin=3cm]{geometry}
\usepackage[P,R]{dblfont}

\newcommand{\be}{\begin{equation*}}
\newcommand{\ee}{\end{equation*}}
\newcommand{\bcr}{\begin{center}}
\newcommand{\ecr}{\end{center}}
\newcommand{\beq}{\begin{equation}}
\newcommand{\eeq}{\end{equation}}
\newcommand{\Frac}{\displaystyle\frac}

\newcommand{\fl}{\textrm{fl}}
\newcommand{\s}{\textrm{s}}

\newcommand{\n}{\textrm{n}}
\newcommand{\q}{\textrm{q}}
\renewcommand{\v}{\textrm{v}}
\newcommand{\ECM}{\textrm{ECM}}
\newcommand{\cells}{\textrm{cells}}

\newcommand{\qui}{\textrm{qui}}

\newcommand{\cell}{\textrm{cell}}

\newcommand{\unit}[1]{\mathrm{#1}}

\definecolor{darkmagenta}{rgb}{0.55, 0.0, 0.55}
\newcommand{\rs}[1]{\textcolor{black}{#1}}

\begin{document}

\title[A Poroelastic Mixture Model, Part II]
{A Poroelastic Mixture Model of Mechanobiological
Processes in Tissue Engineering. \\
Part II: Numerical Simulations}

\author{Chiara Lelli$^{1}$ \and Riccardo Sacco$^{1}$ 
\and Paola Causin$^{2}$
\and Manuela T. Raimondi$^{3}$}

\address{$^{1}$ Dipartimento di Matematica, Politecnico di Milano, \\
Piazza Leonardo da Vinci 32, 20133 Milano, Italy \\
$^{2}$ Dipartimento di Matematica \lq\lq F. Enriques\rq\rq,
Universit\`a degli Studi di Milano, \\
via Saldini 50, 20133 Milano, Italy \\
$^{3}$ Dipartimento di Chimica, Materiali e Ingegneria
Chimica \\ \lq\lq Giulio Natta\rq\rq, Politecnico di Milano, 
Piazza Leonardo da Vinci 32, 20133 Milano, Italy}

\date{\today}

\begin{abstract}
In Part I of this article we have developed a novel
mechanobiological model of a \rs{Tissue Engineering process} 
hat accounts for the mechanisms through which an isotropic or anisotropic
adherence condition regulates the active functions of the cells
in the construct.
The model expresses mass balance and force equilibrium balance for
a multi-phase
mixture in a 3D computational domain and in time dependent conditions.
In the present Part II, we study
the mechanobiological model in a simplified 1D geometrical setting
with the purpose of highlighting the ability of the formulation to
represent the influence of force isotropy and nutrient availability on the \rs{growth} of the tissue construct. In particular,
an example of isotropy estimator is proposed and
coded within a fixed-point solution map that is used
at each discrete time level for system linearization
and subsequent finite element approximation of the linearized equations.
% Extensive simulations are conducted to investigate the temporal evolution
% of tissue production and growth under several working conditions
% characterized by rather different values of system input parameters.
% \rs{Results show that: 1) the spatial and temporal evolution of
% the various cellular populations (solid biomass components) are
% in good agrement with the local growth/production conditions predicted by
% the mechanobiological stress-dependent model; and 2)
% the isotropy indicator, which in our model is responsible of
% the biological fate of the each cellular population,
% is strongly influenced by both the maximum cell
% specific growth rate and the mechanical boundary conditions
% enforced at the interface between the biomass construct
% and the interstitial fluid.}
Extensively conducted simulations show that:
% , conducted to investigate the evolution
% of the construct under several working conditions, show that:
1) the spatial and temporal evolution of
the cellular populations are
in good agrement with the local growth/production conditions predicted by
the mechanobiological stress-dependent model; and 2)
the isotropy indicator and all
model variables are strongly influenced by both maximum cell
specific growth rate and mechanical boundary conditions
enforced at the interface between the biomass construct
and the interstitial fluid.
%solid-fluid interface.}}
% \rs{Results show that the behavior of the solid biomass components
% is coherent with the evolution of the isotropic indicator, which, in turn,
% is strongly influenced by both the maximum cell
% specific growth rate and the mechanical boundary conditions
% enforced at the interface between the biomass construct
% and the interstitial fluid.
% %solid-fluid interface.}}
\end{abstract}

\maketitle


\vspace*{5pt}

{\bf Abbreviations:} TE (tissue engineering); ECM (extracellular matrix);
ACC (articular \rs{chondrocyte} cell).

\section{Introduction}
In the mathematical model proposed and illustrated in
Part I of the present research, the physical problem of tissue growth in a scaffold-based bioreactor is described by a system of PDEs
constituted by: 1) the balance of mass
for the solid and fluid phases of the growing mixture; 2) the continuity equation for the nutrient concentration; 3) the
linear momentum balance equation for the mixture components.

In this second part we recast the mathematical picture in a simplified geometrical one-dimensional setting and we perform an accurate
numerical simulation of the engineered construct with a twofold purpose.
First, because of the complexity of the problem, we
assess the reliability of model predictions by verifying
that numerical results are biophysically reasonable and consistent with experimental measurements.
Second, we single out the presence of some critical parameters
in the mathematical formulation and investigate their role and
quantitative influence on the evolution of mixture components.

An extensive set of numerical simulations,
conducted to study the evolving construct under
different working conditions, show that:
1) the computed cellular populations are
in good agrement with the local growth/production conditions predicted by
the mechanobiological stress-dependent model; and 2)
the isotropy indicator, which in our model is responsible of
the biological fate of the each cellular population, and consequently all the other model variables are strongly influenced by both maximum cell
specific growth rate and mechanical boundary conditions
enforced at the interface between the biomass construct
and the interstitial fluid.

An outline of this Part II is as follows. In Sect.~\ref{sec:mechanobio_1D}
we describe in detail the mechanobiological model in the one-dimensional
spatial geometrical configuration. Sect.~\ref{sec:definition_of_r} is
devoted to defining the stress indicator parameter $r$ and to
providing a biophysical interpretation.
In Sect.~\ref{sec:numerical_approximation_1D}
we illustrate the computational algorithm that is used to
solve the equation system of Sect.~\ref{sec:mechanobio_1D}.
Sections~\ref{sec:simulations} and~\ref{sec:discussion}
contain the description and discussion
of the various test cases numerically studied to validate
the proposed model, while in Sect.~\ref{sec:conclusions}
\rs{we address some future research perspectives.}

\section{The mechanobiological model in 1D}\label{sec:mechanobio_1D}
In this section we formulate the mechanobiological model proposed
in Part I in a one-dimensional (1D) geometrical configuration (1D).
\rs{In the remainder of the article we shortly write BVP and IV-BVP 
to denote "boundary value problem" and 
"initial value/boundary value problem", respectively.}

\subsection{The one-dimensional computational domain}\label{sec:1D_domain}
\begin{figure}[h!]
\centerline{
\includegraphics[width=0.6\textwidth]{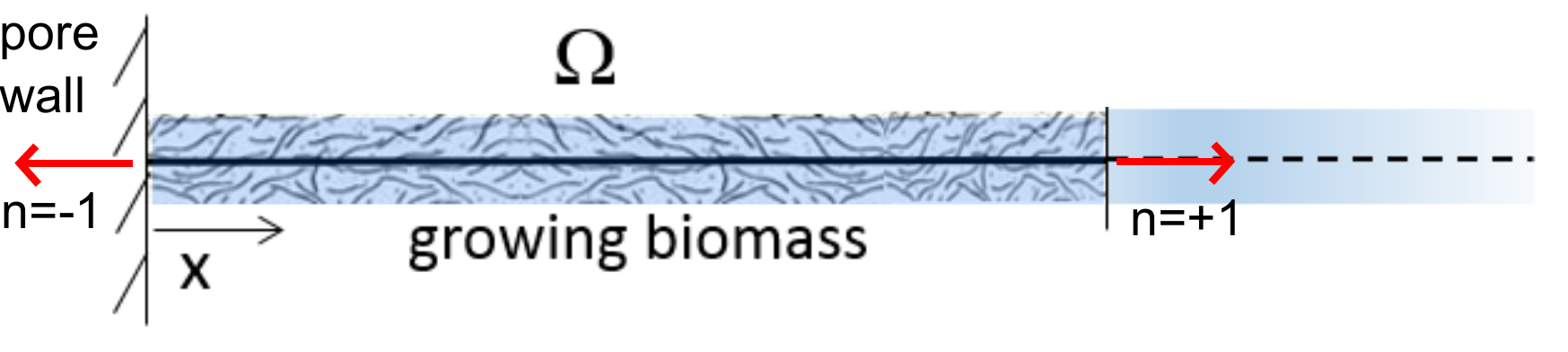}} \null
\centerline{
\includegraphics[width=0.6\textwidth]{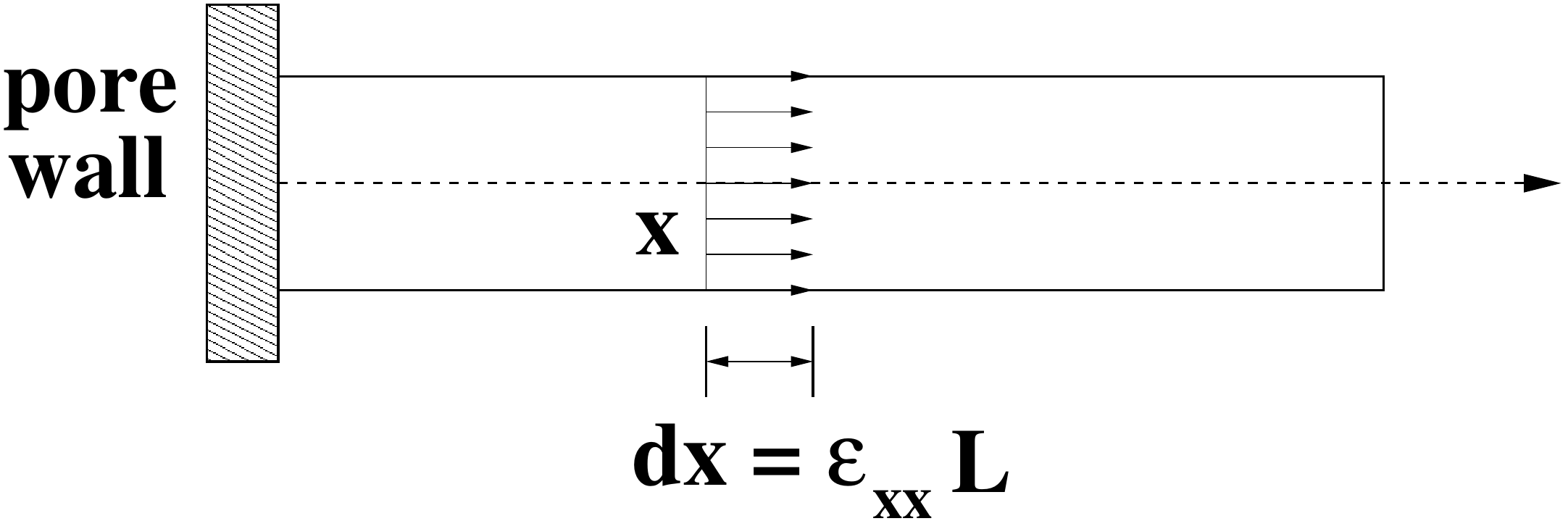}
}
\caption{Top panel: schematic representation of a 
1D scaffold-based bioreactor. Bottom panel: uniaxial stress state.}
\label{fig:1D_domain}
\end{figure}

Fig.~\ref{fig:1D_domain} \rs{(top panel)}
shows a (rather) simplified representation of
the 3D scaffolded bioreactor used in the experimental analysis discussed in~\cite{Raimondi2}. Denoting by $x$ the spatial coordinate, the
region $x < 0$ represents the scaffold wall, the open interval
$\Omega = (0, L)$ is the growing tissue
whereas the region $x > L$ corresponds to the interstitial fluid that
brings nutrient to the growing construct.
We denote by $\partial \Omega = \left\{ 0, \, L \right\}$ the boundary
of the computational domain and by $n$ the outward unit normal vector
on $\partial \Omega$. We have $n=-1$ at $x=0$ and $n=+1$ at $x=L$.

\subsection{The one-dimensional model}\label{sec:1D_model}
The geometrical reduction introduced in the previous section
implies that:
\begin{itemize}
\item all model variables depend on the sole spatial coordinate $x$
and on the time variable $t$;
\item the solid displacement field $\mathbf{u}$ has only one nonvanishing
component, that is $\mathbf{u} = [u, 0, 0]^T$ with $u=u(x,t)$;
\item the strain tensor has only one nonvanishing component,
that is
$$
\boldsymbol{\varepsilon}(\mathbf{u}(\mathbf{x},t)) =
\left[
\begin{array}{ccc}
\Frac{\partial u(x,t)}{\partial x} & 0 & 0 \\[3mm]
0 & 0 & 0 \\
0 & 0 & 0
\end{array}
\right].
$$
%A similar expression can also be written for the total stress $\mathbf{T}$.
\end{itemize}

The above analysis allows us to conclude that the simplified model
of biomass growth considered in the present Part II can be regarded
as a nonhomogeneous bar (fixed at one endpoint) subject to a
uniaxial state of mechanical stress in such a way that each point
$P = P(x,y,z)$ of every cross-section of the bar undergoes the
\emph{same} deformation 
$\rs{\varepsilon_{xx}(\mathbf{u},t)} = \partial u(x,t) /\partial x$
(see Fig.~\ref{fig:1D_domain}, bottom panel).

% At each point $x \in \Omega$ of Fig.~\ref{fig:1D_domain}
% we can identify a REV as in Fig~1 of Part I where all the components of the mixture coexist in an intermingled manner. The associated
% volumetric fractions are regulated by the 1D evolution equations
% presented below.

\subsubsection{Poroelastic \rs{IV-BVP for biomass}}
For given $\phi_\eta$ and $g_\eta$,
$\eta=\cells, \ECM, \fl$,
find the solid displacement $u:\mathcal{Q}_{T_{end}}
\rightarrow\mathbb{R}$
%,the fluid friction velocity $w:\mathcal{Q}_{T_{end}}\rightarrow\mathbb{R}$
and the pressure
$p:\mathcal{Q}_{T_{end}}\rightarrow\mathbb{R}$ that satisfy the
following system of partial differential equations in balance form:
\begin{subequations}\label{mixture_momentum_balance_1D}
\begin{align}
& \Frac{\partial T_{xx}}{\partial x} = 0 &
\label{eq:momentum_balance_1D}\\[2mm]
%(u,p,\widetilde{\boldsymbol{\phi}})}
& \Frac{\partial}{\partial t}\Frac{\partial u}{\partial x} +
\Frac{\partial V}{\partial x} = 0
& \label{eq:continuity_2_1D}\\
& V=-K(\phi_\fl)\frac{\partial p}{\partial x} &
\label{eq:darcy_velocity_1D}\\[2mm]
&T_{xx} = H_\mathrm{A}\phi_\s\Frac{\partial u}{\partial x} - p
- H_\mathrm{A}\phi_\n g_n
- H_\mathrm{B}\sum_{\eta=\v,\q,\ECM}
\phi_{\eta}g_\eta & \label{eq:stress_Txx_1D}
\end{align}
% \end{subequations}
% \bibliographystyle{plain}      % mathematics and physical sciences
% \bibliography{PAPERMB}   % name your BibTeX data base
% 
% \end{document}
% 
%
where
\begin{align}
& K(\phi_\fl)=K_\text{ref}\Psi(\phi_\fl)
&\label{eq:darcy_permeability_1D}
\end{align}
% \begin{align}
% & \mathbf{T}(u, p,\widetilde{\boldsymbol{\phi}}) =
%  \Sigma_\e(u,p,\widetilde{\boldsymbol{\phi}})-
%  \Sigma_\g(\widetilde{\boldsymbol{\phi}})
% &\label{total_stress_1D}\\[2mm]
% &\Sigma_\e(u,p,\widetilde{\boldsymbol{\phi}})=\text{diag}
% \left(H_\mathrm{A}\phi_\s \Frac{\partial u}{\partial x}
%   - p,\,\phi_\s \lambda\Frac{\partial u}{\partial x}
%   - p,\,\phi_\s \lambda\Frac{\partial u}{\partial x}
%   - p\right)&\label{eq:stress_solid_1D} \\[2mm]
%  &\Sigma_\g(\widetilde{\boldsymbol{\phi}})=\text{diag}
%  \left(H_\mathrm{A}\phi_\n
%  g_\n+H_\mathrm{B}
%  \sum_{\eta=\v,\q,\ECM}\phi_{\eta}g_\eta,\right.&\notag\\
% &\quad\left.\phi_\n\lambda\g_n+H_\mathrm{B}
% \sum_{\eta=\v,\q,\ECM}\phi_{\eta}g_\eta,
% \quad\phi_\n\lambda g_n+H_\mathrm{B}
% \sum_{\eta=\v,\q,\ECM}\phi_{\eta}g_\eta\right)
% &\label{eq:stress_growth_1D} \\[2mm]
% & v_{\s}(u)
% = \frac{\partial u}{\partial t} \qquad w=\frac{V}{\phi_\fl} &\label{eq:velocity_1D}
% %&r=\phi_\s\dfrac{\partial u}{\partial x}-\phi_\n g_\n
% \end{align}
is tissue permeability with $K_\text{ref}$ defined in Eq.~(10g) of Part I
and
\begin{align}
& \Psi(\phi_\fl)= \Frac{\phi_\fl^2}{1-\phi_\fl},
&\label{eq:psi_function}
\end{align}
while
$H_\mathrm{A}=\lambda+2\mu$ is the so-called aggregate
modulus~\cite{ateshian} and $H_\mathrm{B}=3\lambda+2\mu$.
To close the \rs{problem,} we specify the following initial and boundary
conditions:
\begin{align}
&u(x,0)=u^{0}(x) &\qquad\text{in}\quad\Omega
&\label{eq:displacement_ic1_1D}\\[2mm]
&u(0,t)=0        &\qquad\forall t \in (0,T)&\label{eq:displacement_bc1_1D}\\[2mm]
&T_{xx}(L,t) \cdot n = T_\text{b}(t) & \qquad\forall t \in (0,T)& \label{eq:displacement_bc2_1D}\\[2mm]
&p(0,t)=0 & \qquad\forall t \in (0,T)&\label{eq:Darcy_bc1_1D}\\[2mm]
&V(L,t) \cdot n = V_\text{b}(t) &\qquad\forall t \in (0,T).&\label{eq:Darcy_bc2_1D}
\end{align}
\end{subequations}

For all times $t$, at $x=0$ we set an homogeneous Dirichlet boundary condition for the solid phase. This expresses the fact the solid phase is constrained
at the scaffold wall. At $x=0$ we also enforce an homogeneous Neumann boundary condition
for the fluid phase expressing hydraulic impermeability of the scaffold wall.
At the interface with the interstitial fluid, $x=L$, we set
a Neumann boundary condition for the solid phase. This expresses the fact
that the fluid exerts a stress on the biomass. At $x=L$ we also enforce
a Neumann boundary condition for the fluid phase by equating the
Darcy flux to a given velocity of the external fluid.

\subsubsection{Mass balance \rs{IV-BVP} for nutrient concentration}
For given
$u$, $V$ and $\phi_\eta$, $\eta = \n, \v, \q, \fl$,
find the oxygen nutrient concentration
$c: \mathcal{Q}_{T_{end}} \rightarrow \mathbb{R}^+$
that satisfies the following system of partial differential equations
in balance form:
\begin{subequations}\label{oxygen_mass_balance_1D}
\begin{align}
& \frac{\partial c}{\partial t}+
\frac{\partial J_\mathrm{c}}{\partial x}= Q_\mathrm{c}(\widetilde{\boldsymbol{\phi}},
c) & \label{eq:continuity_oxygen_1D} \\[2mm]
& J_c = v_{\fl} c-D_\mathrm{c} \frac{\partial c}{\partial x}  &
\label{eq:flux_oxygen_1D}
\end{align}
where:
\begin{align}
& w=\frac{V}{\phi_\fl} & \\
& v_{\fl} = w + \dfrac{\partial u}{\partial t}
& \label{eq:v_fl_1D}
\end{align}
and
\begin{align}
& D_\mathrm{c}=D_{\mathrm{c},\fl}\dfrac{3k-2\phi_\fl(k-1)}{3+\phi_\fl(k-1)},
\qquad k:=K_\mathrm{eq}\dfrac{D_{\mathrm{c},\s}}{D_{\mathrm{c},\fl}}&
\label{eq:nutrient_diffusion_1D}\\[2mm]
&Q_\mathrm{c}(\widetilde{\boldsymbol{\phi}},c) =
- (R_\n \phi_\n + R_\v \phi_\v + R_\q \phi_\q)
\Frac{c}{c+K_{1/2}}.& \label{eq:MM_1D}
\end{align}
To close the \rs{problem,} 
we specify the following initial and boundary conditions:
\begin{align}
& c(x,0) = c^0(x) &\qquad\text{in } \Omega& \label{eq:nutrient_ic_1D}\\[2mm]
&\left.\dfrac{\partial c}{\partial x}\right|_{x=0}=0 &\qquad\forall\,t& \label{eq:nutrient_bc1_1D}\\[2mm]
&c(L,t)=c_\text{ext}(t).& \label{eq:nutrient_bc2_1D}
\end{align}
\end{subequations}
As in the case of the cellular phase model, the Neumann boundary condition at $x=0$ represents the fact that the scaffold is impermeable to oxygen flow,
while the Dirichlet boundary condition at $x=L$ indicates that we assume
the interstitial fluid to deliver \rs{to} the construct a prescribed amount
of nutrient.

\subsubsection{Mass conservation \rs{IV-BVP} for cellular populations}
For given $u$, $p$ and $c$, find the volume fractions
$\widetilde{\boldsymbol{\phi}} =
\left[\phi_\n, \, \phi_\v, \, \phi_\q, \, \phi_{\ECM} \right]^T
: (\mathcal{Q}_{T_{end}})^4 \times (\mathbb{R}^+)^4$
%and $\phi_{\fl}: \mathcal{Q}_{T_{end}} \rightarrow \mathbb{R}^+$
that satisfy the following system of partial differential equations
in balance form:
\begin{subequations}\label{mixture_mass_balance_1D}
\begin{align}
& \frac{\partial \widetilde{\boldsymbol{\phi}}}{\partial t}+
\Frac{\partial J_{\widetilde{\boldsymbol{\phi}}}}{\partial x}=
\mathbf{Q}(\widetilde{\boldsymbol{\phi}},c,
\mathbf{T})= (\mathbf{P}(\widetilde{\boldsymbol{\phi}},c,
\mathbf{T})-\mathbf{C}(c,\mathbf{T}))\widetilde{\boldsymbol{\phi}}
&\label{eq:continuity_1D} \\[2mm]
& \left(J_{\widetilde{\boldsymbol{\phi}}}\right)_\eta
= \phi_\eta v_{\s}-D_\eta\Frac{\partial \phi_\eta}{\partial x}
& \eta = \cells, \ECM \label{eq:flux_1D}
\end{align}
\rs{where $v_{\s} = \partial u/\partial t$
is the solid phase velocity,} the production terms are:
\begin{align}
& \mathbf{P}(\boldsymbol{\phi},c,\mathbf{T}) =
\Big[\phi_\fl\dfrac{c}{K_\text{sat}+c} k_\text{g},\quad0,\quad
\beta_{\q\n}H_r, \quad 0; \notag\\
&\qquad \qquad \qquad
\qquad\;0,\qquad\;0,\quad\beta_{\q \v}\left(1-H_r\right),
\quad 0;\notag\\
&\qquad\qquad\qquad\qquad\;\dfrac{1}{\tau_\text{m}},\quad\beta_{\v\q}H_r,
\quad0,\quad0;\notag\\
&\qquad\qquad\;0,\quad\dfrac{1}{V_{\cell}} \, c \,E\, k_\text{GAG} \max\Big[0,1-\dfrac{\phi_\text{ECM}}{\phi_\text{ECM,max}}\Big],
\quad0,\quad0\Big]&
\label{eq:production_matrix_1D}\\[2mm]
&\mathbf{C}(c,\mathbf{T}) =\mathrm{diag}\Big(\dfrac{1}{\tau_\text{m}}+k_\qui
(1-H_c),\; \beta_{\v \q}H_r+k_\qui(1-H_c)+ k_\text{apo},\notag\\
&\qquad \quad \beta_{\q\n}H_r+
\beta_{\q\v}\left(1-H_r \right)+k_\qui\left(1-H_c\right)+ k_\text{apo},\;
k_\text{deg}\Big).&\label{eq:consumption_matrix_1D}
\end{align}
and the fluid fraction is computed as
\begin{align}
&\phi_{\fl} = 1 - \sum_{\eta=\n, \v, \q, \ECM} \phi_{\eta}.
&\label{eq:phi_fl_evaluation_1D}
\end{align}
To close the \rs{problem,} we specify the following initial and boundary
conditions:
\begin{align}
&\widetilde{\boldsymbol{\phi}}(x,0)=\widetilde{\boldsymbol{\phi}}^0(x)
& \qquad\text{in } \Omega&\label{eq:phi_ic_1D}\\[2mm]
&\left.\dfrac{\partial \phi_\eta}{\partial x}\right|_{x=0}=0 & \quad\qquad\forall\,t\qquad\eta=\n,\v,\q,\ECM&\label{eq:phi_bc1_1D}\\[2mm]
&\left.\dfrac{\partial \phi_\eta}{\partial x}\right|_{x=L}=0&\qquad\forall\,t\qquad\eta=\n,\v,\q,\ECM.&
\label{eq:phi_bc2_1D}
\end{align}
\end{subequations}

The boundary conditions~\eqref{eq:phi_bc1_1D}-~\eqref{eq:phi_bc2_1D}
express the fact that cellular phases can flow out of the biomass
only because of the presence of an advective field.

\section{Indicator of the isotropy of the local stress state}
\label{sec:definition_of_r}
This section is devoted to the definition of
the indicator $r$ of the isotropy of the local stress state.
To this purpose, we observe that in the 1D configuration of
Fig.~\ref{fig:1D_domain}, the total stress
tensor $\mathbf{T}$ can be decomposed into the
sum of isotropic and anisotropic components as
\begin{subequations}
%\label{eq:mohr}
\begin{align}
& \mathbf{T} = \mathbf{T}_\textrm{iso} + \mathbf{T}_\textrm{aniso}
\label{eq:stress_total}
\end{align}
with:
\begin{align}
&\mathbf{T}_\textrm{iso}=\Big[\lambda\left(\phi_\s\dfrac{\partial u}{\partial x}-g_\n\phi_\n\right)-\left(\dfrac{2}{3}\mu+\lambda\right)
\sum_{\eta=\v,\q,\ECM}g_{\eta}\phi_\eta-p\Big]\mathbf{I} & \label{eq:Tiso}\\
&\mathbf{T}_\textrm{aniso}=2\mu\left(\phi_\s\dfrac{\partial u}{\partial x}-g_\n\phi_\n\right)\mathbf{d}_\textrm{pol}\otimes\mathbf{d}_\textrm{pol},
& \label{eq:Taniso}
\end{align}
and where $\mathbf{d}_\textrm{pol}$ is the unit vector $[1,0,0]^T$.
The decomposition~\eqref{eq:stress_total} suggests that
$\mathbf{T}_\textrm{aniso}$ can be used to
measure the degree of anisotropicity
of the stress state at any point $x$ of the mixture and at any time $t$.
With this aim, we define
the parameter $r$ introduced in Section 5 of Part I as
\beq\label{eq:r}
r(x,t)=\Frac{\|\mathbf{T}_\textrm{aniso}(x,t)\|_F}{2 \mu}
= \Big|\phi_\s(x,t) \Frac{\partial u(x,t)}{\partial x}-g_\n(x,t)
\phi_\n(x,t)\Big|
\eeq
where $\| \mathbf{A} \|_F$ is the Frobenius norm of a matrix
$\mathbf{A} \in \mathbb{R}^{3 \times 3}$. We can give a mechanical
interpretation of~\eqref{eq:r} by studying the Mohr circle at
point $(x,t)$. The principal components of $\mathbf{T}$ are:
\begin{align}
&\sigma_\textrm{I}=H_\textrm{A}\phi_\s\dfrac{\partial u}{\partial x}-p-H_\textrm{A}g_\n\phi_\n-H_\textrm{B}
\sum_{\eta=\v,\q,\ECM}g_{\eta}\phi_\eta\\
&\sigma_\textrm{II}=\sigma_\textrm{III}=\lambda\phi_\s\dfrac{\partial u}{\partial x}-p-\lambda g_\n\phi_\n-H_\textrm{B}\sum_{\eta=\v,\q,\ECM}g_{\eta}\phi_\eta,
\end{align}
from which it follows that the Mohr circle at $(x,t)$
has center $C=\dfrac{\sigma_\textrm{I}+\sigma_\textrm{II}}{2}$ and radius
equal to the maximum total shear stress at $(x,t)$
$$
\tau_\textrm{max}(x,t) =\dfrac{\sigma_\textrm{I}(x,t) -
\sigma_\textrm{II}(x,t)}{2}
= \mu \left(\phi_\s(x,t) \Frac{\partial u(x,t)}{\partial x}-g_\n(x,t)
\phi_\n(x,t) \right).
$$
Comparing this latter relation with~\eqref{eq:r} we conclude that
the indicator of the local stress state anisotropy can be written as
\beq\label{eq:r_2}
r(x,t) = \Frac{|\tau_\textrm{max}(x,t)|}{\mu}.
\eeq
According to Eq.~(7a) of Part I, we need also
characterize an appropriate value for the threshold parameter
$\bar{r}$ representing the level of hydrodynamic shear stress that induces metabolic activity of the cell population
$\n$ and therefore separates the isotropic regime from the anisotropic regime. In~\cite{Raimondi2008,Raimondi2006a} it is shown that hydrodynamic
shear below 10 mPa may promote GAG synthesis, so that,
coherently with~\eqref{eq:r_2}, we assume
\beq\label{eq:rbar}
\bar{r}=\dfrac{\unit{10mPa}}{\mu}.
\eeq
\end{subequations}

\section{Numerical approximation of the 1D mechanobiological model} \label{sec:numerical_approximation_1D}
In this section we focus on the numerical approximation of the
model illustrated in Sect.~\ref{sec:mechanobio_1D}. With this purpose,
we illustrate in Sect.~\ref{sec:algorithm} the computational algorithm
to iteratively solve the coupled systems of equations and in
Sect.~\ref{sec:fem1d} we shortly discuss the finite element
discretization scheme used to numerically solve the linearized equations.
%
%
% in Sect.~\ref{sec:computational_methods} we give some details of the
% computational algorithms used to translate the solution of
% the model system into the successive solution of systems
% of linearized algebraic equations.

\subsection{\rs{Computational algorithm}}\label{sec:algorithm}
Solving in closed form the mechanobiological model illustrated in
Sect.~\ref{sec:1D_model} is a very difficult task because of
the strong nonlinear nature of the problem. Therefore, a numerical
treatment is in order. Prior to discretization we need 
to reduce the solution of the whole coupled
system to the solution of a sequence of linearized 
equations of simpler form. For this purpose we set
\begin{subequations}
\begin{align}
& \mathbf{U} = \left[ u, p, \phi_\n, \phi_\v, \phi_\q, \phi_\ECM, 
c \right]^T & \label{eq:solution_U} 
\end{align}
and subdivide the time interval $[0, T_{end}]$ into
$N_T \geq 1$ uniform subintervals of length $\Delta t = T_{end}/N_T$,
in such a way that the discrete time levels $t^n= n \Delta t$, 
$n= 0, \ldots, N_T$, are obtained. 
Then, for each $n = 0, \ldots, N_T-1$, we set 
% 
% we perform the temporal semidiscretization at $t=t^{n+1}$ of the three continuity equations~\eqref{eq:continuity_2_1D},~\eqref{eq:continuity_1D} and~\eqref{eq:continuity_oxygen_1D} using the Backward Euler difference scheme (see~\cite{QSS2007}, Chapt. 11), and set
$$
\mathbf{U}^{(0)}:=\mathbf{U}^n,
$$
and for all $m \geq 0$ until convergence we perform the following fixed point iteration:
\begin{enumerate}
\item solve the linear poroelastic system:
\begin{align}
& \Frac{\partial T_{xx}^{(m+1)}}{\partial x} = 0 &
\label{eq:momentum_balance_1D_iter}\\[2mm]
%(u,p,\widetilde{\boldsymbol{\phi}})}
& \Frac{1}{\Delta t}\Frac{\partial u^{(m+1)}}{\partial x} +
\Frac{\partial V^{(m+1)}}{\partial x} = \Frac{1}{\Delta t}
\Frac{\partial u^{n}}{\partial x}
& \label{eq:continuity_2_1D_iter}\\
& V^{(m+1)}=-K(\phi_\fl^{(m)})\frac{\partial p^{(m+1)}}{\partial x} &
\label{eq:darcy_velocity_1D_iter}\\[2mm]
&T_{xx}^{(m+1)} = H_\mathrm{A}\phi_\s^{(m)}
\Frac{\partial u^{(m+1)}}{\partial x} - p^{(m+1)}
- H_\mathrm{A}\phi_\n^{(m)} g_n^{(m)}
- H_\mathrm{B}\sum_{\eta=\v,\q,\ECM}
\phi_{\eta}^{(m)}g_\eta^{(m)} & \label{eq:stress_Txx_1D_iter}
\end{align}
supplied by the boundary 
conditions~\eqref{eq:displacement_bc1_1D}-~\eqref{eq:Darcy_bc2_1D}.
This yields the updated solid displacement $u^{(m+1)}$ and fluid
pressure $p^{(m+1)}$;
\item solve the linear advection-diffusion-reaction (ADR) system:
\begin{align}
& \frac{c^{(m+1)}}{\Delta  t} +
\frac{\partial J_\mathrm{c}^{(m+1)}}{\partial x} - 
\widehat{Q}_\mathrm{c}^{(m)}(\widetilde{\boldsymbol{\phi}}^{(m)}, 
c^{(m)}) c^{(m+1)}
= \frac{c^n}{\Delta  t} & \label{eq:continuity_oxygen_1D_iter} \\[2mm]
& J_c^{(m+1)} = v_{\fl}^{(m+1)} c^{(m+1)} - 
D_\mathrm{c}^{(m)} \frac{\partial c^{(m+1)}}{\partial x}  &
\label{eq:flux_oxygen_1D_iter}
\end{align}
where:
\begin{align}
& w^{(m+1)} = \frac{V^{(m+1)}}{\phi_\fl^{(m)}} & \\
& v_{\fl}^{(m+1)} = w^{(m+1)} + \dfrac{u^{(m+1)}-u^n}{\Delta t}
& \label{eq:v_fl_1D_iter}
\end{align}
and
\begin{align}
& D_\mathrm{c}^{(m)} =D_{\mathrm{c},\fl}\dfrac{3k-2\phi_\fl^{(m)}(k-1)}{3+\phi_\fl^{(m)}(k-1)},
\qquad k:=K_\mathrm{eq}\dfrac{D_{\mathrm{c},\s}}{D_{\mathrm{c},\fl}}&
\label{eq:nutrient_diffusion_1D_iter}\\[2mm]
&\widehat{Q}_\mathrm{c}^{(m)}(\widetilde{\boldsymbol{\phi}}^{(m)},
c^{(m)}) = - (R_\n \phi_\n^{(m)} + R_\v \phi_\v^{(m)} + R_\q \phi_\q^{(m)})
\Frac{1}{c^{(m)}+K_{1/2}}& \label{eq:MM_1D_iter}
\end{align}
supplied by the boundary 
conditions~\eqref{eq:nutrient_bc1_1D}-~\eqref{eq:nutrient_bc2_1D}. 
This yields the updated oxygen concentration $c^{(m+1)}$;
\item solve the linear advection-diffusion-reaction system:
\begin{align}
& \frac{\widetilde{\boldsymbol{\phi}}^{(m+1)}}{\Delta t} +
\Frac{\partial J_{\widetilde{\boldsymbol{\phi}}}^{(m+1)}}{\partial x} 
+ \mathbf{C}(c^{(m+1)},\mathbf{T}^{(m+1)}))
\widetilde{\boldsymbol{\phi}}^{(m+1)} = & \nonumber \\
& \mathbf{P}(\widetilde{\boldsymbol{\phi}}^{(m)},c^{(m+1)},
\mathbf{T}^{(m+1)}) \widetilde{\boldsymbol{\phi}}^{(m)} 
&\label{eq:continuity_1D_iter} \\[2mm]
& \left(J_{\widetilde{\boldsymbol{\phi}}}\right)_\eta^{(m+1)}
= \phi_\eta^{(m+1)}  v_{\s}^{(m+1)} - 
D_\eta \Frac{\partial \phi_\eta^{(m+1)}}{\partial x}
\qquad \eta = \cells, \ECM & \label{eq:flux_1D_iter}
\end{align}
%linearized mass balance system~\eqref{mixture_mass_balance_1D}.
supplied by the boundary 
conditions~\eqref{eq:phi_bc1_1D}-~\eqref{eq:phi_bc2_1D}. 
This yields the updated cellular volume fractions
$\boldsymbol{\phi}^{(m+1)}$.
\end{enumerate}
Should the above fixed point iteration 1.-3. reach convergence
at a certain value $m^\ast \geq 0$, then we set
$$
\mathbf{U}^{n+1}:= \mathbf{U}^{(m^\ast)}
$$
and we increment the time loop counter by setting
$$
n \leftarrow n + 1.
$$
until conclusion of the time advancement loop.
\end{subequations}

Two remarks are in order about the above described solution map.
The first remark concerns with the linear poroelastic 
system~\eqref{eq:momentum_balance_1D_iter}-~\eqref{eq:stress_Txx_1D_iter}.
The weak formulation of this problem leads to solving a saddle-point
problem in block symmetric form to which the abstract analysis 
of~\cite{QuarteroniValli94}, Chapt.~7 
and~\cite{BrezziFortin1991} can be applied to prove existence and uniqueness
of the solution pair $u^{(m+1)}, \, p^{(m+1)}$.
The second remark concerns with the two linear ADR problems.
The splitting of the source term ${Q}_\mathrm{c}$ in~\eqref{eq:continuity_oxygen_1D} and $\mathbf{Q}$ in~\eqref{eq:continuity_1D}
gives rise to two BVPs to which the application
of the maximum principle (see~\cite{RoosStynesTobiska1996}) 
allows to prove nonnegativity of the 
solutions $c^{(m+1)}$ and $\phi_\eta^{(m+1)}$, $\eta = \cells, \ECM$.

\subsection{Finite element discretization}\label{sec:fem1d}
The computational procedure described in Sect.~\ref{sec:algorithm}
leads to solving two kinds of BVPs:
\rs{(i) a saddle-point problem; (ii) two ADR equations.}
We numerically solve (i) and (ii) using the Galerkin finite element approximation scheme on a family of partitions
$\left\{\mathcal{T}_h \right\}_{h>0}$ of
the computational domain, $h$ being the discretization parameter (see~\cite{QuarteroniValli94}).
In the case of the saddle point problem (i) we
employ piecewise linear finite elements on $\mathcal{T}_h$
for both solid
displacement and fluid pressure. Equal-order interpolation
for $u$ and $p$ does not give rise to numerical instabilities as
it would be the case if the Stokes equations for an incompressible fluid
were to be solved (cf.~\cite{QuarteroniValli94}, Chapt. 9), because
in the present model the variable $p$ is {\em not}
a Lagrange multiplier (as in the Stokes system), rather,
it is the solution of the elliptic Darcy
problem~\eqref{eq:continuity_2_1D}-\eqref{eq:darcy_velocity_1D}.
In the case of the ADR equation we employ for the approximation
of the concentration and of the cellular volume fractions
the primal-mixed finite element discretization scheme
with exponential fitting
stabilization proposed and investigated in~\cite{SaccoABB2014}.
This choice is taken because it ensures that the computed numerical
solutions satisfy a strict positivity property even in the case of
a strongly advective regime. Moreover, it can be checked that, if
advective terms do not play a major role compared to oxygen molecular
diffusion in the biomass, then the effect of the stabilization
introduced by the primal-mixed method of~\cite{SaccoABB2014} becomes
negligible so that the accuracy of the scheme is not spoiled.
This, instead,
would not be the case if the classic upwind stabilization were
adopted (see~\cite{Brooks1982} for a discussion of this important
issue).

\section{Simulation tests}\label{sec:simulations}

% The main goal is to assess the bio-physical accuracy
% of model predictions by addressing in a self-consistent coupled manner
% the principal biophysical events that ultimately drive biomass
% production under the action of a perfusion fluid.
%
% In this section we conduct a series of simulations
% of the 1D bioengineered tissue growth under different working
% conditions compatible with realistic approaches.
%
%
% The remarkable conclusion of these simulations is that, despite
% the simplicity of the mechanical configuration, the obtained results
% fit very well with the experimental outcomes illustrated in~\cite{Nikolaev,Novakovic,Obradovic}.

In this section we show the numerical results obtained
by solving the 1D problem with the computational algorithm described
in Sect.~\ref{sec:numerical_approximation_1D}.
% With this aim, we show two groups of numerical simulations.
Two sets of simulation tests are performed. In the first set
of simulations we set $T_\text{b} = V_\text{b}=0$ in Eqns.~\eqref{eq:displacement_bc2_1D}-\eqref{eq:Darcy_bc2_1D}.
This corresponds to investigating a static culture  environment (see~\cite{Raimondi2008,Raimondi2006a}). In the second set
of simulations we set $T_\text{b} = 100 \unit{mPa}$ and
$V_\text{b} = 50 \unit{\mu m s^{-1}}$ as in~\cite{Causin}.
These values are characteristic of a culture in a perfusion
bioreactor where an external hydrodynamic shear stress is applied~\cite{Sacco,Raimondi2004,Raimondi2006a}.

The principal scope of the numerical experiments is to perform
a sensitivity analysis of the model and its effect on the predicted
profiles of cellular components. This analysis allows a biophysical
validation of the proposed model and may provide
a useful guideline for an optimal design and calibration
of a bioreactor in \rs{TE} applications.

The first investigated input model parameter is the amount $A$
of cell density at the beginning of the culture process ($t=0$)
and at the pore wall ($x=0$). We adopt the
following exponentially decaying profiles for cells and ECM at $t=0$:
\begin{subequations}
\begin{align}\label{exponential_initial_condition}
\phi_\n&=A_\n\exp(-x/L_d)\\
\phi_\eta&=A_\eta\exp(-x/L_d)\qquad \eta=\v, \q, \ECM,
\end{align}
with $L_d = L/5$, and
for each set of simulations we use the following values of $A$:
\begin{center}
\begin{description}
\item[(IC1)] $A_n=0.005$, $A_\eta=0.001 \qquad\eta=\v,\q,\ECM$;
\item[(IC2)] $A_n=0.05$, $A_\eta=0.01 \qquad \quad\eta=\v,\q,\ECM$.
\end{description}
\end{center}
The above values of $A_n$ and $A_{\eta}$
agree with the biophysical
evidence that at the beginning of the growth process,
proliferating cells are present in larger amount than
the other cellular populations.
% of the initial data represent a compromise between accuracy of the
% numerical model and computation time and take account of the grater amount,
% at the beginning of the experiment, of proliferating cell with respect
% to the other populations.

The second investigated input model parameter is the cellular growth rate
$k_\text{g}$. In our computations we use two values of this parameter,
$k_\text{g1}$ and $k_\text{g2}$ (cf. Tab.~\ref{parameters}).
These two values are selected by comparison with the
maximum specific cell growth rate $k_\text{g0}$ used in~\cite{Sacco}
in such a way that $k_\text{g} < k_\text{g0}$ corresponds to
\rs{"low growth regime" whereas $k_\text{g} > k_\text{g0}$ corresponds to 
"high growth regime".} 
% . The threshold separating the two levels of growth of
% the cellular construct is represented by a critical value of the growth rate,
% denoted $k_\text{g,cr}$. This value can be estimated by a linear stability
% analysis of a reduced model represented by the sole
% system~\eqref{eq:continuity_1D} under the simplified assumption of homogenized spatial conditions of cellular populations and given oxygen concentration
% and stress distribution. This stability
% analysis, that will be reported in detail in a forthcoming publication,
% furnishes the following estimate of the critical growth rate
% \begin{align}\label{eq:kgCR}
% & k_{\text{g,cr}}  =
% \dfrac{1+k_\qui\tau_\m(1-H_c)-\tau_\m\beta_{qn}\delta_qH_r}{\tau_\m}
% \cdot \Big(\dfrac{c}{c+K_{\text{sat}}} \Big)^{-1} \\
% & \delta_q  =\dfrac{1}{\tau_\m} \cdot
% \dfrac{1}{\beta_{q-n}H_r+\beta_{q-v}
% (1-H_r)+k_{\qui}(1-H_c)+k_\apo-\beta_{v-q}\alpha_{q}H_r}.
% \end{align}
% In the numerical experiments we set $k_\text{g,1} = k_\text{g,cr}/????$
% and $k_\text{g,2} = k_\text{g,cr}*????$.

The third investigated input model parameter is the maximum
value of the external oxygen concentration $c_\text{ext}$ in Eq.~\eqref{eq:nutrient_bc1_1D} that is supplied
to the growing structure by the surrounding environment.
To determine the effect of oxygen availability on biomass growth
we set $c_\text{ext}(t) = c_\text{sat}$ and $c_\text{ext}(t) =
c_\text{thr}$ for all $t \in [0, T_{end}]$, $c_\text{sat}$ and
$c_\text{thr}$ being the saturation and threshold oxygen concentration,
respectively (see Sect.~5 of Part I and Tab.~\ref{parameters}).

% In each simulation set we also investigate the influence of the growth rate
% parameter $k_\text{g}$ on system functional behavior. As a matter of fact,
% a theoretical stability analysis of the equilibrium points of the
% mechanobiological model of Sect.~\ref{sec:mechanobio_1D}
% (not reported here and object of a forthcoming publication),
% shows that there exists a critical value $k_{\text{g,cr}}$
% of the growth rate $k_\text{g}$, defined as
% \begin{align}\label{eq:kgCR}
% & k_{\text{g,cr}}  =
% \dfrac{1+k_\qui\tau_\m(1-H_c)-\tau_\m\beta_{qn}\delta_qH_r}{\tau_\m}
% \cdot \Big(\dfrac{c}{c+K_{\text{sat}}} \Big)^{-1} \\
% & \delta_q  =\dfrac{1}{\tau_\m} \cdot
% \dfrac{1}{\beta_{q-n}H_r+\beta_{q-v}
% (1-H_r)+k_{\qui}(1-H_c)+k_\apo-\beta_{v-q}\alpha_{q}H_r},
% \end{align}
% in correspondence of which the local stress state
% switches from isotropic to anisotropic.
% This drastic transitional behavior of the model is well documented by
% the simulation results discussed later in this section and
% is responsible for a drastic change in the cell adhesion mechanisms and,
% consequently, in the functional response of the cellular growing system.

We conclude this introduction to numerical simulations by
defining the following new (equivalent) parameter $\xi = \xi(r)$ for
a synthetic representation of the isotropy indicator $r$
\beq\label{eq:csi}
\xi (r): =
\begin{cases}
1 & \quad \text{if} \, r < \bar{r}\\
0 & \quad \text{if} \, r > \bar{r}.
\end{cases}
\eeq
\end{subequations}
\rs{In the remainder of the discussion,} 
no plot is reported for the fluid volume
fraction $\phi_{\fl}$ because this variable can be computed
by post-processing using~\eqref{eq:phi_fl_evaluation_1D}.
Simulations are run over the time interval
$[0, T_{end}]$, with $T_{end}=30$ days, and the one-dimensional
plots show the time evolution of the solid and fluid mixture components
at the spatial coordinate $x=L/2$.
The values of model parameters used in the numerical
experiments are reported in Tab.~\ref{parameters}.

\subsection{Static culture}
In this section the boundary data of the poroelastic
system ~\eqref{mixture_momentum_balance_1D} are $V_b = T_b =0$ while
the external nutrient concentration $c_{\rm ext}$ is equal to
$c_{\rm sat}$ and $c_{\rm thr}$, respectively.
% This choice does not represent an approximation, since we verify that the evolution of the volumetric fractions
% of cells and ECM is in general independent from the spatial coordinate.

\subsubsection{Initial condition IC1}\label{sec:Initial_condition_IC1}

\begin{figure}[h!]
\begin{minipage}[c]{1\textwidth}
\centering
{\includegraphics[width=0.5\textwidth]{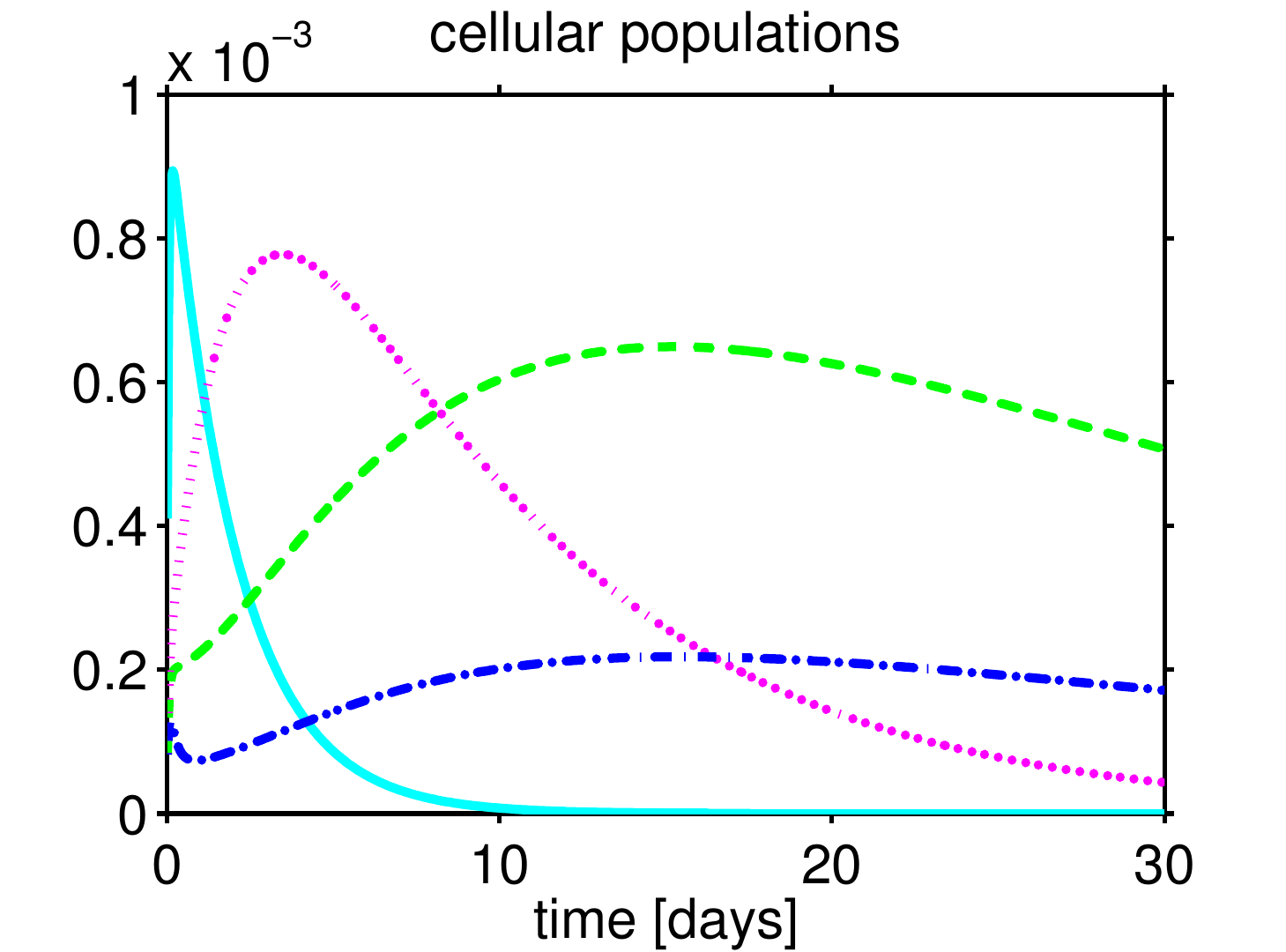}}
\\[0.3cm]
{\includegraphics[width=1\textwidth]{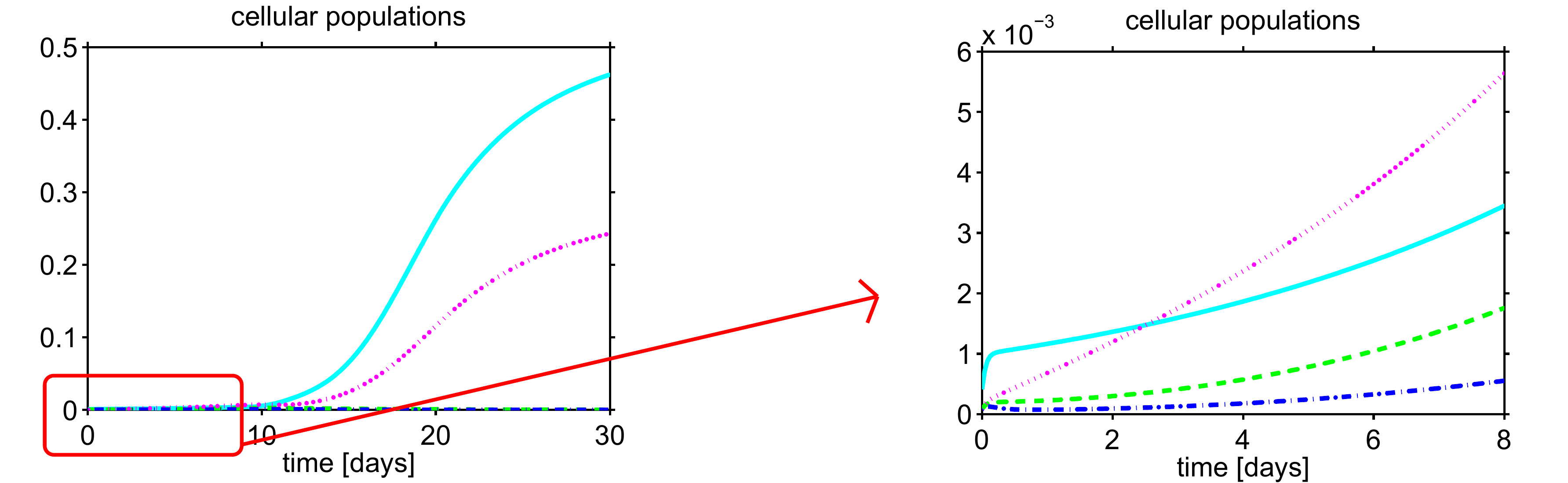}}
\end{minipage}
\caption{Temporal evolution of cellular populations
and ECM in the static culture for $c_\text{ext}=c_\text{sat}$.
Initial condition IC1. Top: $k_\text{g} = k_\text{g1}$.
Bottom left: $k_\text{g} = k_\text{g2}$. Bottom right: $k_\text{g} = k_\text{g2}$, zoom of the first eight days of culture.
\rs{Solid line: $\phi_\n$; dashed line: $\phi_\v$;
dotted line: $\phi_\q$; dash-dot line: $\phi_\ECM$}.
}
\label{fig:Static_phi}
\end{figure}

\begin{figure}[h!]
\begin{minipage}[c]{1\textwidth}
{\includegraphics[width=0.5\textwidth]{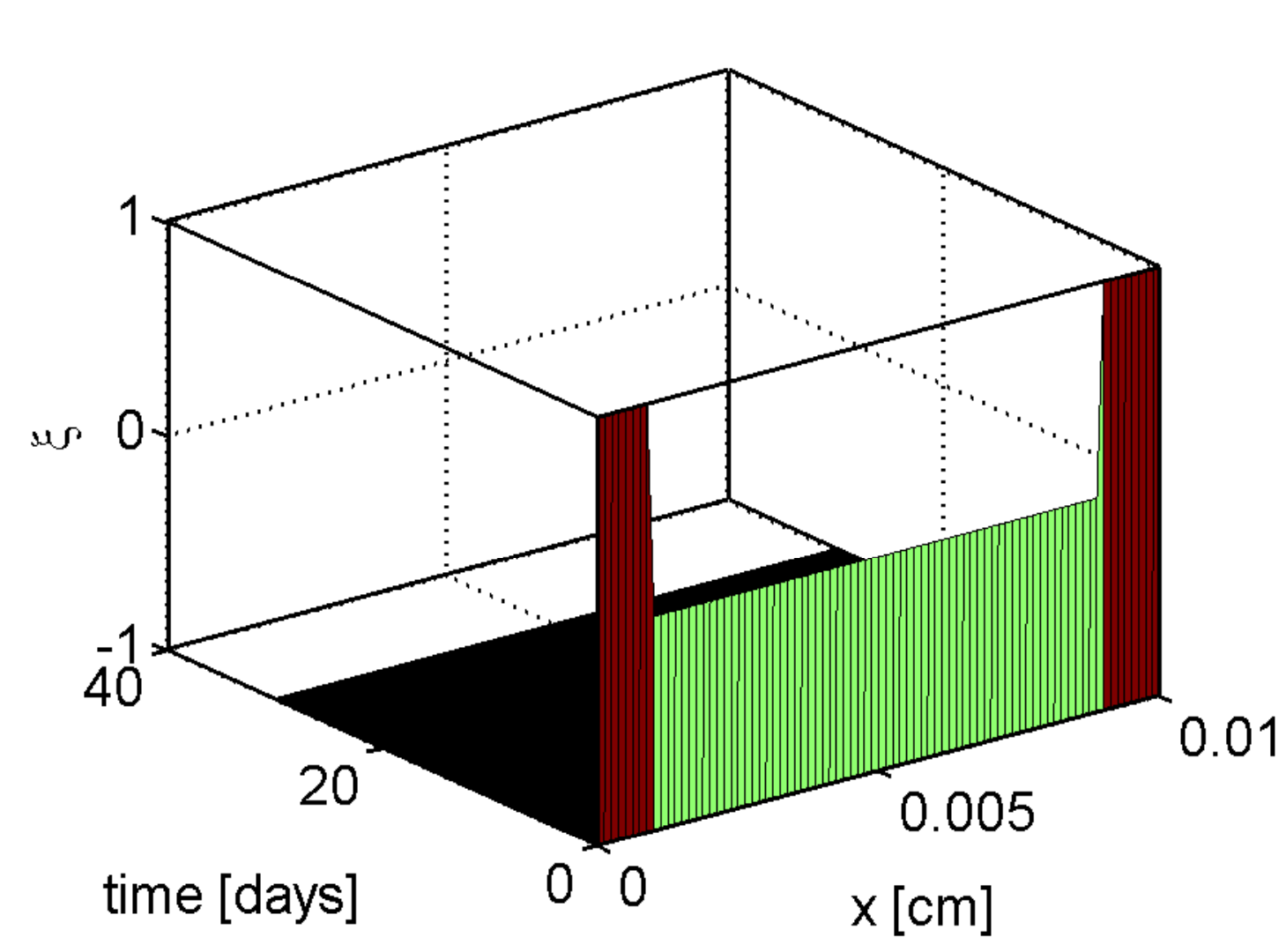}}
{\includegraphics[width=0.5\textwidth]{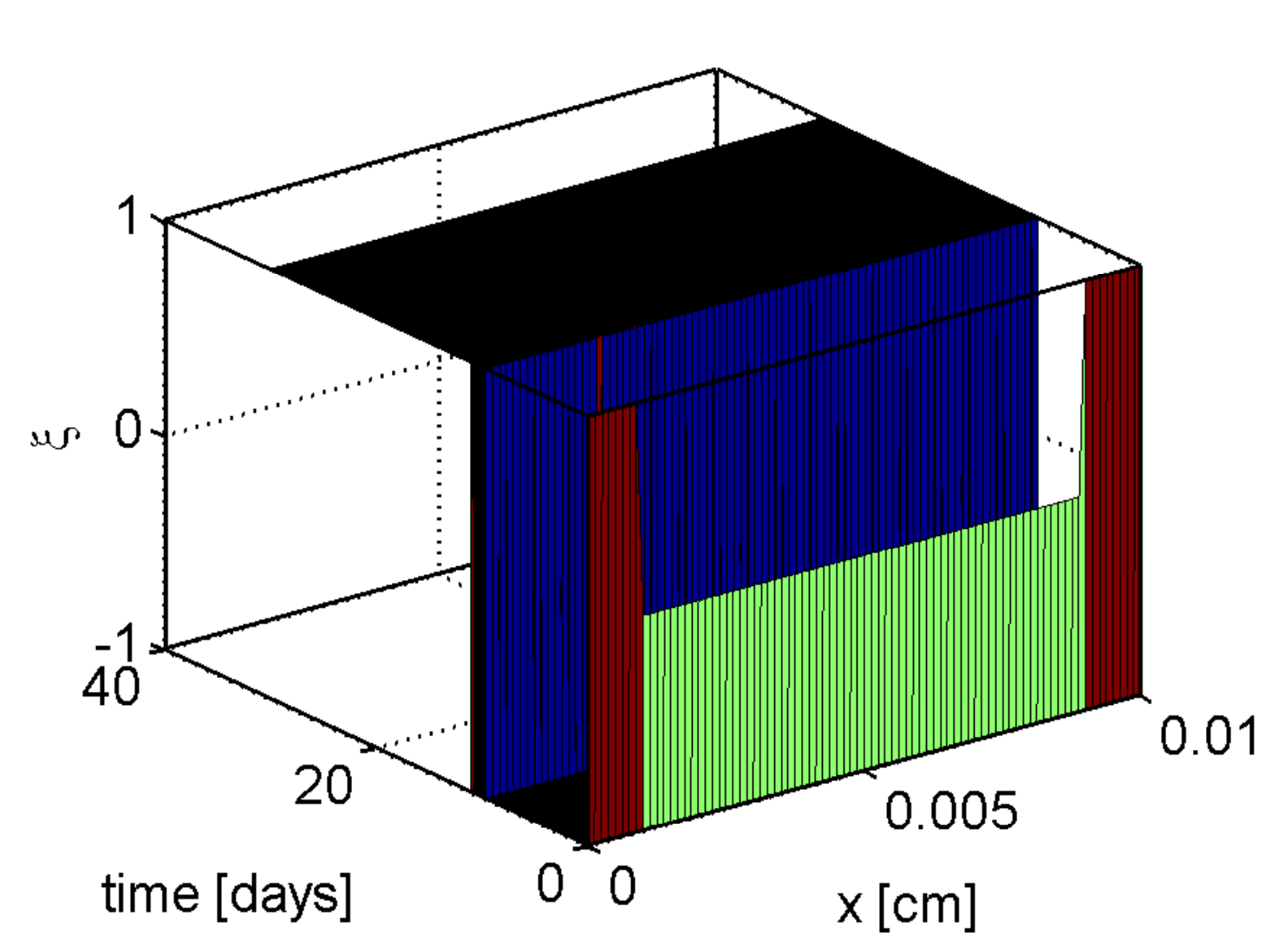}}
\end{minipage}
\caption{Spatial and temporal evolution of parameter $\xi$ in the static culture for
$c_\text{ext}=c_\text{sat}$. Initial condition IC1.
Left: $k_\text{g} = k_\text{g1}$. Right: $k_\text{g} =k_\text{g2}$.}
\label{fig:Static_csi}
\end{figure}

Figure~\ref{fig:Static_phi} illustrates the comparison between the behavior of cell populations and ECM for two different values
of the cell growth rate and in correspondence of a very
high level of nutrient concentration ($c_\text{ext}=c_\text{sat}$).
In the low growth regime ($k_\text{g} = k_\text{g,1}$), cell mitosis
experiences a drastic increase during the first days of culture, then the volumetric fraction of proliferating cells diminishes and rapidly tends to
zero (Fig.~\ref{fig:Static_phi}, top panel).
On the other hand synthesizing cells slowly increase and reach
the maximum value after about 10 days of culture.
Due to cell apoptosis and ECM degradation, all the solid volumetric fractions tend to vanish as culture time increases.
This behavior is in accordance with
the situation described by the experimental set-up,
where, during the first days of culture, cells proliferate
intensively, and then after 10-12 days, the cell-polymer construct starts to increase in size mainly due to ECM deposition~\cite{Obradovic,Nikolaev}.
The above described scenario is also consistent with the
behavior of the parameter $\xi$ that predicts an isotropic
adherence state at each time $t>0$
and at each point $x$ of the growing biomass (Fig.~\ref{fig:Static_csi},
left panel).

In the high growth regime ($k_\text{g} = k_\text{g,2}$) the mitotic profile
apparently evolves in a complimentary way with respect to the low growth case
(Fig.~\ref{fig:Static_phi}, bottom left panel).
The $\n-$cell density markedly increases after 10 days of culture,
whereas the evolution of the ECM density is negligible.
Actually, Fig.~\ref{fig:Static_phi}, bottom right panel,
shows that, in the very initial phase of the cultivation, the growth of proliferating cells is similar to that in the low growth situation.
\textcolor{black}
{However, once the corresponding maximum in the low growth regime
has been reached, mitotic activity persists because 
of the up-regulated growth rate and,} after about 30 days of culture, 
stabilizes around a \textcolor{black}{point of equilibrium.} Note that in Fig.~\ref{fig:Static_phi}, bottom left panel,
the increase of $\phi_\n$ becomes significant
in correspondence of the manifestation of the anisotropic region
in Fig.~\ref{fig:Static_csi} (right panel).

\begin{figure}[h!]
\begin{minipage}[c]{1\textwidth}
\begin{center}
{\includegraphics[width=0.45\textwidth]{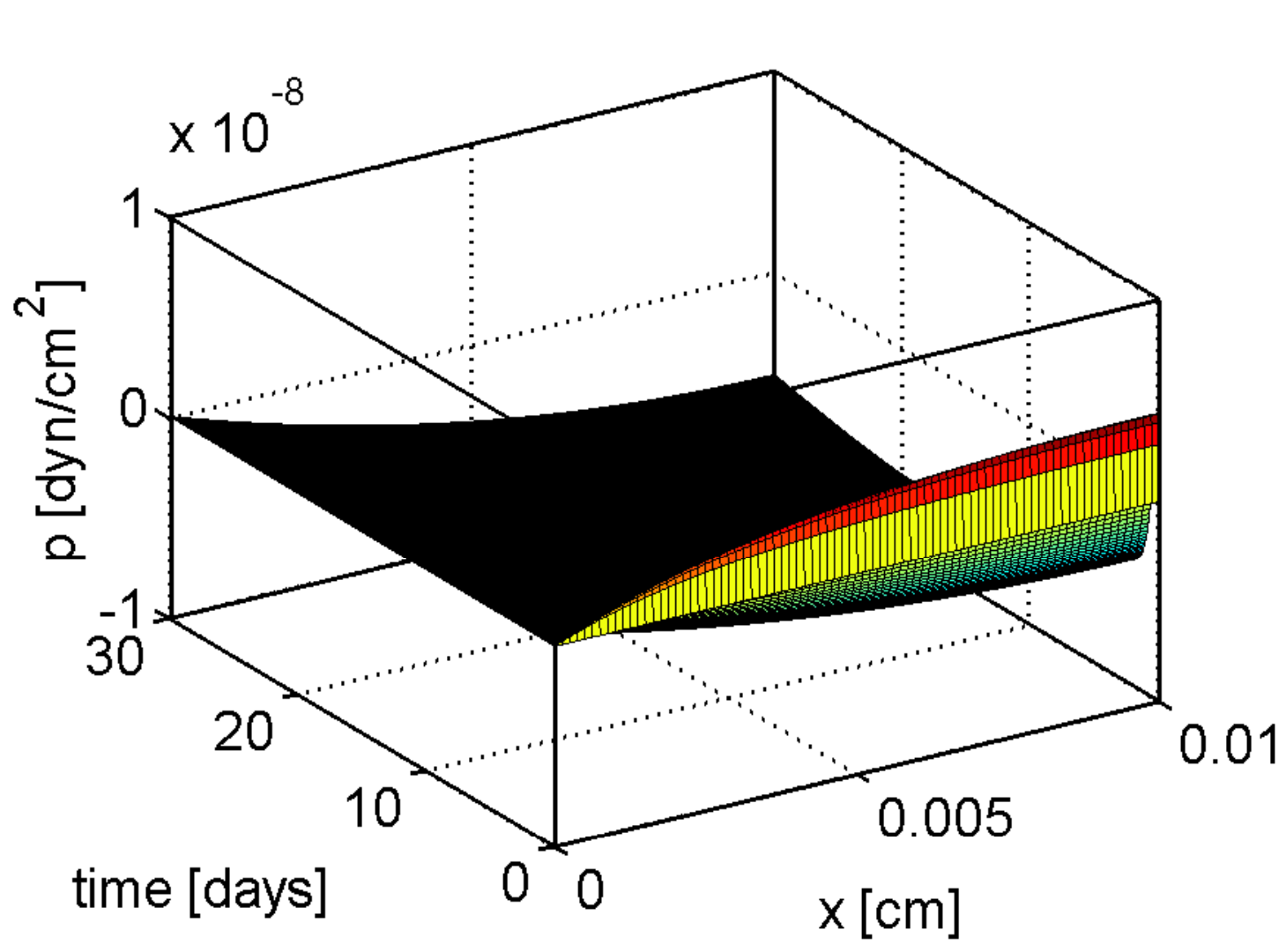}}\qquad
{\includegraphics[width=0.45\textwidth]{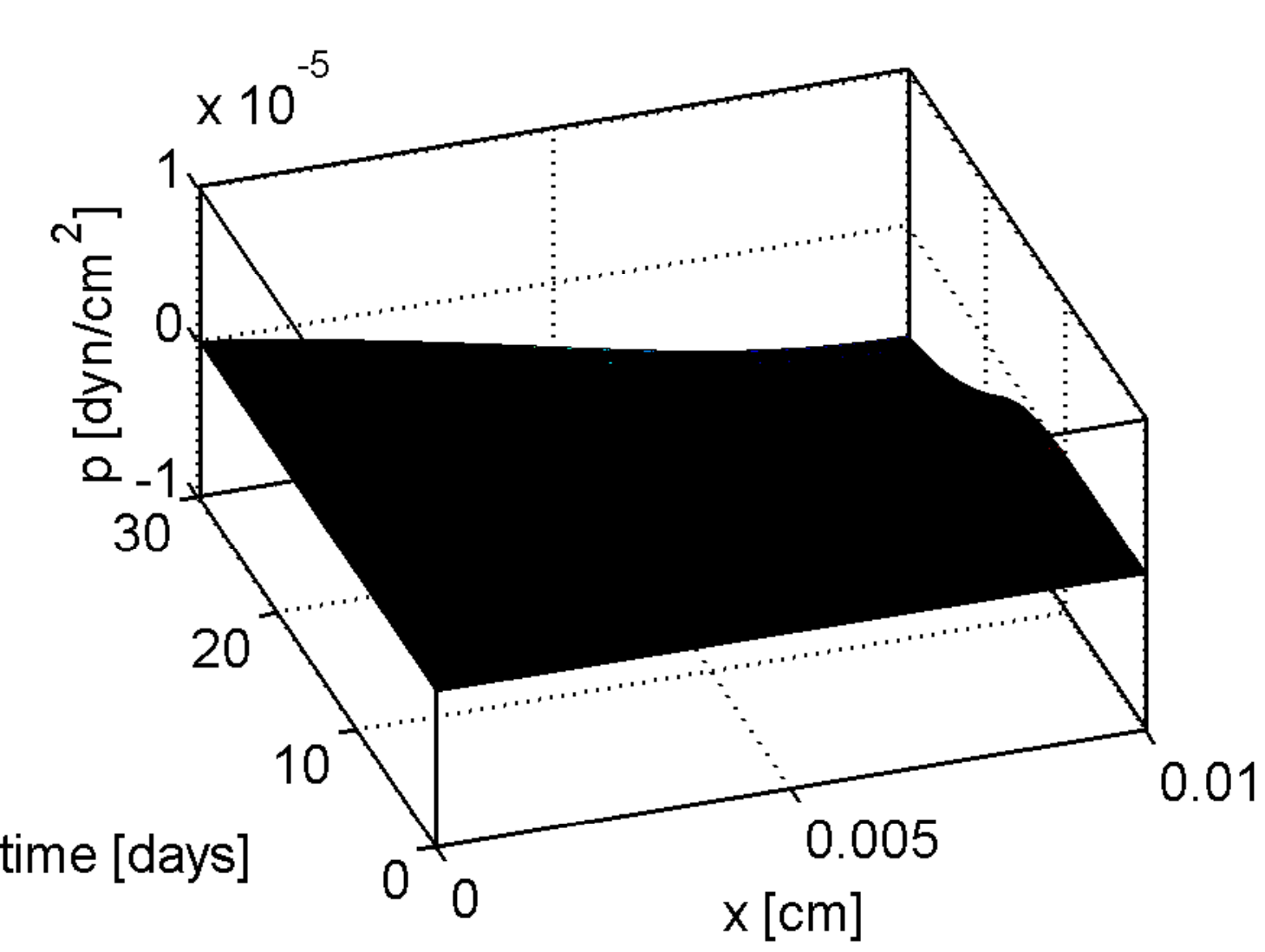}}
\end{center}
\end{minipage}
\caption{Spatial and temporal evolution of fluid pressure in
the static culture for $c_\text{ext}=c_\text{sat}$.
Initial condition IC1. Left: $k_\text{g}=k_\text{g1}$. Right: $k_\text{g}=k_\text{g2}$.}
\label{fig:Static_fl}
\end{figure}

The spatial and temporal distribution of fluid pressure
is shown in Fig.~\ref{fig:Static_fl}.
We notice that the fluid pressure drop is significantly
larger in the high growth regime (right panel) than in the
lower growth regime (left panel).
This different quantitative response of the fluid phase agrees with the
above predicted cellular populations and expresses
the biophysical fact that the larger the cellular phase, the larger
the forces exerted by this latter on the fluid. We also notice that
the spatial distribution in the low growth regime is remarkably
different from that in the high growth regime. In the former case
(left panel) the simulated working conditions practically
correspond to a linear Darcy hydraulic model with constant permeability,
subject to mixed homogeneous boundary conditions and
to a constant strain rate spatial
distribution from which it follows that the fluid pressure attains
a parabolic profile. In the latter case (right panel)
the simulated working conditions significantly deviate from a linear
Darcy hydraulic model because of the increase of cellular populations
around the time level $t = 12$ days which corresponds to a decrease
of the hydraulic permeability of the mixture in accordance
with~\eqref{eq:psi_function} and, consequently,
to an increase of the pressure drop compared to the low growth condition.

\begin{figure}[h!]
\begin{minipage}[c]{1\textwidth}
{\includegraphics[width=0.5\textwidth]{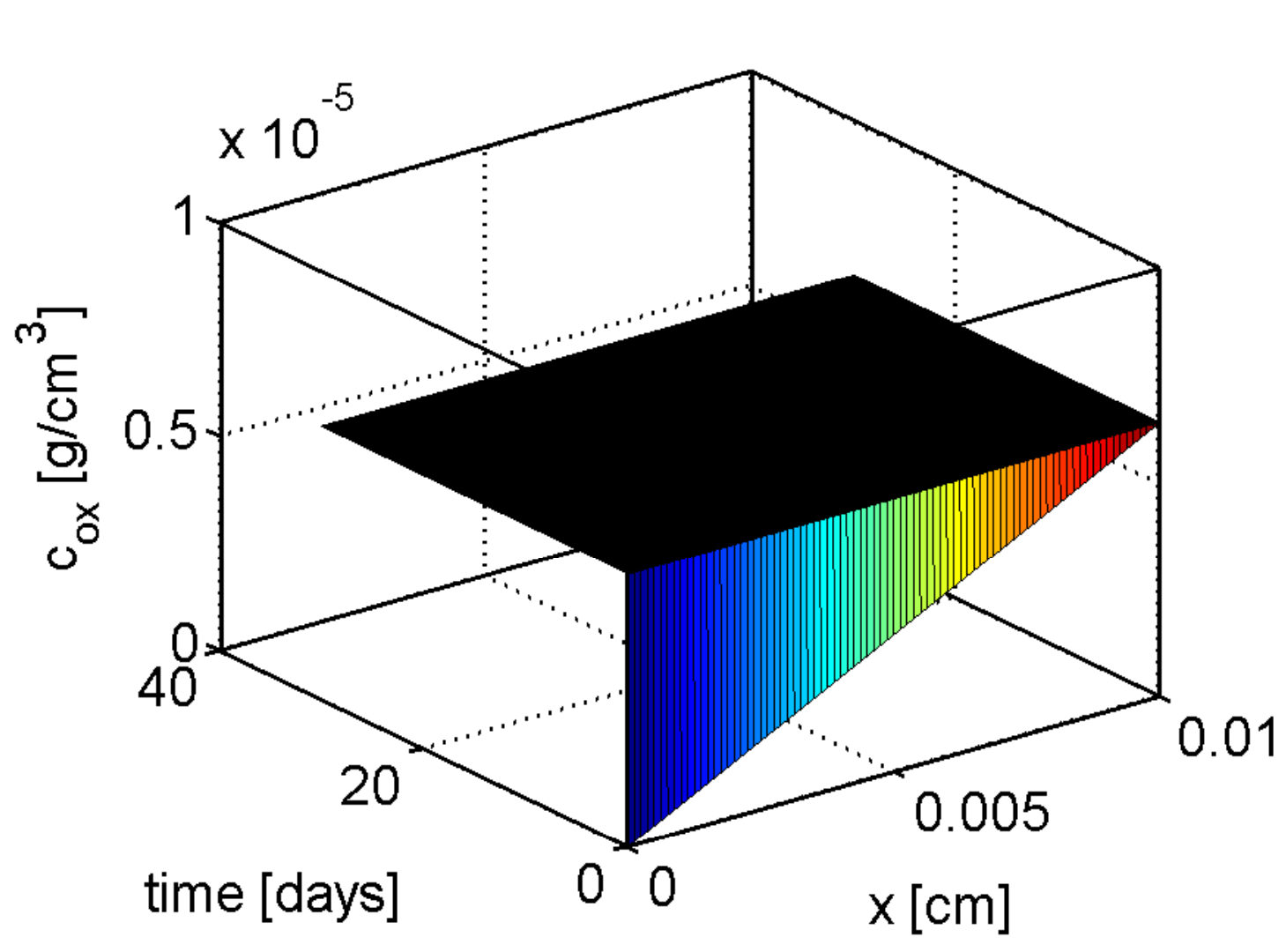}}
{\includegraphics[width=0.5\textwidth]{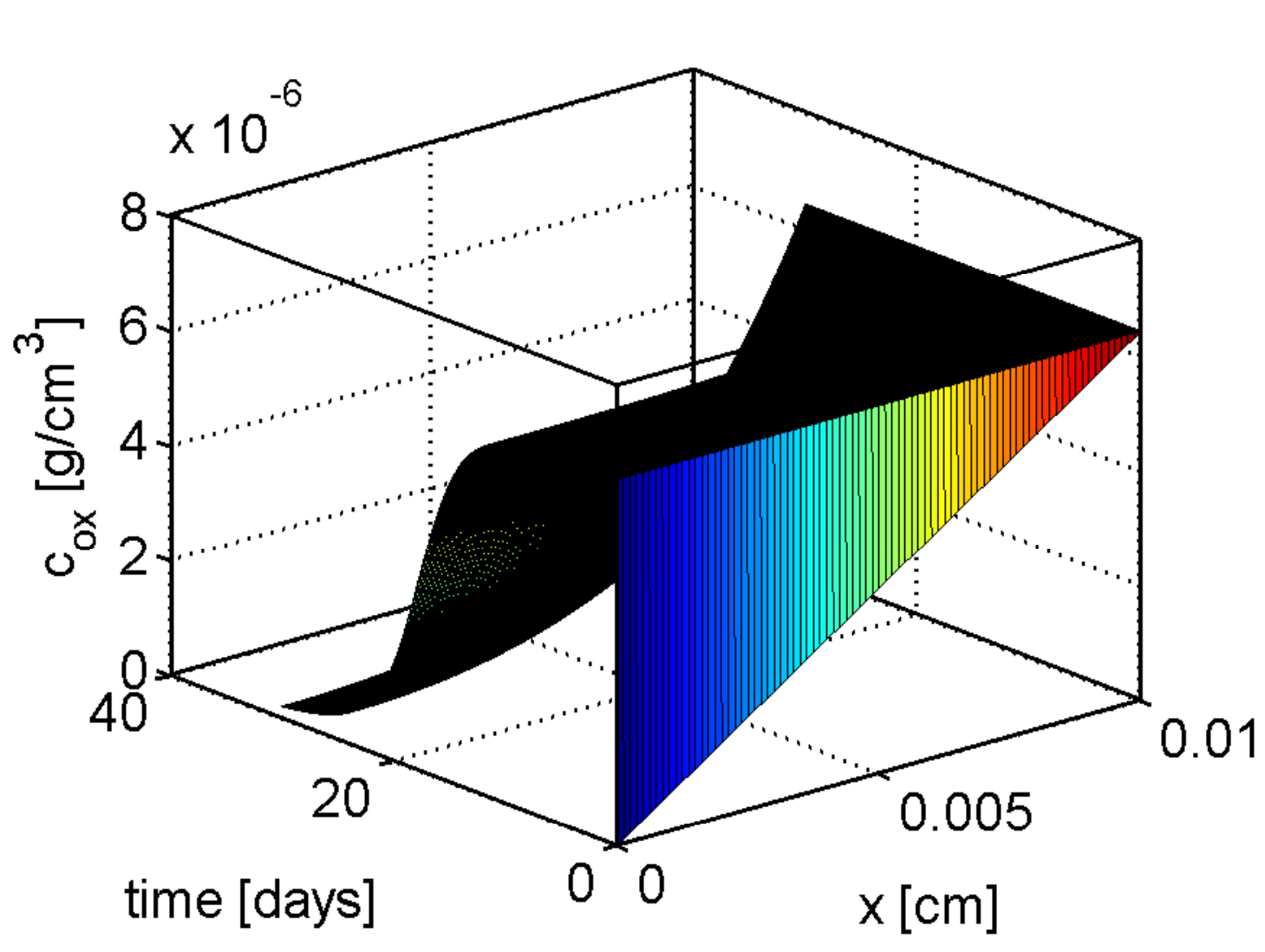}}
\end{minipage}
\caption{Spatial and temporal evolution of oxygen concentration in the static culture for $c_\text{ext}=c_\text{sat}$. Initial condition IC1. Left:
$k_\text{g} = k_\text{g1}$. Right: $k_\text{g} = k_\text{g2}$.}
\label{fig:Static_cox}
\end{figure}

\begin{figure}[h!]
\begin{minipage}[c]{1\textwidth}
{\includegraphics[width=0.5\textwidth]{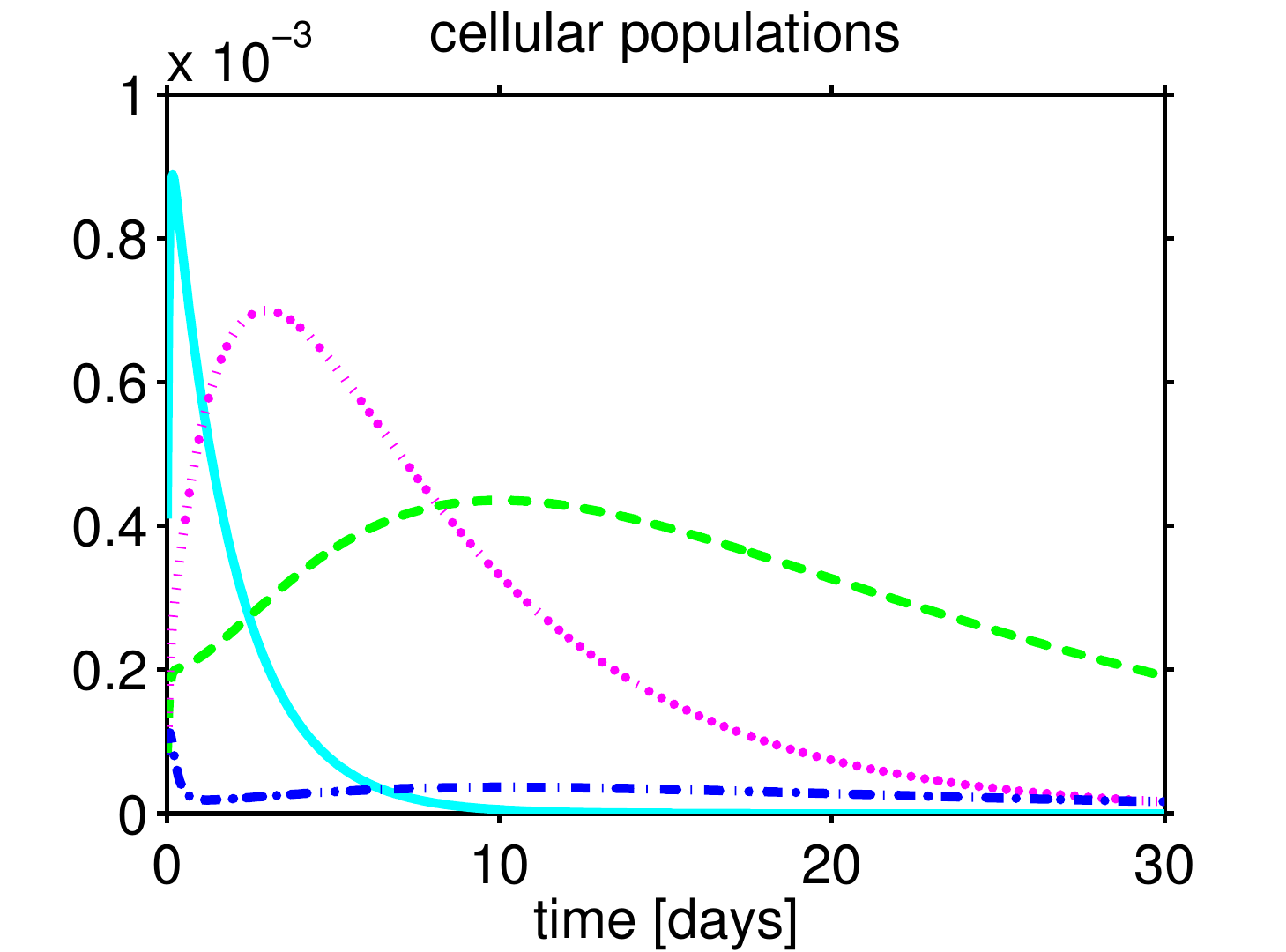}}
{\includegraphics[width=0.5\textwidth]{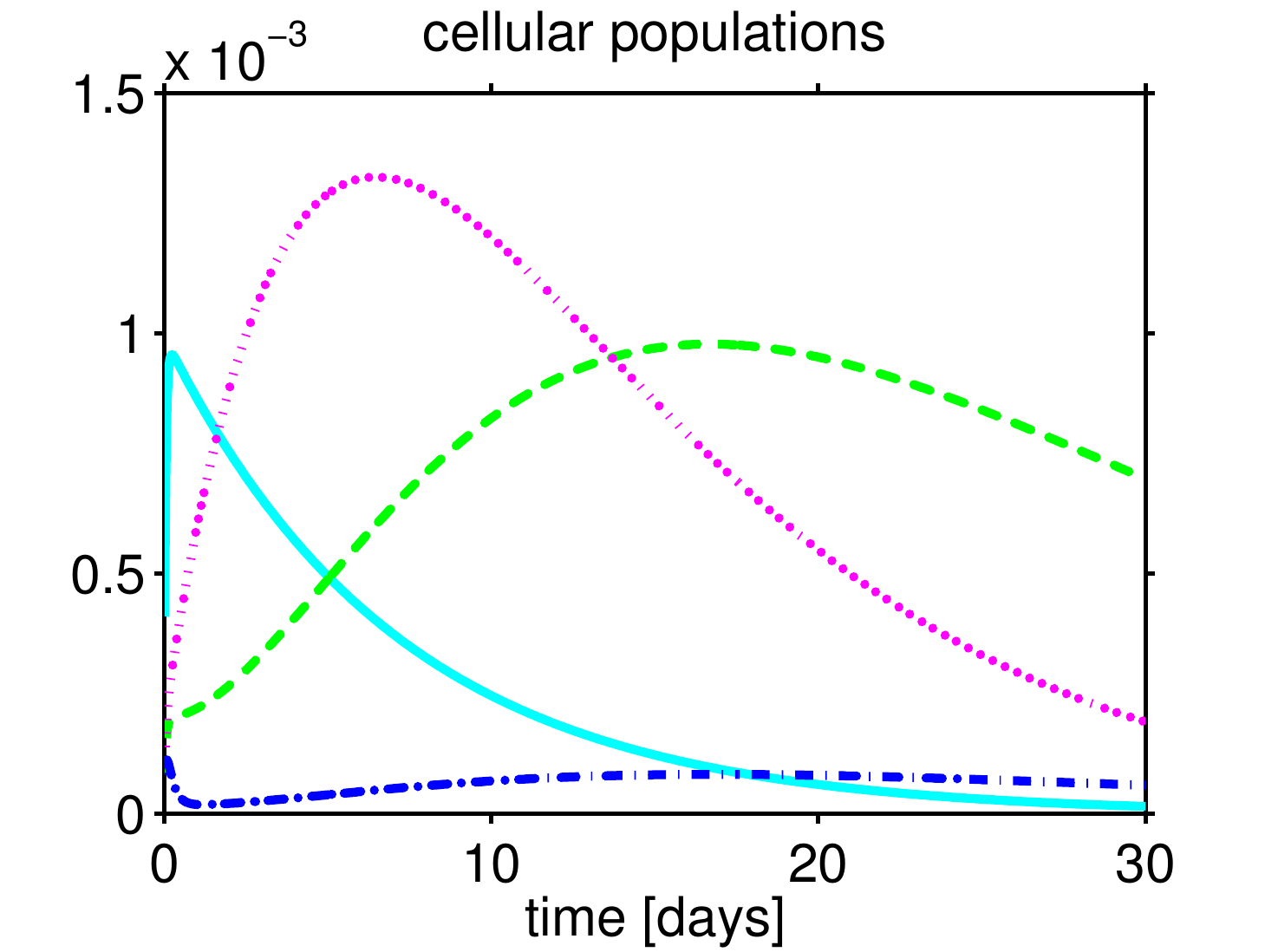}}
\end{minipage}
\caption{Temporal evolution of cellular populations and ECM in the static culture for $c_\text{ext}=c_\text{thr}$.
 Initial condition IC1. Left: $k_\text{g} = k_\text{g1}$. Right: $k_\text{g} = k_\text{g2}$. \rs{Solid line: $\phi_\n$; dashed line: $\phi_\v$;
dotted line: $\phi_\q$; dash-dot line: $\phi_\ECM$.}}
\label{fig:Static_phi_cthr}
\end{figure}

The spatial and temporal distribution of delivered nutrient concentration
is shown in Fig.~\ref{fig:Static_cox}. In the low growth regime
(left panel) oxygen consumption is practically equal to zero
whereas in the high growth regime (right panel) it becomes
clearly visible
in correspondence of the onset of the anisotropic stress state (Fig.~\ref{fig:Static_csi}, right panel)
and of the drastic increase of proliferating cells
(Fig.~\ref{fig:Static_phi}, bottom left). The marked nonuniformity
in the spatial distribution of oxygen in Fig.~\ref{fig:Static_cox}, right
panel, is to be ascribed to the rather different
boundary conditions at $x=0$ (null diffusive flux) and at $x=L$
(nutrient concentration in local equilibrium with the externally
supplied value).

Fig.~\ref{fig:Static_phi_cthr} shows the time evolution
of cell populations and ECM in correspondence of a reduced
(but still high) level of nutrient concentration
($c_\text{ext}=c_\text{thr}$).
In the low growth regime the evolution of cell populations
is very similar to the case
where $c_\text{ext}=c_\text{sat}$ (Fig.~\ref{fig:Static_phi_cthr},
left) while the ECM concentration undergoes a clear
decrease due to the fact that the
supply of growth factors is not sufficient to
promote both cell mitosis and biosynthesis.
In the high growth regime
the $\n$-cell profile is very similar to the profile in the
low growth regime, except
for a delayed decrease, while
synthesizing cells exhibit a markedly visible increase (Fig.~\ref{fig:Static_phi_cthr}, right panel),
reaching approximately the same maximum value as that of
proliferating cells. This interesting biophysical result
may be ascribed to the increased growth rate $k_\text{g2}$ that
promotes mitotic activity at the expense of a very low oxygen consumption.
As a consequence, the nutrient contained in the fluid is made
completely available for secreting cells that
increase their density within the growing tissue and, accordingly,
their biomass production.

\subsubsection{Initial condition IC2}

In this second set of simulations, the initial
seeding density of cells increase by an order of magnitude
with respect \textcolor{black}{to the initial density in condition IC1.}
\begin{figure}[h!]
\begin{minipage}[c]{1\textwidth}
{\includegraphics[width=0.5\textwidth]{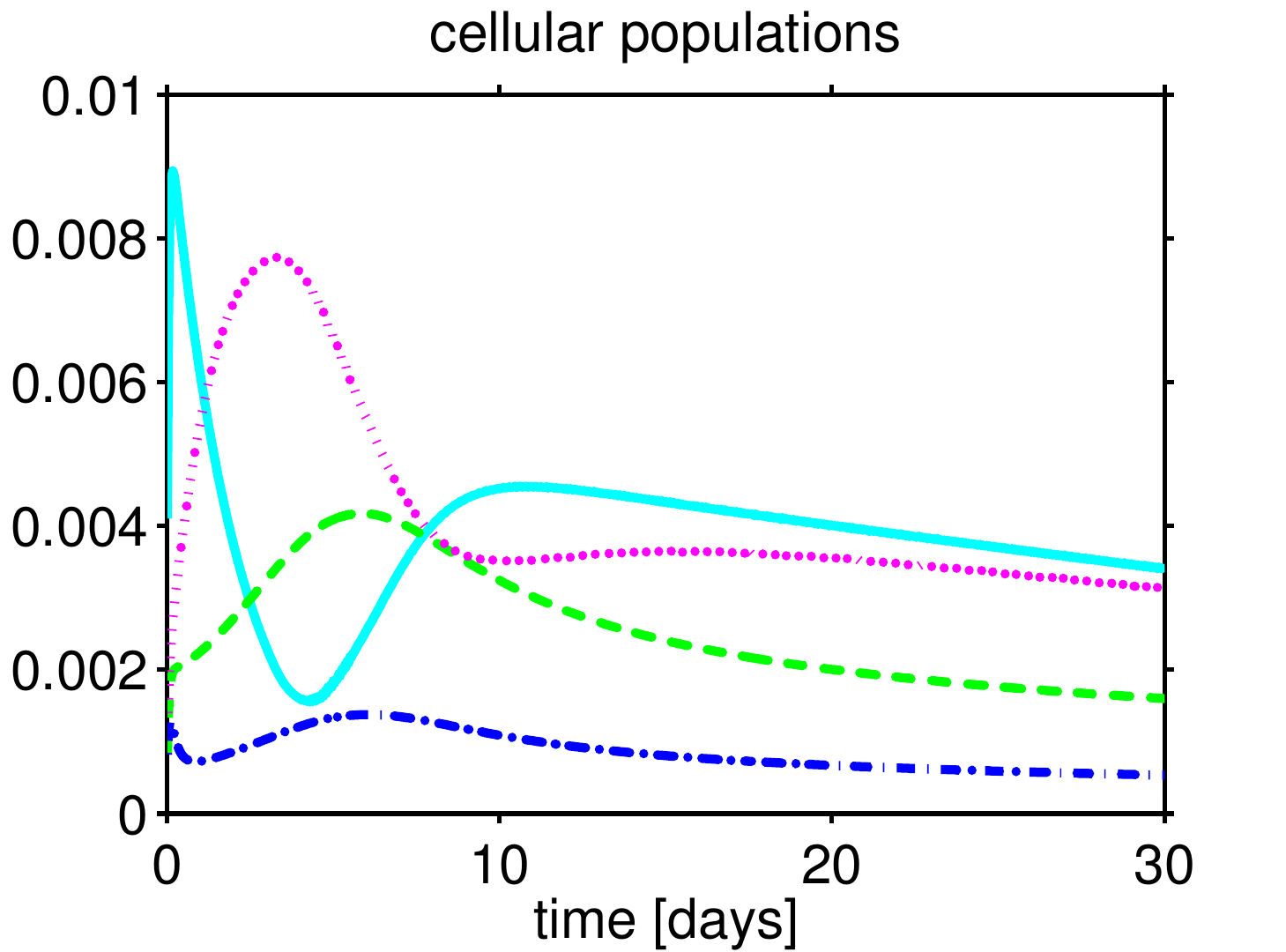}}
{\includegraphics[width=0.5\textwidth]{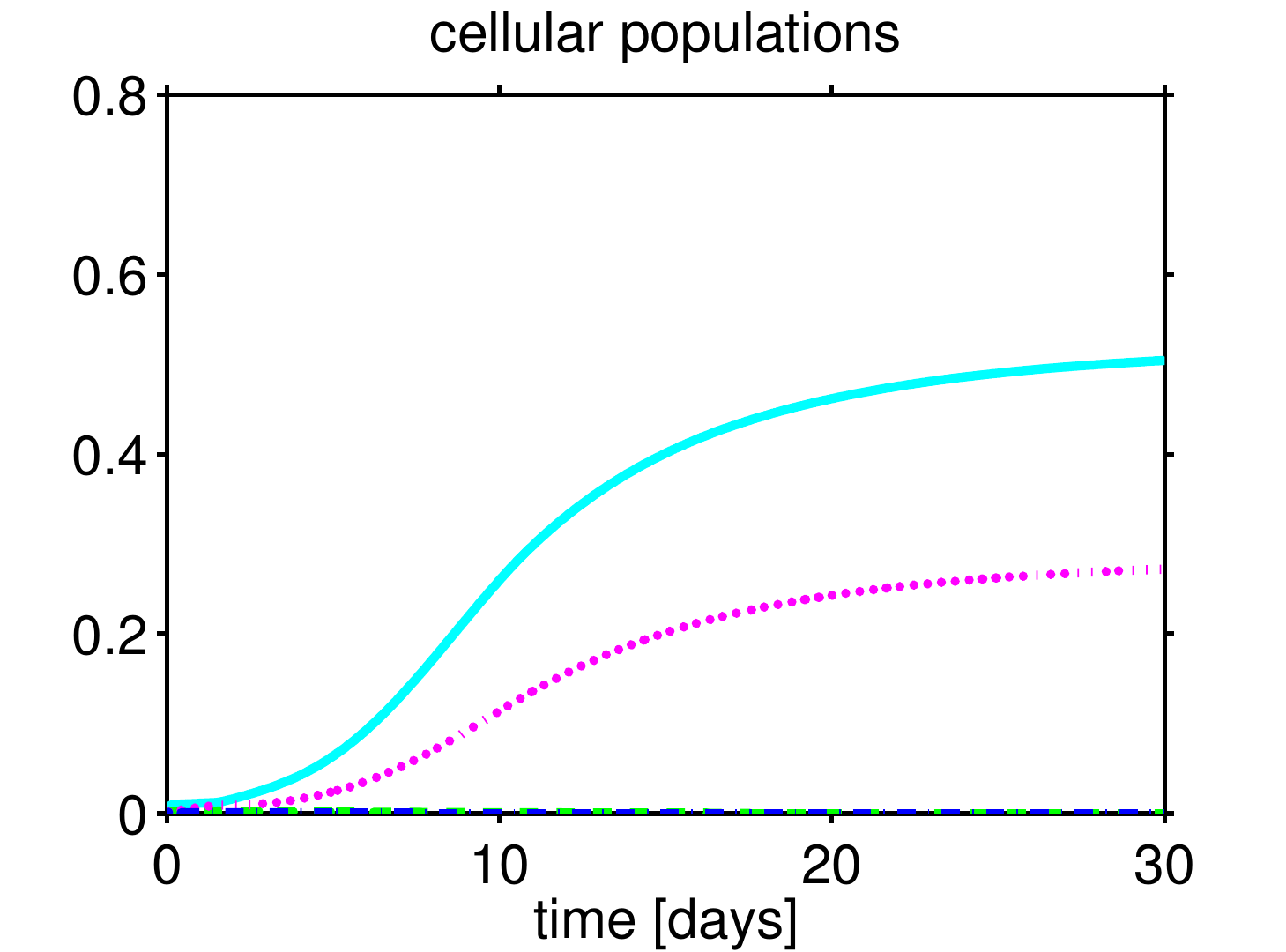}}
\end{minipage}
\caption{Temporal evolution of cellular populations and ECM in the static culture for $c_\text{ext}=c_\text{sat}$.
Initial condition IC2. Left: $k_\text{g} = k_\text{g1}$. Right: $k_\text{g} = k_\text{g2}$. \rs{Solid line: $\phi_\n$; dashed line: $\phi_\v$;
dotted line: $\phi_\q$; dash-dot line: $\phi_\ECM$.}}
\label{fig:Static_phi_IC2}
\end{figure}

The resulting response of the system 
\textcolor{black}{in the low growth regime} displays
an oscillatory evolution of the solid volumetric fractions
before approaching the equilibrium stable state  $(\phi_\n,\phi_\v,\phi_\q,\phi_\ECM)=(0,0,0,0)^T$.
In particular, as shown in Fig.~\ref{fig:Static_phi_IC2}, left panel, proliferating cells, in addition to the initial spike 
\textcolor{black}{already visible in the case IC1}, 
exhibit a second (smaller) growth peak at about 10 days of culture, and then start again to slowly decrease.
This fluctuation may be interpreted as an oscillation around a mean value, that represents the average level of proliferation
measured in the first two weeks of the experiments.
As a matter of fact, as shown in~\cite{Issa}, at high seeding density no significant change in cell number was measured after 2 weeks of culture
(see Fig.~3c of~\cite{Issa}).
\textcolor{black}{Model predictions} suggest
that high initial seeding densities might negatively affect biomass growth, because, on the one hand,
increased cell-cell contact might inhibit the formation of new colonies, and, on the other hand, because nutrient availability, although being abundant at the beginning of the culture process, becomes insufficient so that
cellular metabolism looses its functionality.
%analyzed the effect of using different seeding densities in mechanostimulated bioreactors and

\begin{figure}[h!]
\begin{minipage}[c]{1\textwidth}
{\includegraphics[width=0.5\textwidth]{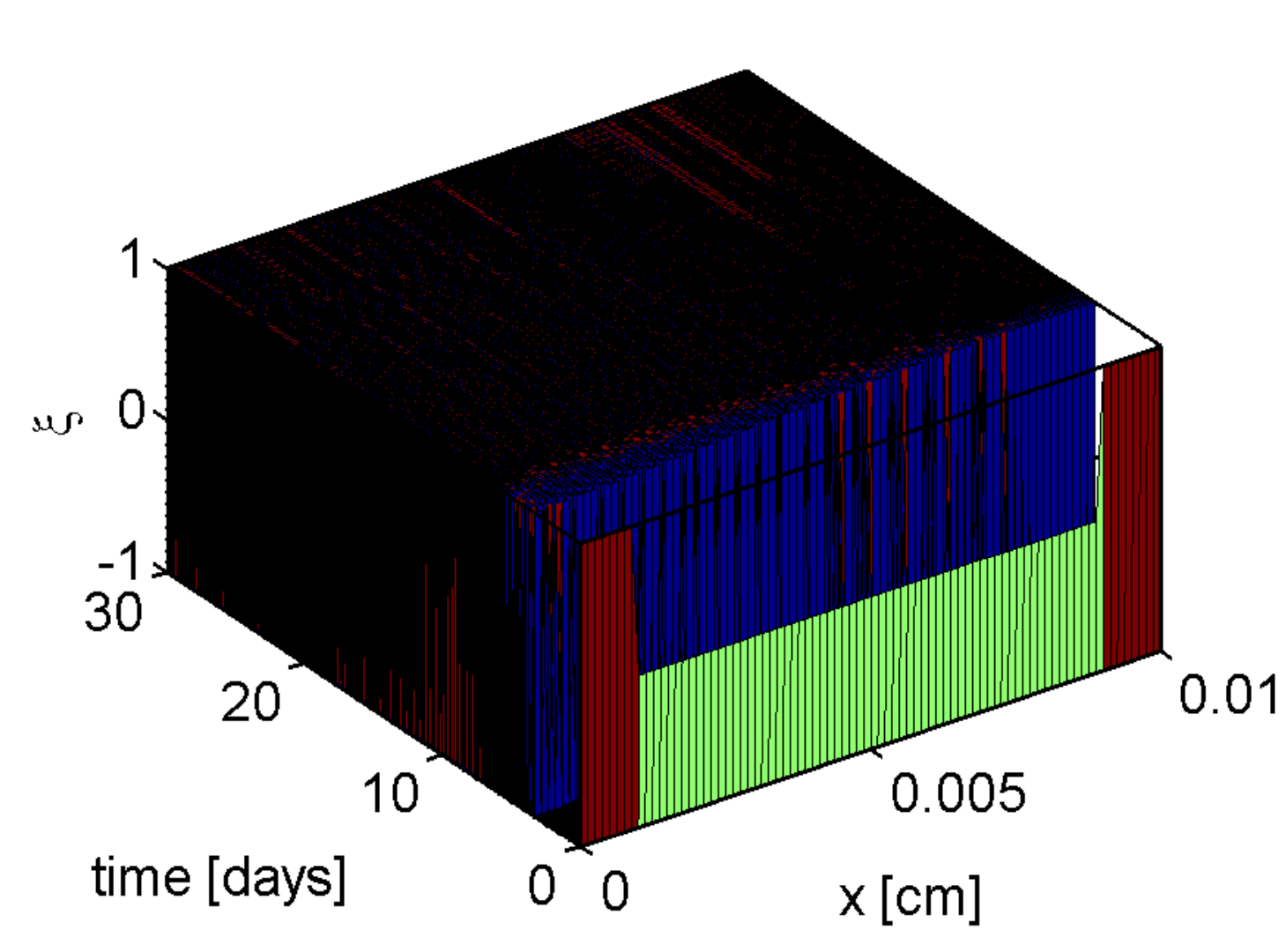}}
{\includegraphics[width=0.5\textwidth]{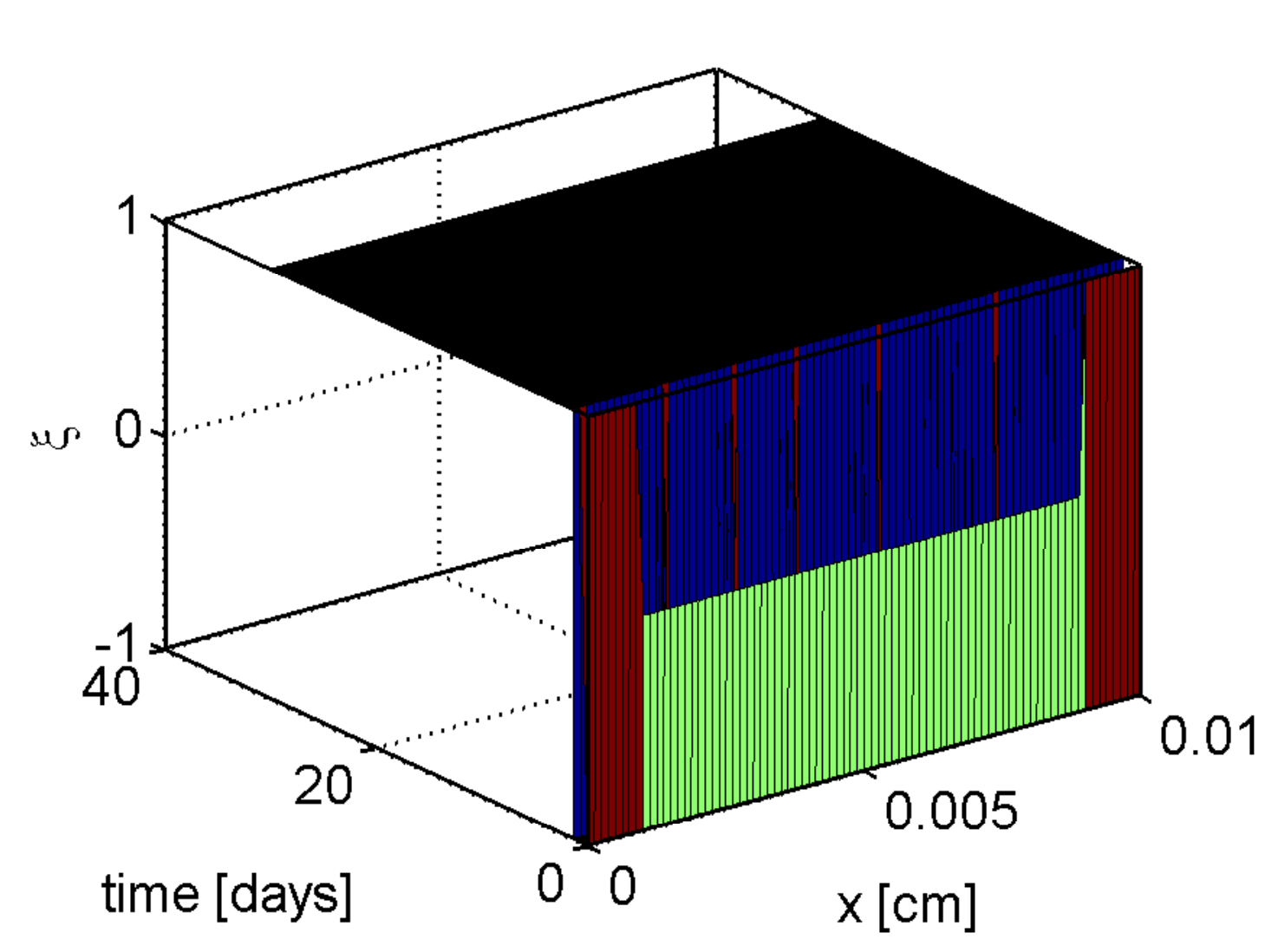}}
\end{minipage}
\caption{Spatial and temporal evolution
of parameter $\xi$ in the static culture for $c_\text{ext}=c_\text{thr}$. Initial condition IC2. Left: $k_\text{g} = k_\text{g1}$. Right: $k_\text{g} = k_\text{g2}$.}
\label{fig:Static_csi_IC2}
\end{figure}

\begin{figure}[h!]
\begin{minipage}[c]{1\textwidth}
\begin{center}
{\includegraphics[width=0.45\textwidth]{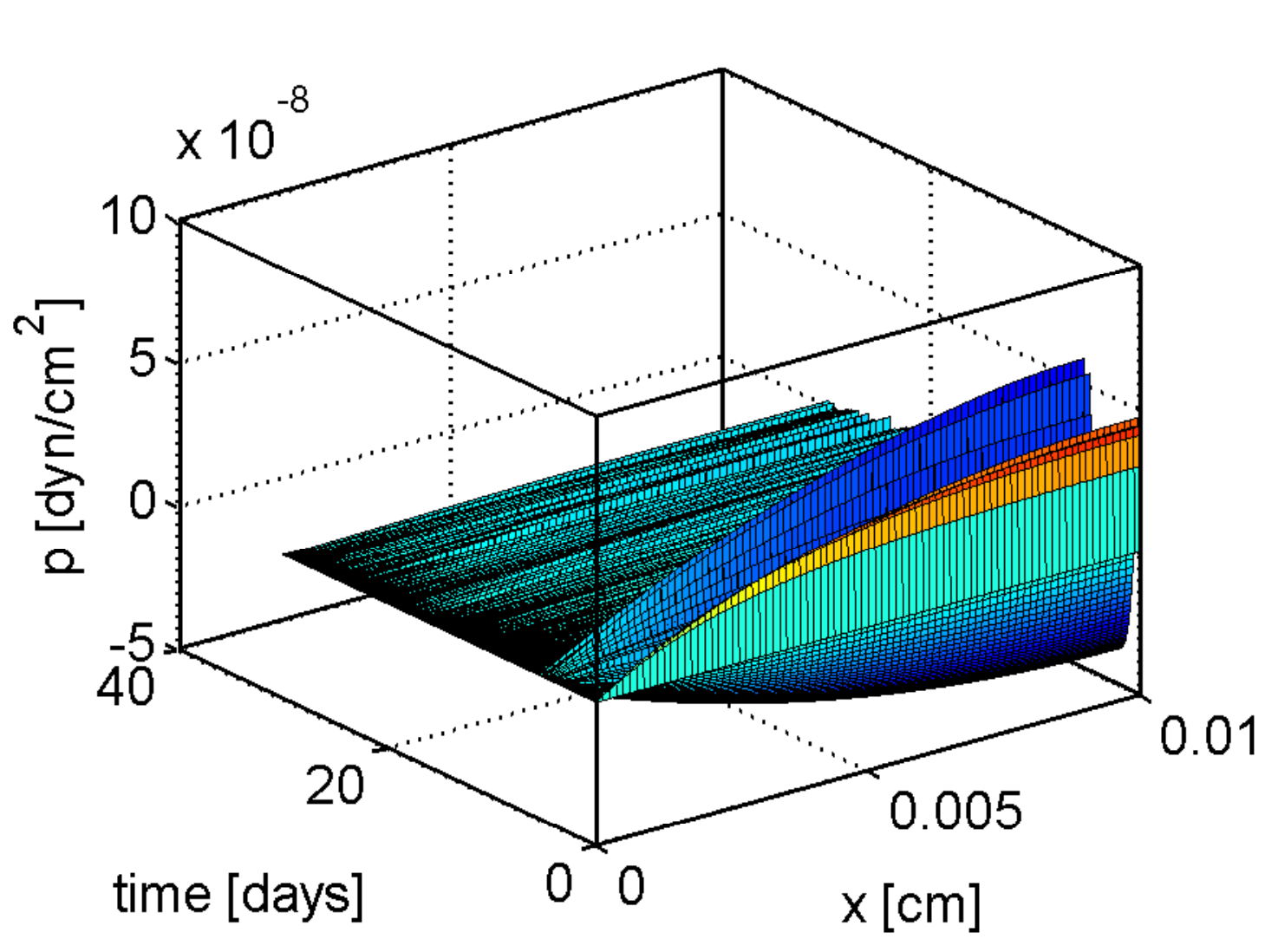}}\qquad
{\includegraphics[width=0.45\textwidth]{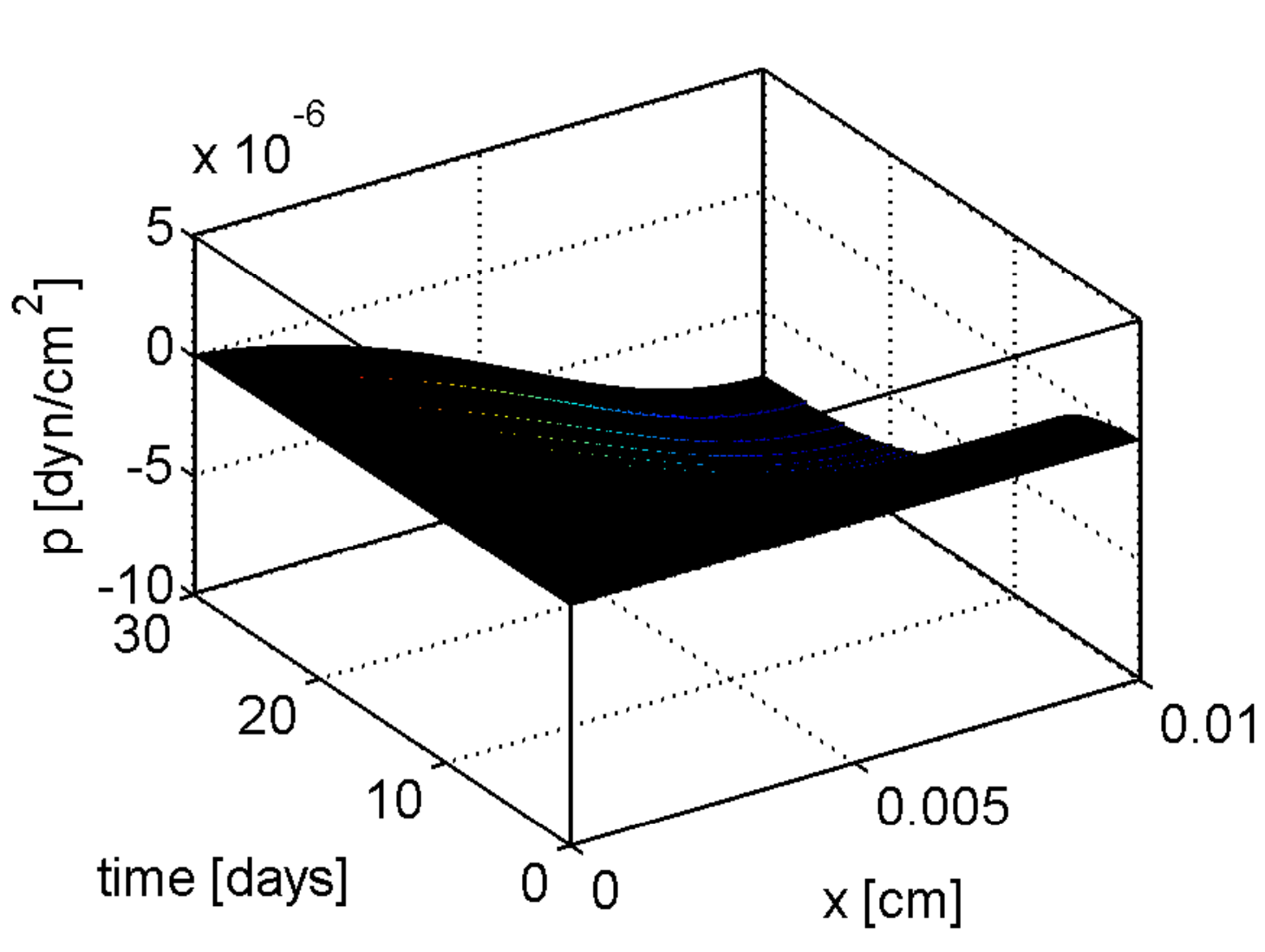}}
\end{center}
\end{minipage}
\caption{Spatial and temporal evolution of fluid pressure in
the static culture for $c_\text{ext}=c_\text{sat}$.
Initial condition IC2. Left: $k_\text{g}=k_\text{g1}$. Right: $k_\text{g}=k_\text{g2}$.}
\label{fig:Static_fl_IC2}
\end{figure}

The oscillatory behavior of proliferating cells shown in Fig.~\ref{fig:Static_phi_IC2} left is manifested also
in the time evolution of $\xi$ (Fig.~\ref{fig:Static_csi_IC2}, left panel)
and $p$ (Fig.~\ref{fig:Static_fl_IC2}, left panel).
The mechanical stress in the biomass oscillates
around the isotropic state while the fluid pressure oscillates
around a constant value,
and these fluctuations persist until the end of the simulation.
\begin{figure}[h!]
\begin{minipage}[c]{1\textwidth}
{\includegraphics[width=0.5\textwidth]{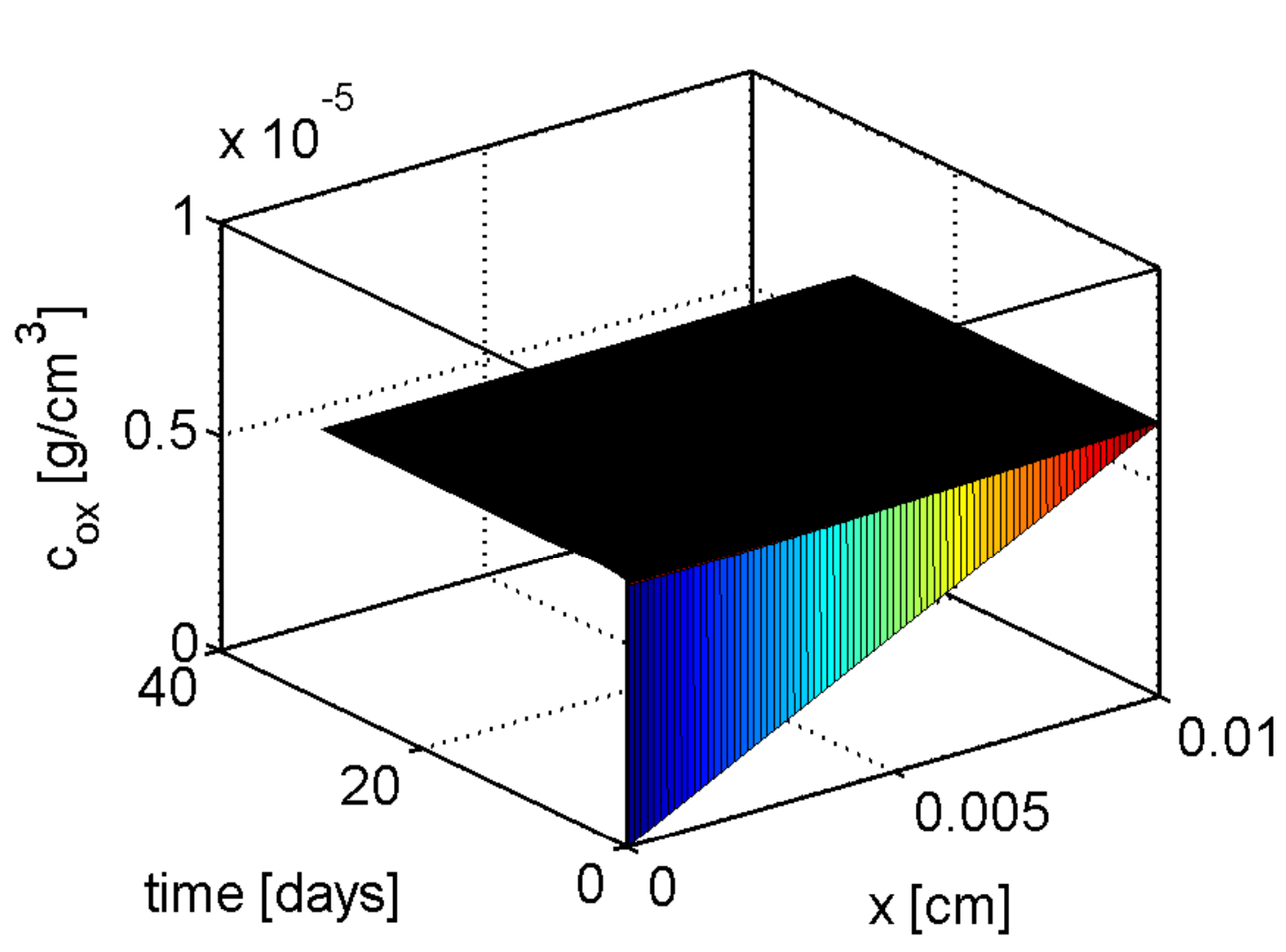}}
{\includegraphics[width=0.5\textwidth]{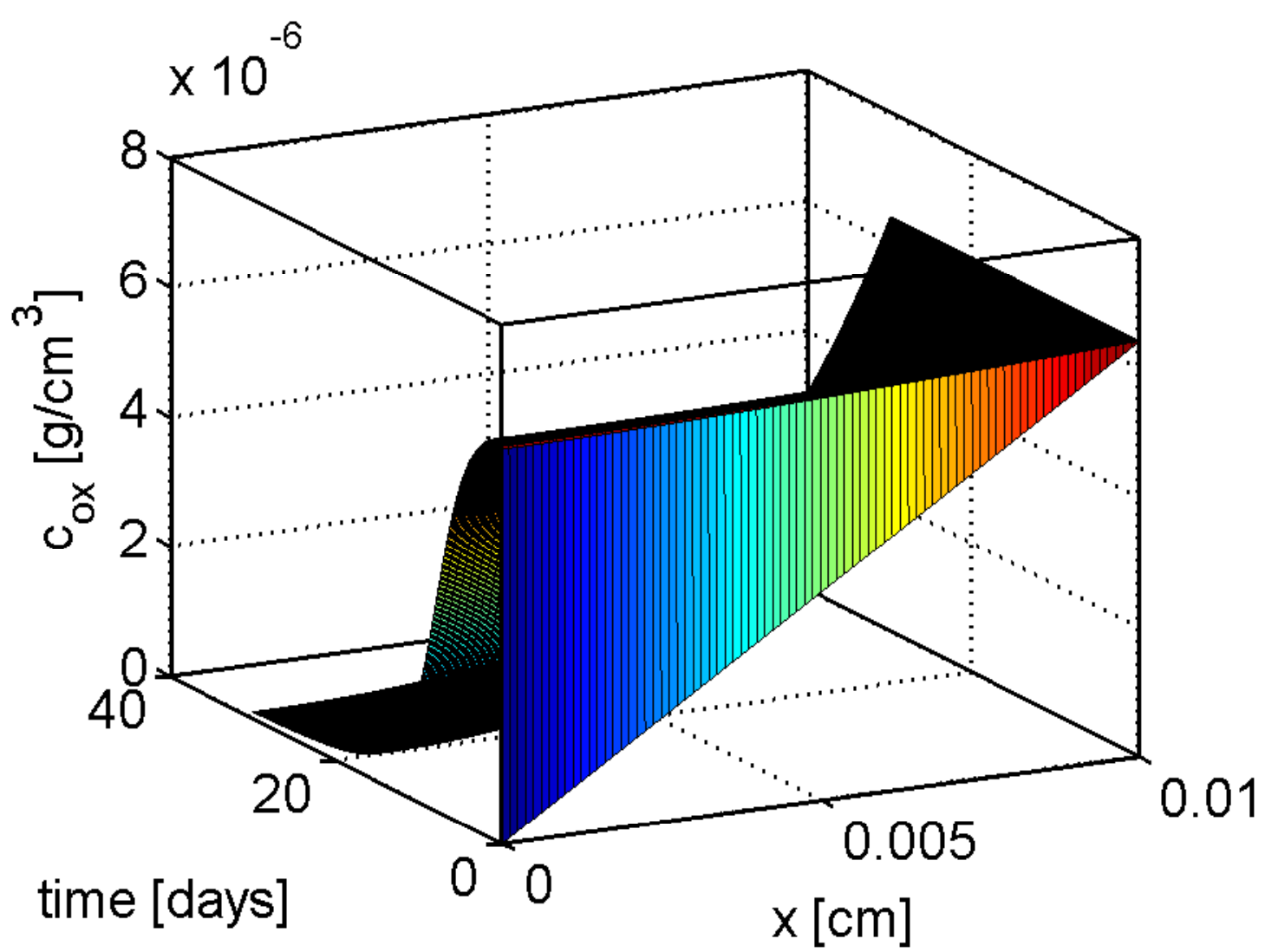}}
\end{minipage}
\caption{Spatial and temporal evolution
of oxygen concentration in the static culture for
$c_\text{ext}=c_\text{sat}$.
Initial condition IC2. Left: $k_\text{g} = k_\text{g1}$. Right: $k_\text{g} = k_\text{g2}$.}
\label{fig:Static_cox_IC2}
\end{figure}

\begin{figure}[h!]
\begin{minipage}[c]{1\textwidth}
{\includegraphics[width=0.5\textwidth]{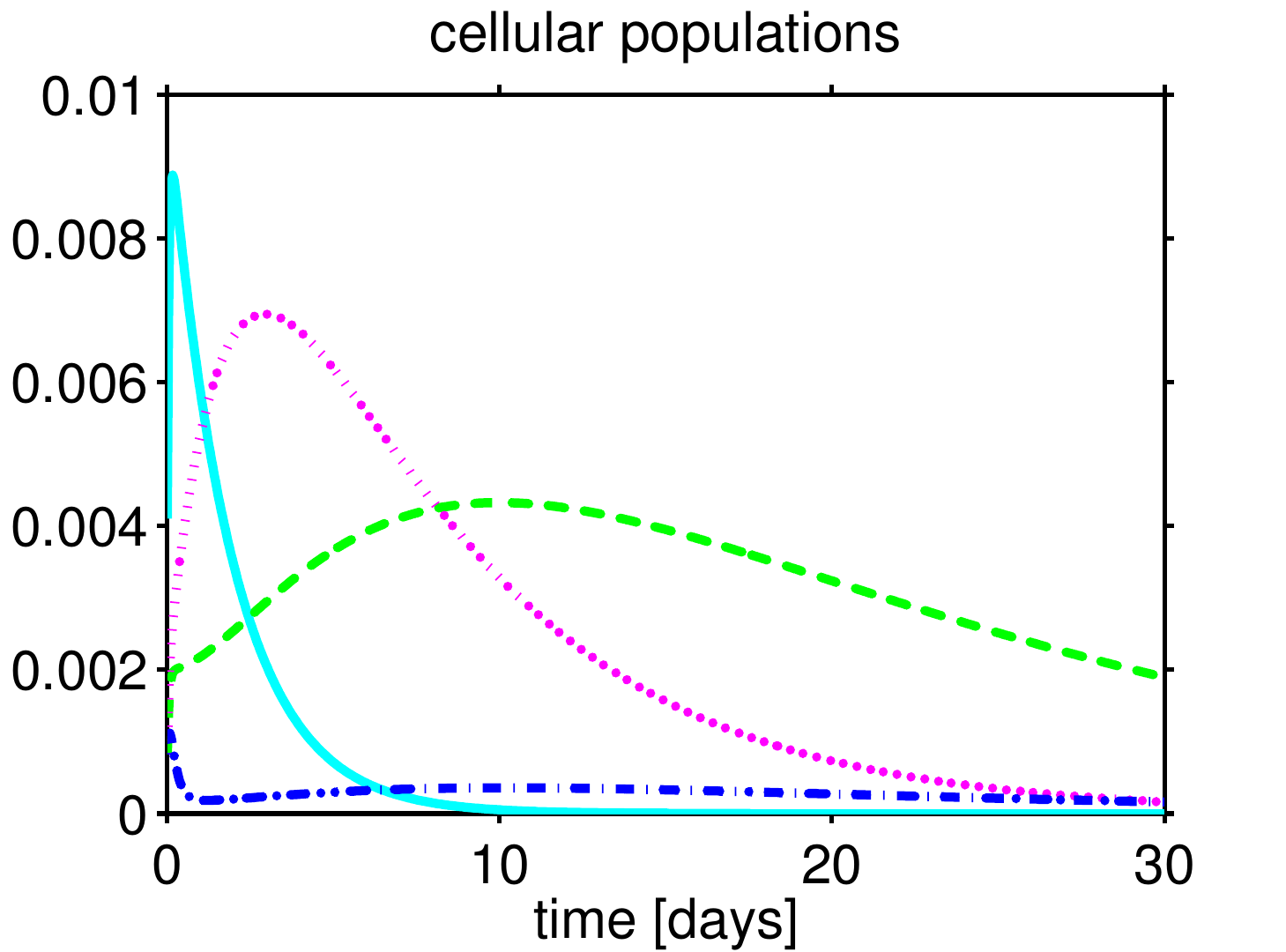}}
{\includegraphics[width=0.5\textwidth]{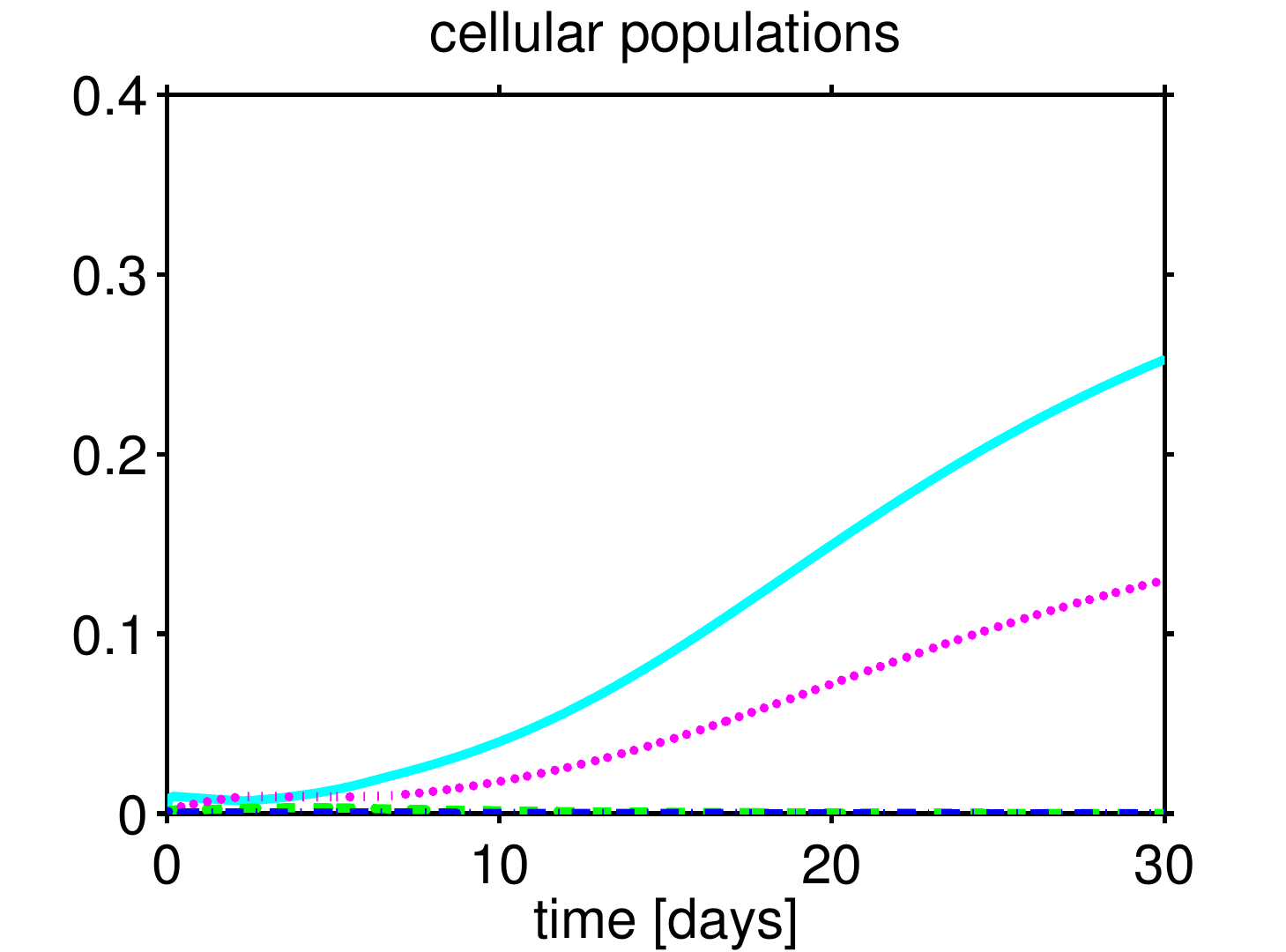}}
\end{minipage}
\caption{Temporal evolution of cellular populations and ECM in the static culture for $c_\text{ext}=c_\text{thr}$. Initial condition IC2.
Left: $k_\text{g} = k_\text{g1}$. Right: $k_\text{g} = k_\text{g2}$.
\rs{Solid line: $\phi_\n$; dashed line: $\phi_\v$;
dotted line: $\phi_\q$; dash-dot line: $\phi_\ECM$.}}
\label{fig:Static_phi_IC2_cthr}
\end{figure}

In the high growth regime, the behavior of cells and ECM
(Fig.~\ref{fig:Static_phi_IC2}, right panel)
is very similar to that of Sect.~\ref{sec:Initial_condition_IC1}
(cf. Fig.~\ref{fig:Static_phi}, left bottom panel).
We notice that the switch from isotropic to anisotropic stress
conditions predicted by parameter $\xi$
(Fig.~\ref{fig:Static_csi_IC2}, right panel) manifests earlier
than in the case of initial conditions IC1
(cf. Fig.~\ref{fig:Static_csi}, right). This different behavior
is probably
to be ascribed to the fact that the initial high seeding density
of proliferating cells, with their spread elongated shape, favors
an anticipated occurrence of the anisotropic adherence state.

The spatial and temporal distribution of fluid pressure in the case $k_\text{g} = k_\text{g2}$ (Fig.~\ref{fig:Static_fl_IC2}, right panel) is
qualitatively similar to the corresponding distribution
in the case of initial conditions IC1 (cf. Fig.~\ref{fig:Static_fl}, right panel). In quantitative terms, the overall fluid pressure drop in the
case of initial conditions IC2 is slightly smaller than in the case IC1
because of the smoother (temporal) transition between a
higher permeability to a lower permeability. Also, it can be noticed
that the (temporal) onset of fluid motion (proportional to fluid pressure
gradient through Darcy's law) occurs much earlier than in the case of IC1
initial conditions because of the same biophysical reasons discussed
above for proliferating cell distribution. Similar arguments apply to
the spatial and temporal distribution of oxygen
concentration in both low and high
growth regimes (Fig.~\ref{fig:Static_cox_IC2}).

We conclude the analysis of the growth process in the case of
a static culture by showing in Fig.~\ref{fig:Static_phi_IC2_cthr}
the evolution of cells and ECM when the
external level of nutrient concentration is reduced.
When $k_\text{g} = k_\text{g1}$ the small nutrient availability
causes the oscillations occurring in the case
$c_\text{ext}=c_\text{sat}$ to vanish and
the $\n$-cell population to rapidly extinguish.
On the other hand, in the high growth regime,
the system is weakly affected by the decreased nutrient availability and
it qualitatively behaves as in the case $c_\text{ext}=c_\text{sat}$, whereas,
from a quantitative point of view, the value reached by $\phi_\n$ at the
end of the simulation is smaller.

\subsection{Perfused culture}
\textcolor{black}{The second group of results refer to the case 
$T_\text{b}=100\, \unit{mPa}$ and $V_\text{b}=50\,\mu\unit{m s^{-1}}$. 
These quantities are the characteristic values of the 
external shear stress and inlet velocity that are found in the 
experimental setup of TE applications (see~\cite{Raimondi2006,Sacco,Causin} and references cited therein). The 
most notable feature that distinguishes biomass growth in such a dynamic culture is that the mitotic function and cellularity are 
in general favored because} the fluid-induced shear stress
strongly contributes to develop anisotropic loads in the construct~\cite{Raimondi2006a,Nava2012}.

\subsubsection{Initial condition IC1}

\begin{figure}[h!]
\begin{minipage}[c]{1\textwidth}
\centering
{\includegraphics[width=0.45\textwidth]{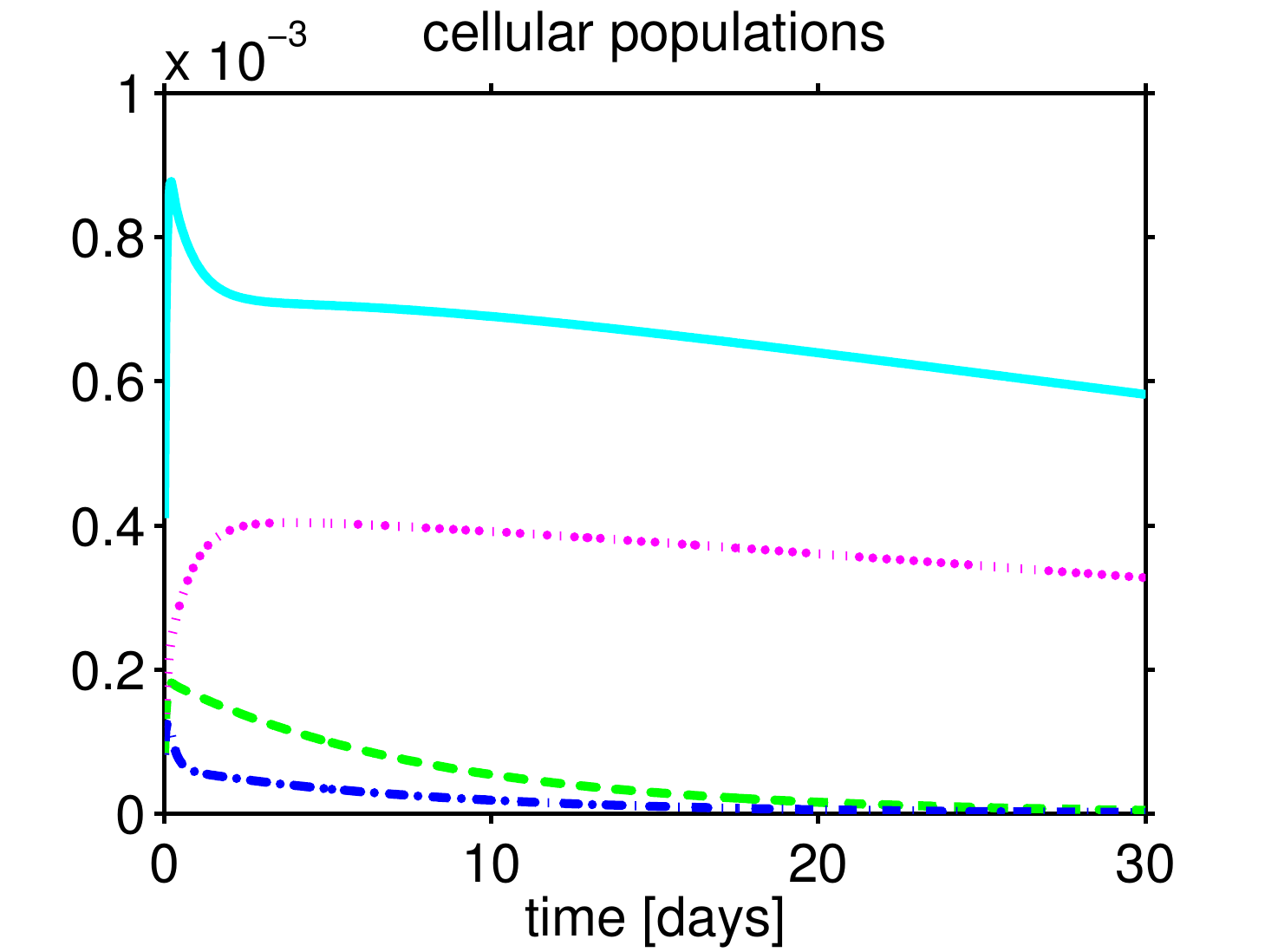}}
{\includegraphics[width=0.45\textwidth]{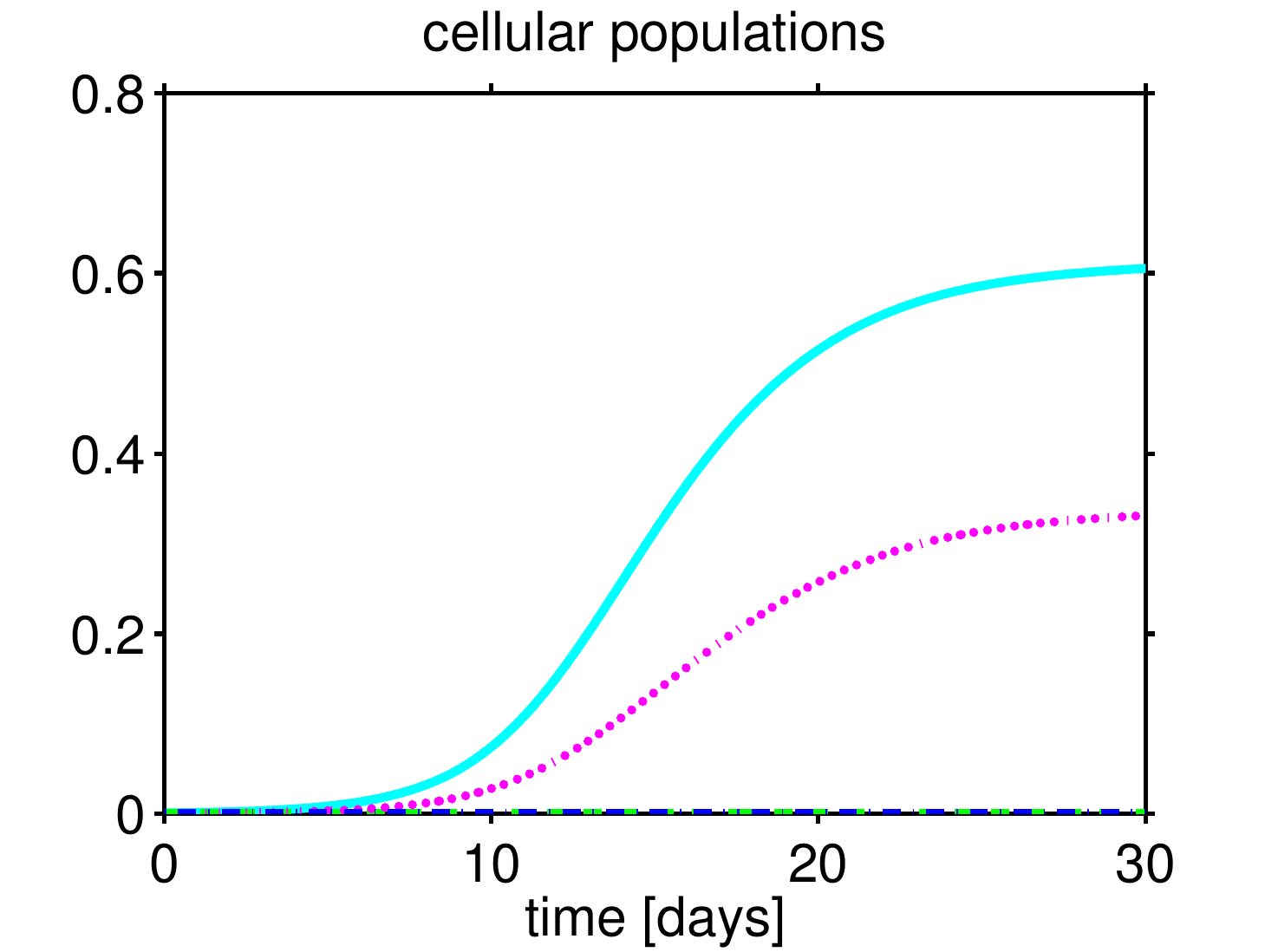}}
\end{minipage}
\caption{Temporal evolution of cellular populations and ECM in the perfused culture for $c_\text{ext}=c_\text{sat}$. Initial condition IC1. Left: $k_\text{g} = k_\text{g1}$. Right: $k_\text{g} = k_\text{g2}$.
\rs{Solid line: $\phi_\n$; dashed line: $\phi_\v$;
dotted line: $\phi_\q$; dash-dot line: $\phi_\ECM$.}
}
\label{fig:Perfused_phi}
\end{figure}

\begin{figure}[h!]
\begin{minipage}[c]{1\textwidth}
\centering
{\includegraphics[width=0.45\textwidth]{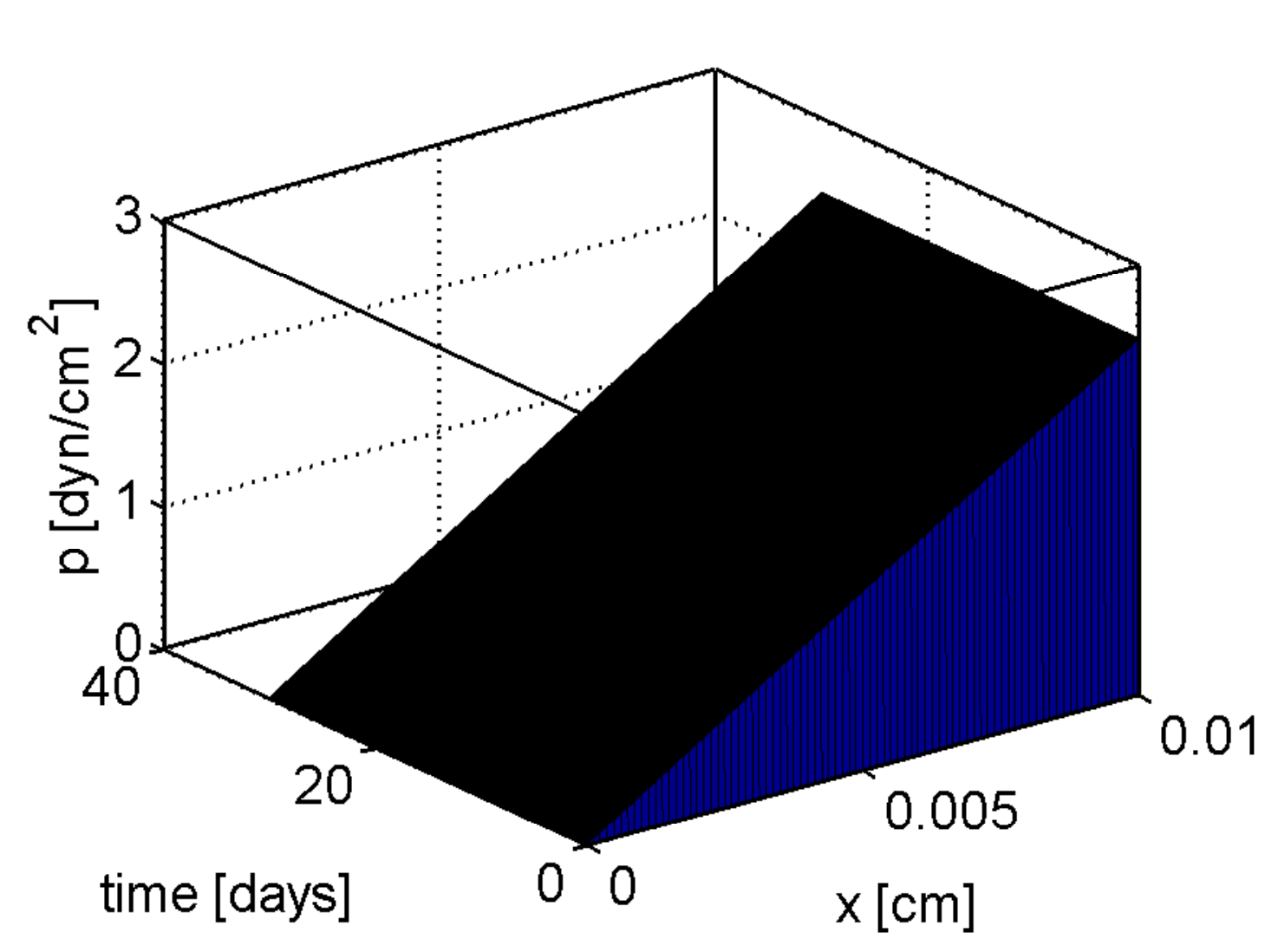}}
{\includegraphics[width=0.45\textwidth]{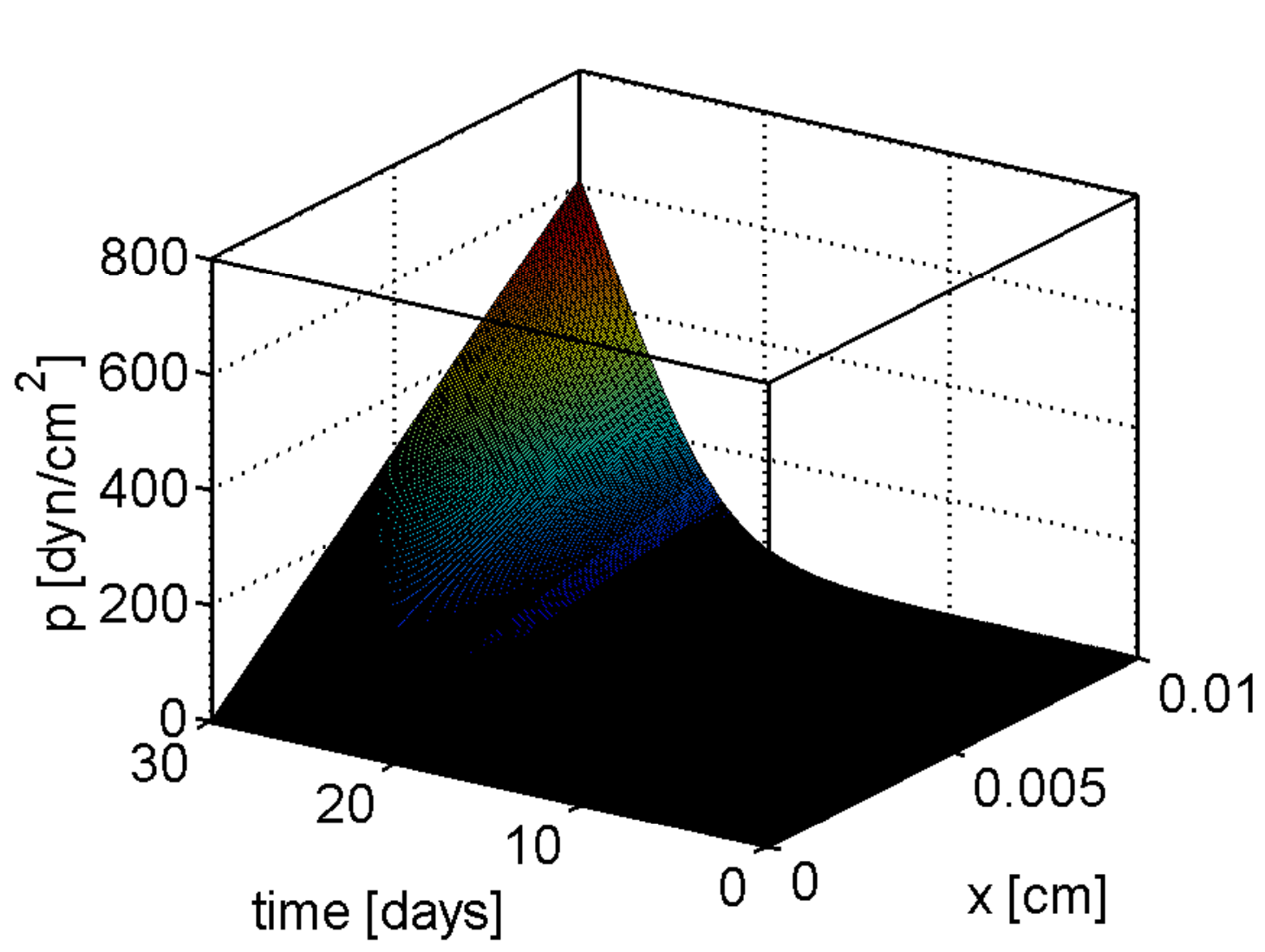}}
\end{minipage}
\caption{Spatial and temporal evolution of fluid pressure in the perfused culture for $c_\text{ext}=c_\text{sat}$. Initial condition IC1.
Left: $k_\text{g}=k_\text{g1}$. Right: $k_\text{g}=k_\text{g2}$.}
\label{fig:Perfused_fl}
\end{figure}

\begin{figure}[h!]
\begin{minipage}[c]{1\textwidth}
{\includegraphics[width=0.5\textwidth]{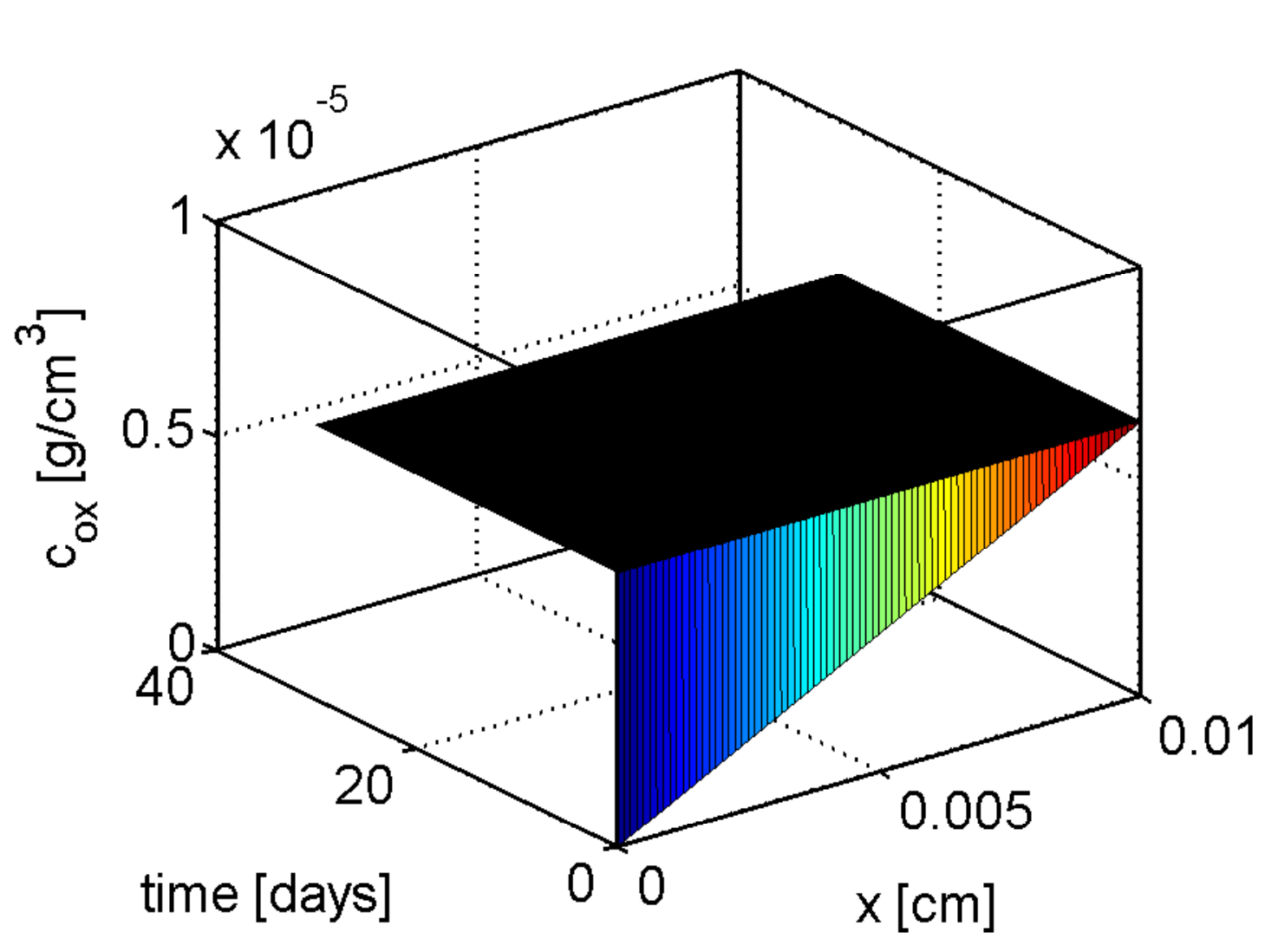}}
{\includegraphics[width=0.5\textwidth]{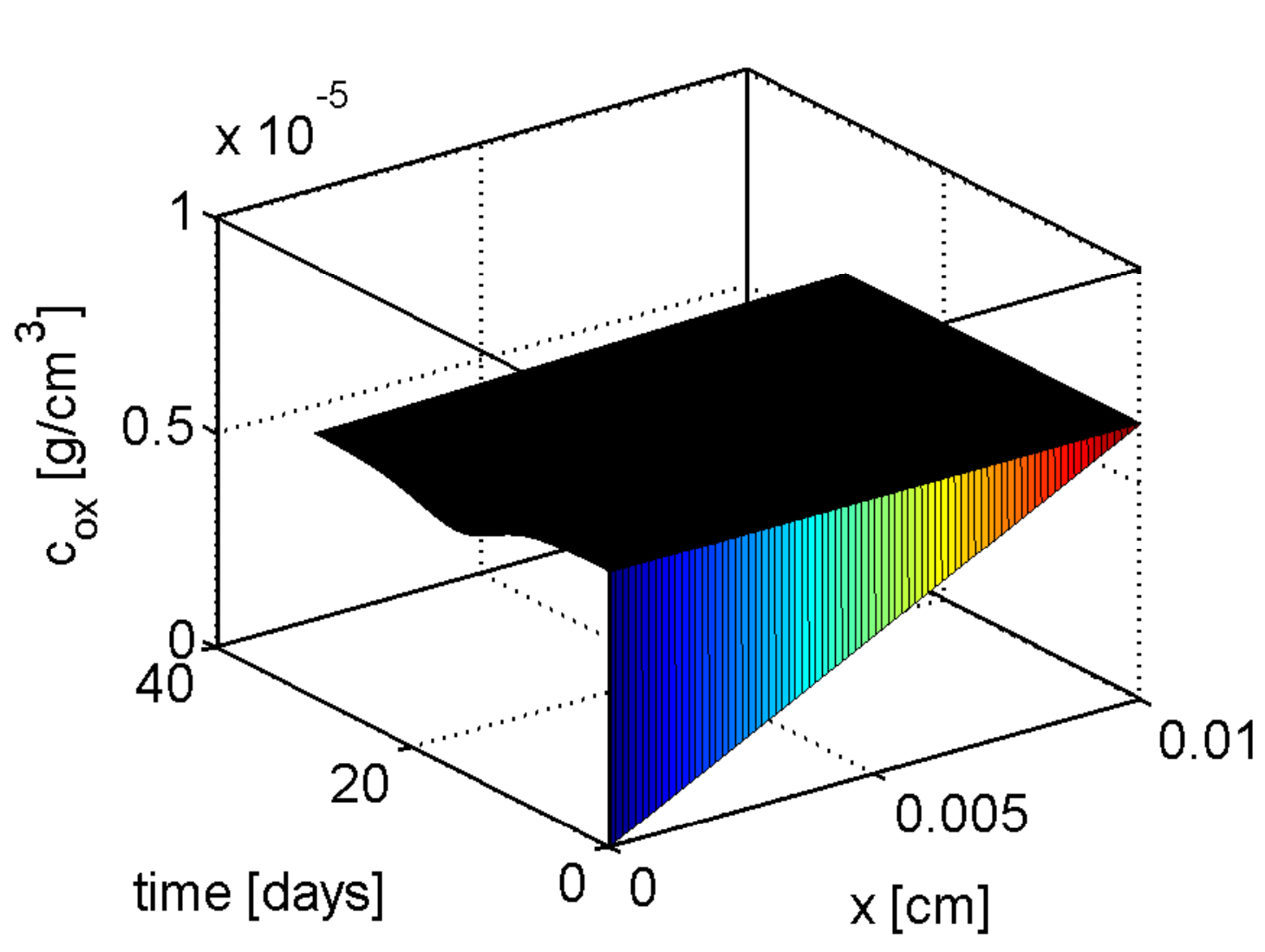}}
\end{minipage}
\caption{Spatial and temporal evolution of oxygen concentration in the perfused culture for $c_\text{ext}=c_\text{sat}$. Initial condition IC1. Left: $k_\text{g} = k_\text{g1}$. Right: $k_\text{g} = k_\text{g2}$.}
\label{fig:Perfused_cox}
\end{figure}

\begin{figure}[h!]
\begin{minipage}[c]{1\textwidth}
\centering
{\includegraphics[width=0.45\textwidth]{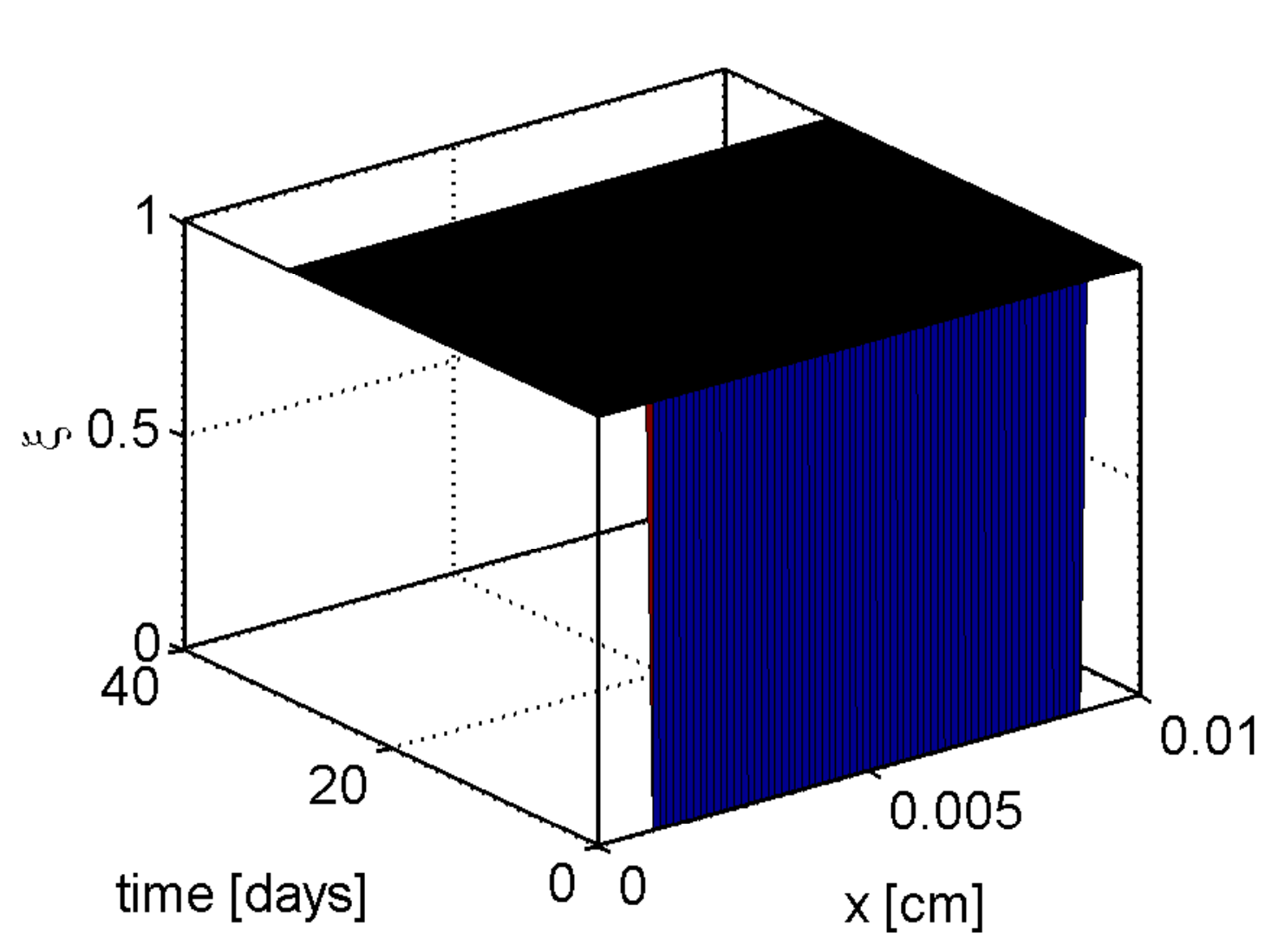}}
{\includegraphics[width=0.45\textwidth]{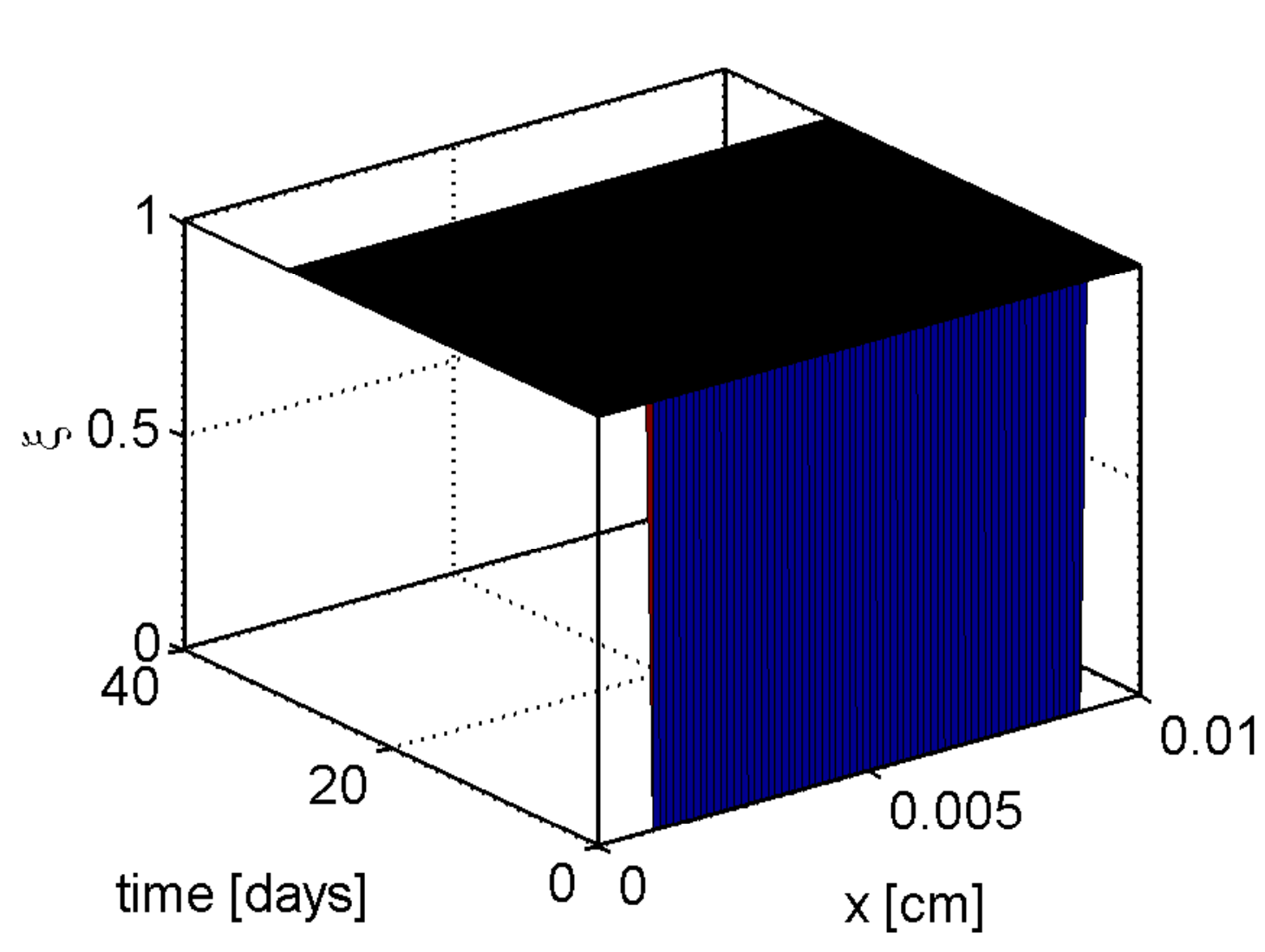}}
\end{minipage}
\caption{Spatial and temporal evolution of parameter $\xi$ in the perfused culture for $c_\text{ext}=c_\text{sat}$. Initial condition IC1. Left:
$k_\text{g} = k_\text{g1}$. Right: $k_\text{g} = k_\text{g2}$.}
\label{fig:Perfused_csi}
\end{figure}

\begin{figure}[h!]
\begin{minipage}[c]{1\textwidth}
\centering
{\includegraphics[width=0.45\textwidth]{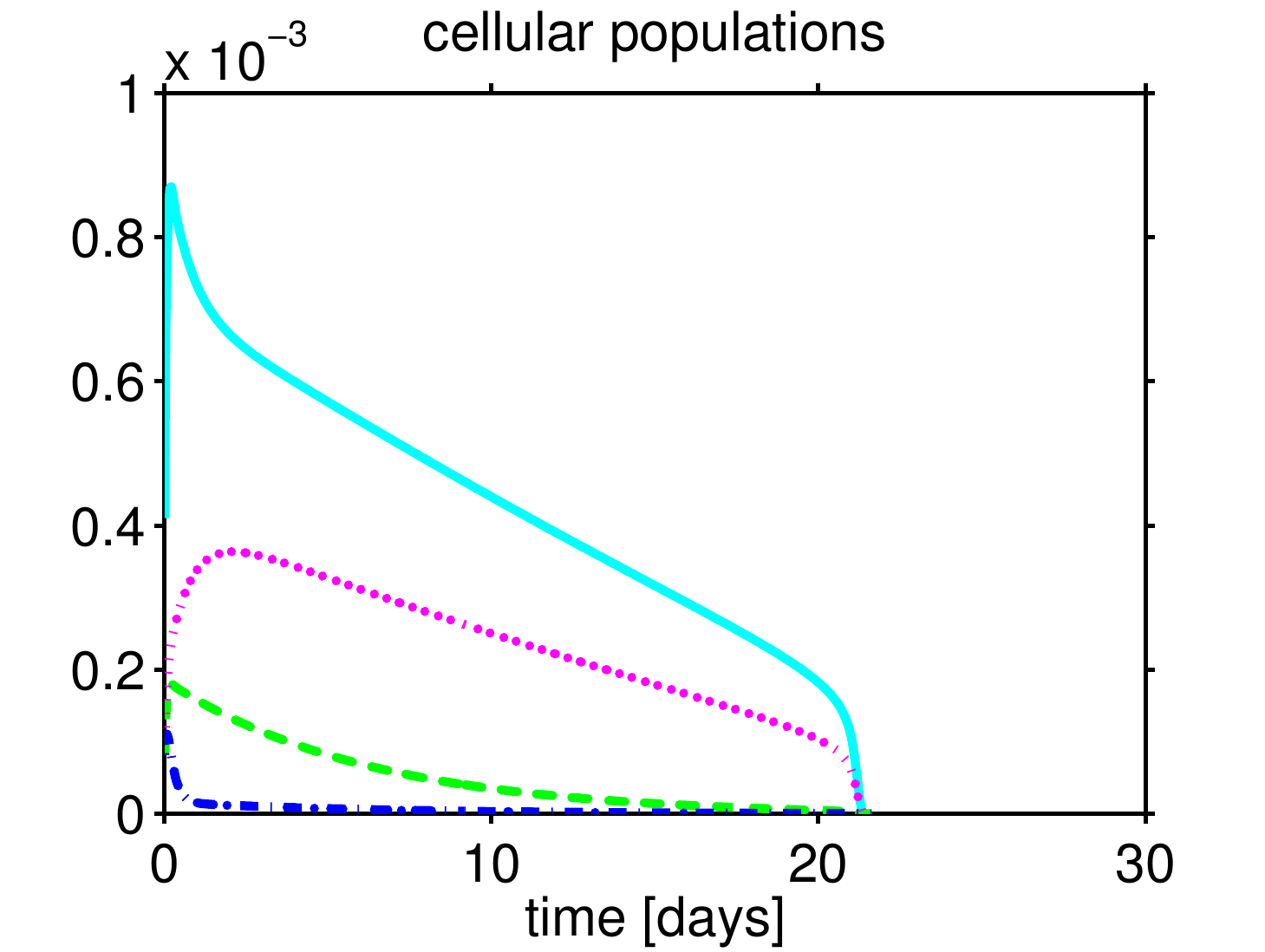}}
{\includegraphics[width=0.45\textwidth]{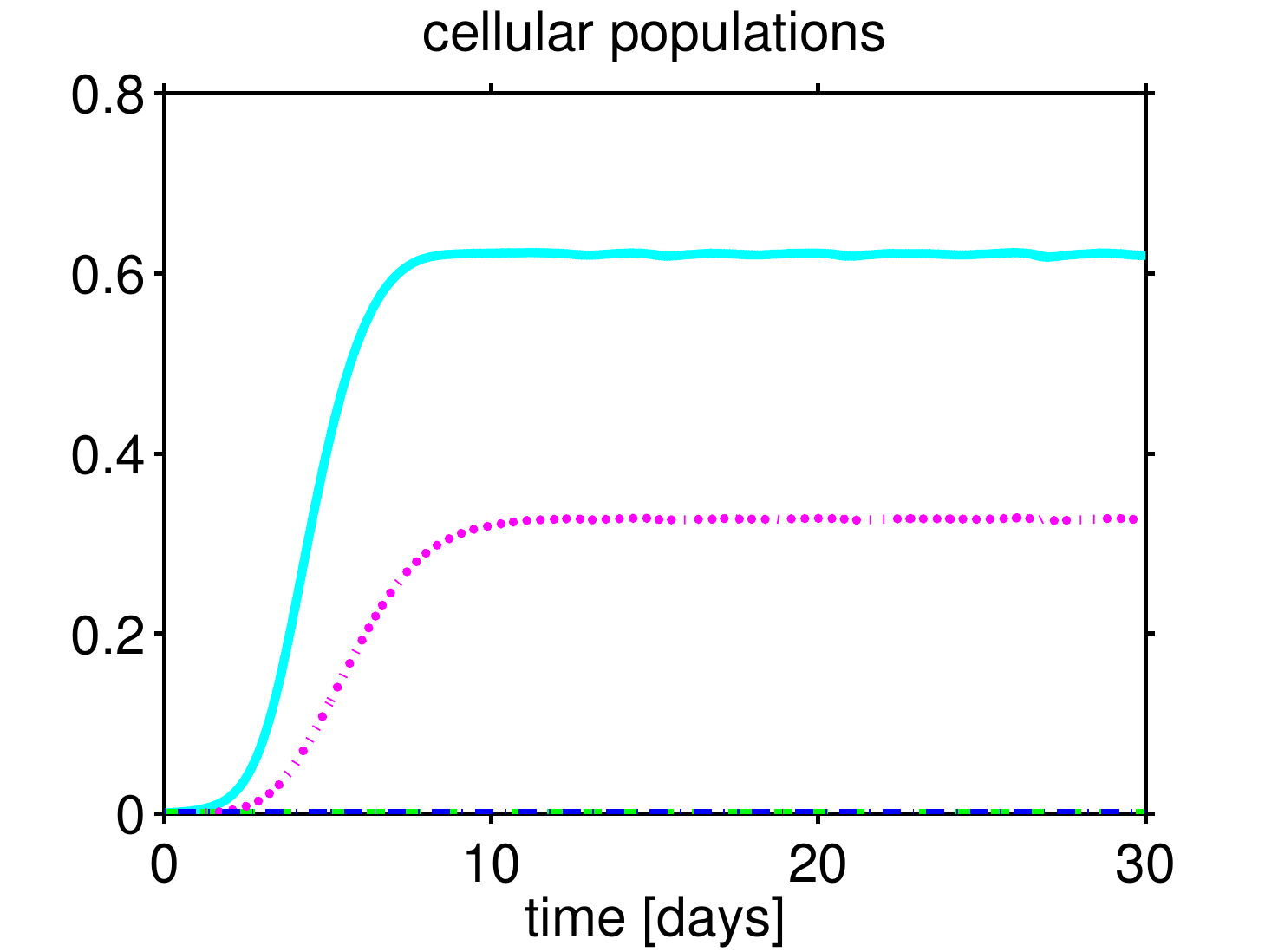}}
\end{minipage}
\caption{Temporal evolution of cellular populations and ECM in the perfused culture for $c_\text{ext}=c_\text{thr}$. Initial condition IC1. Left: $k_\text{g} = k_\text{g1}$. Right: $k_\text{g} = k_\text{g2}$.
\rs{Solid line: $\phi_\n$; dashed line: $\phi_\v$;
dotted line: $\phi_\q$; dash-dot line: $\phi_\ECM$.}
}
\label{fig:Perfused_phi_cthr}
\end{figure}

\textcolor{black}{Fig.~\ref{fig:Perfused_phi} (left panel)}
clarifies how the mitotic activity is preponderant with respect to the other phases, that, differently, very quickly converge to zero. 
As a matter of fact proliferating cells, after the
initial spike typical of the low growth regime, show a plateau-like behaviour that dramatically delays cell apoptosis.
% 
% %As a matter of fact we prove (in the forthcoming paper
% %already recall previously) when the volumetric fraction do not converge to the equilibrium point $(0,0,0,0)$, they tend to another
% %finite steady state.
% 
% %The uncontrolled growth of $\phi_\n$ is probably due to the fact that the anisotropy wave
% %does not arise when the culture has already partially settled, rather it is imposed at the boundary of the domain by the
% %fluid pressure and then induces a strong mechanical stimulus. Consequently a generalized anisotropy load,
% %that holds from the beginning together with the evolution of $\phi_\n$, highlights their exponential increase.
% 
This result highlights the reaction of the system to the mechanical boundary conditions: the external force
exerted by the fluid \textcolor{black}{gives rise to an anisotropic} 
stress state that instantaneously propagates throughout the domain,
as shown in Fig.~\ref{fig:Perfused_csi}, maintaining 
\textcolor{black}{the 
anisotropic mechanical configuration} in both low and high growth regimes.
This mechanical stimulus is the \textcolor{black}{sole} 
responsible of the strong mitotic 
\textcolor{black}{functional activity occurring} 
within the biomass and shown in Fig.~\ref{fig:Perfused_phi} (left panel).
As a matter of fact, as evidenced in Fig.~\ref{fig:Perfused_cox}, oxygen consumption is practically absent
since nutrient concentration remains substantially unchanged 
\textcolor{black}{except the very small sink 
in the case $k_\text{g} = k_\text{g2}$ at $x=0 \, \unit{cm}$ and 
$t = 15$ days.}

The spatial and temporal evolution of the fluid pressure is displayed in Fig.~\ref{fig:Perfused_fl}.
In the low growth regime (Fig.~\ref{fig:Perfused_fl} left panel) 
fluid pressure shows a Darcy-like behavior
and linearly increases along the spatial domain. 
On the contrary, in the high growth regime the
evolution of the pressure becomes nonlinear and, at the end of the simulation
at the \textcolor{black}{rightmost end} 
of the domain, it reaches a spike whose value is about
two orders of magnitude larger than the maximum value in the low growth case.
This happens because of the concurrence of two factors: first, at the right boundary
of the domain an external positive pressure is applied, second, at the end of the simulation
the porosity within the biomass decreases and, consequently, the fluid pressure
dramatically increases.

Fig.~\ref{fig:Perfused_phi_cthr} 
shows the response of cellular populations when the external nutrient
level drops below the saturation threshold. If the maximum cell growth rate is not high enough, all cell populations rapidly decrease and disappear shortly after twenty days of simulation (left panel).
This result proves that in a dynamical culture, when the external oxygen concentration is insufficient, cell survival is ensured only 
if the growth rate \textcolor{black}{is large enough (right panel).}

\subsubsection{Initial condition IC2}

\begin{figure}[h!]
\begin{minipage}[c]{1\textwidth}
\centering
{\includegraphics[width=0.45\textwidth]{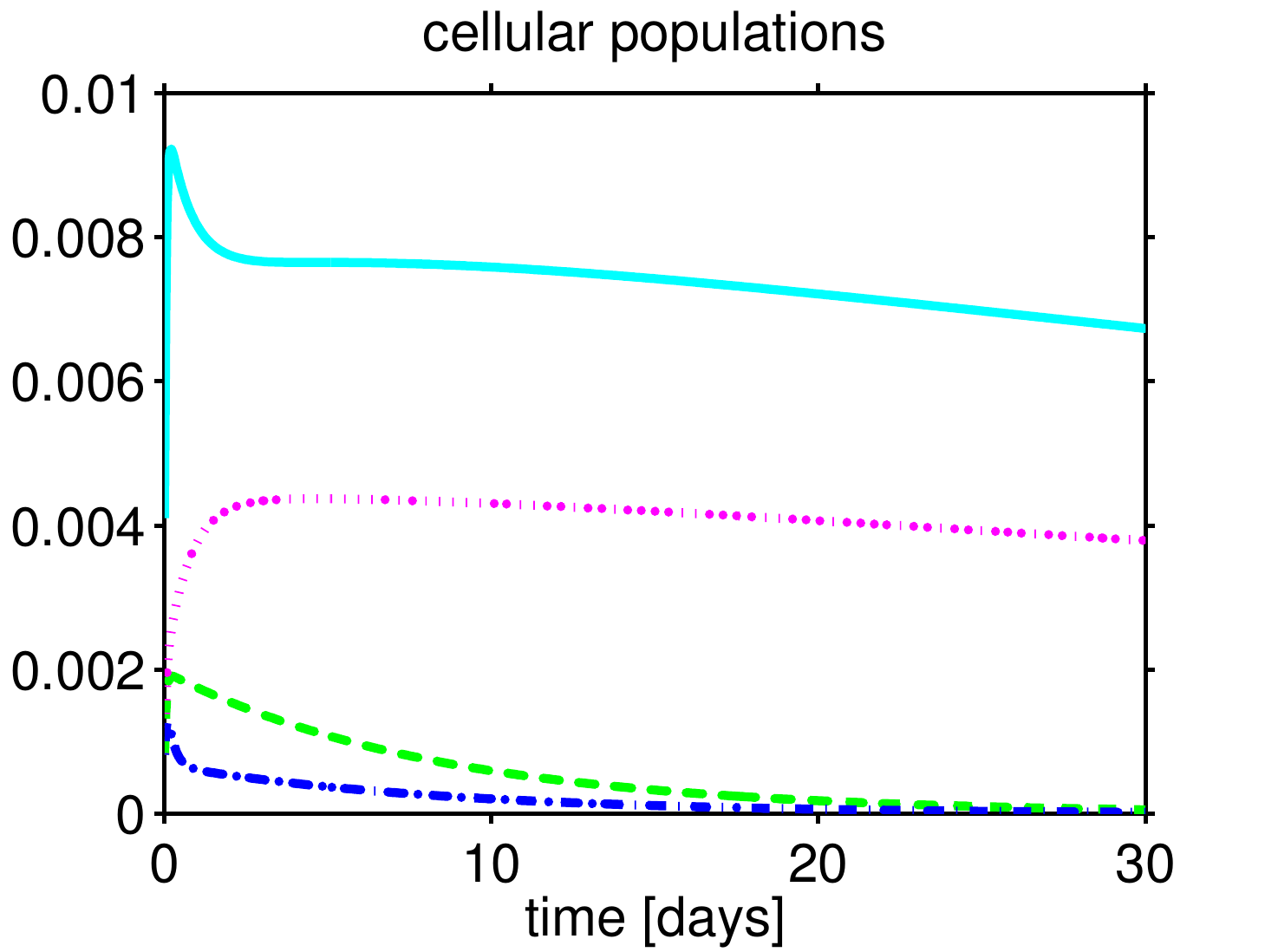}}
{\includegraphics[width=0.45\textwidth]{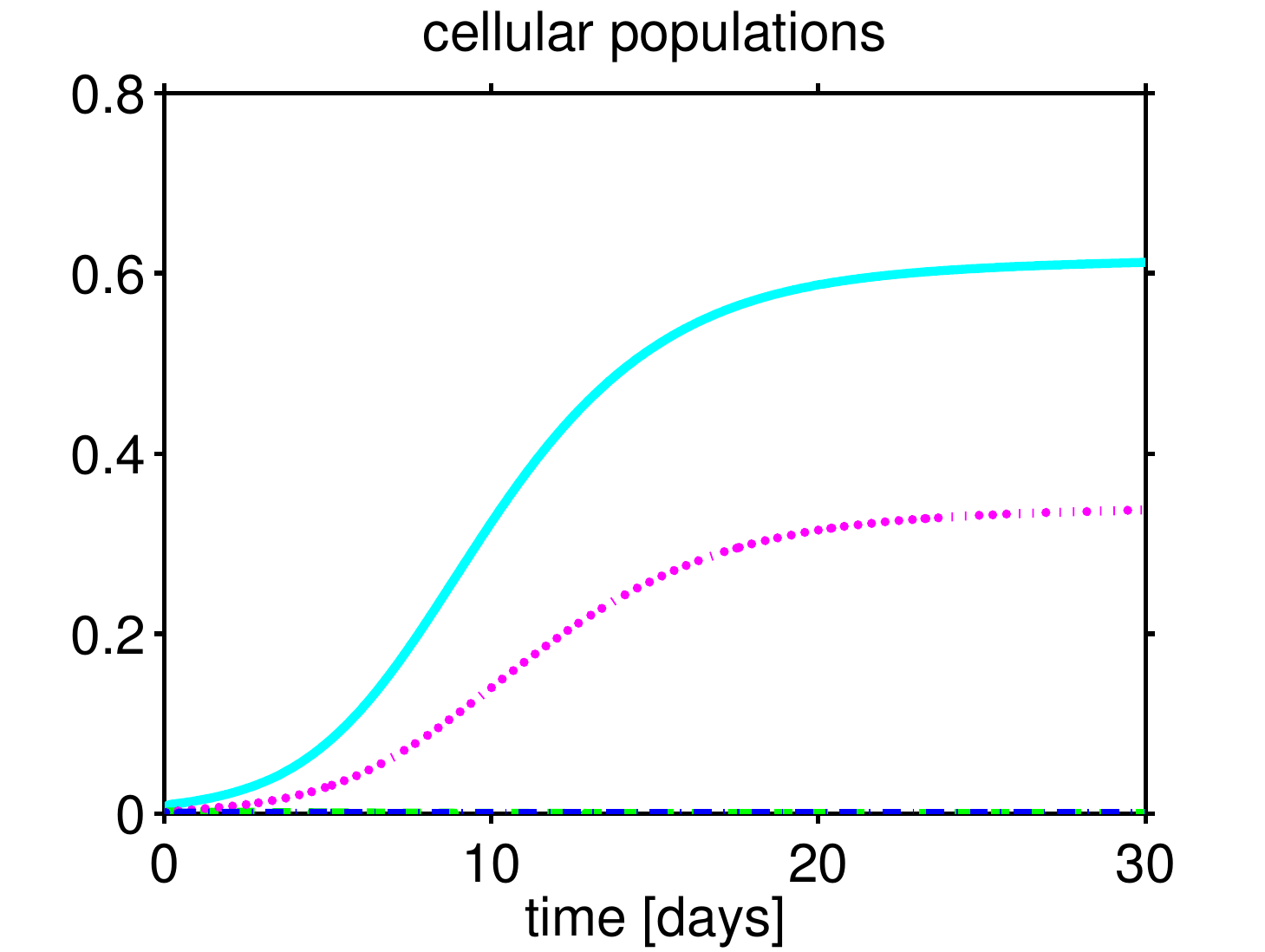}}
\end{minipage}
\caption{Temporal evolution of cellular populations and ECM in the perfused culture for $c_\text{ext}=c_\text{sat}$. Initial condition IC2. Left: $k_\text{g} = k_\text{g1}$. Right: $k_\text{g} = k_\text{g2}$.
\rs{Solid line: $\phi_\n$; dashed line: $\phi_\v$;
dotted line: $\phi_\q$; dash-dot line: $\phi_\ECM$.}
}
\label{fig:Perfused_phi_IC2}
\end{figure}

\begin{figure}[h!]
\begin{minipage}[c]{1\textwidth}
\centering
{\includegraphics[width=0.45\textwidth]{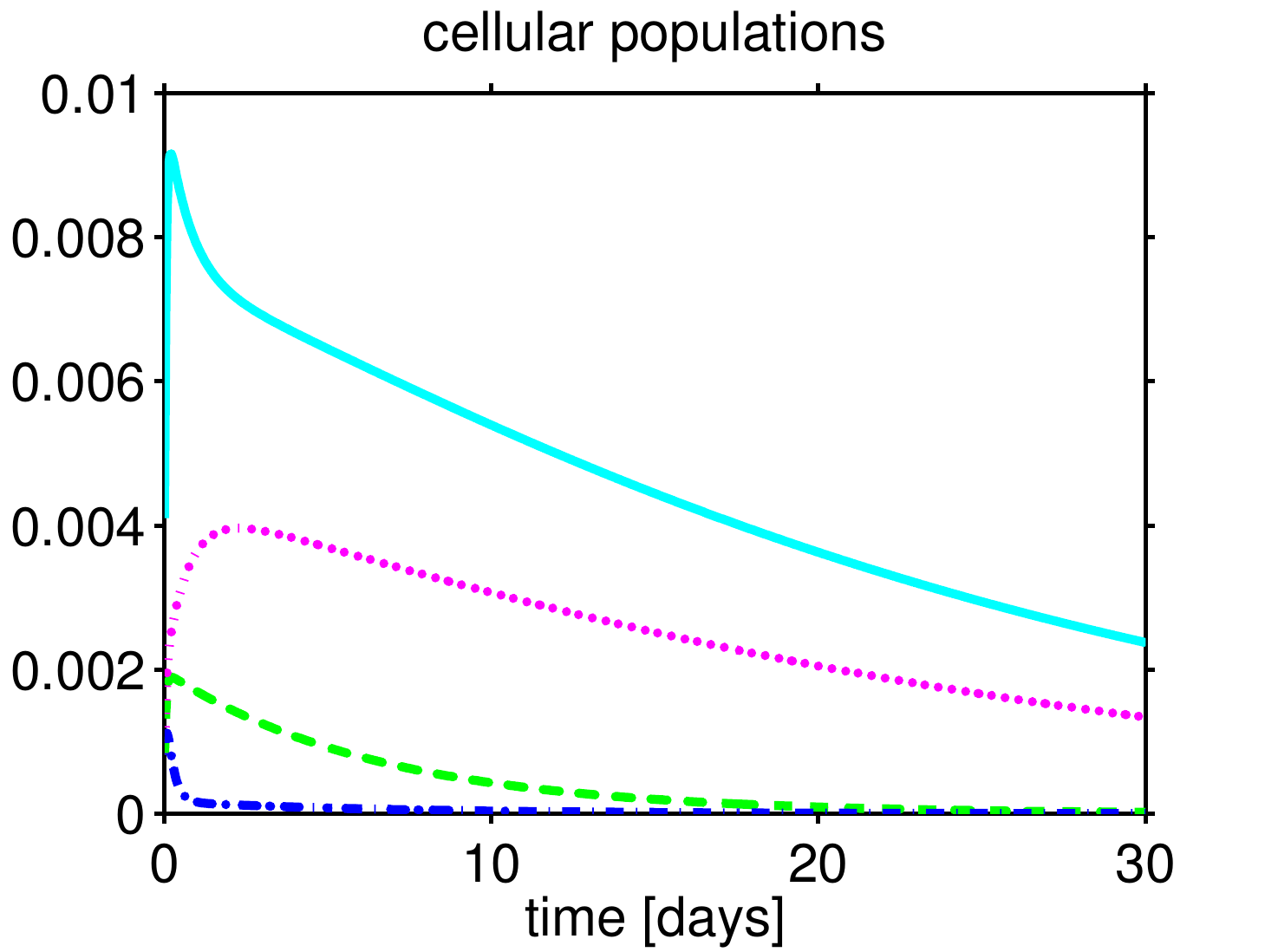}}
{\includegraphics[width=0.45\textwidth]{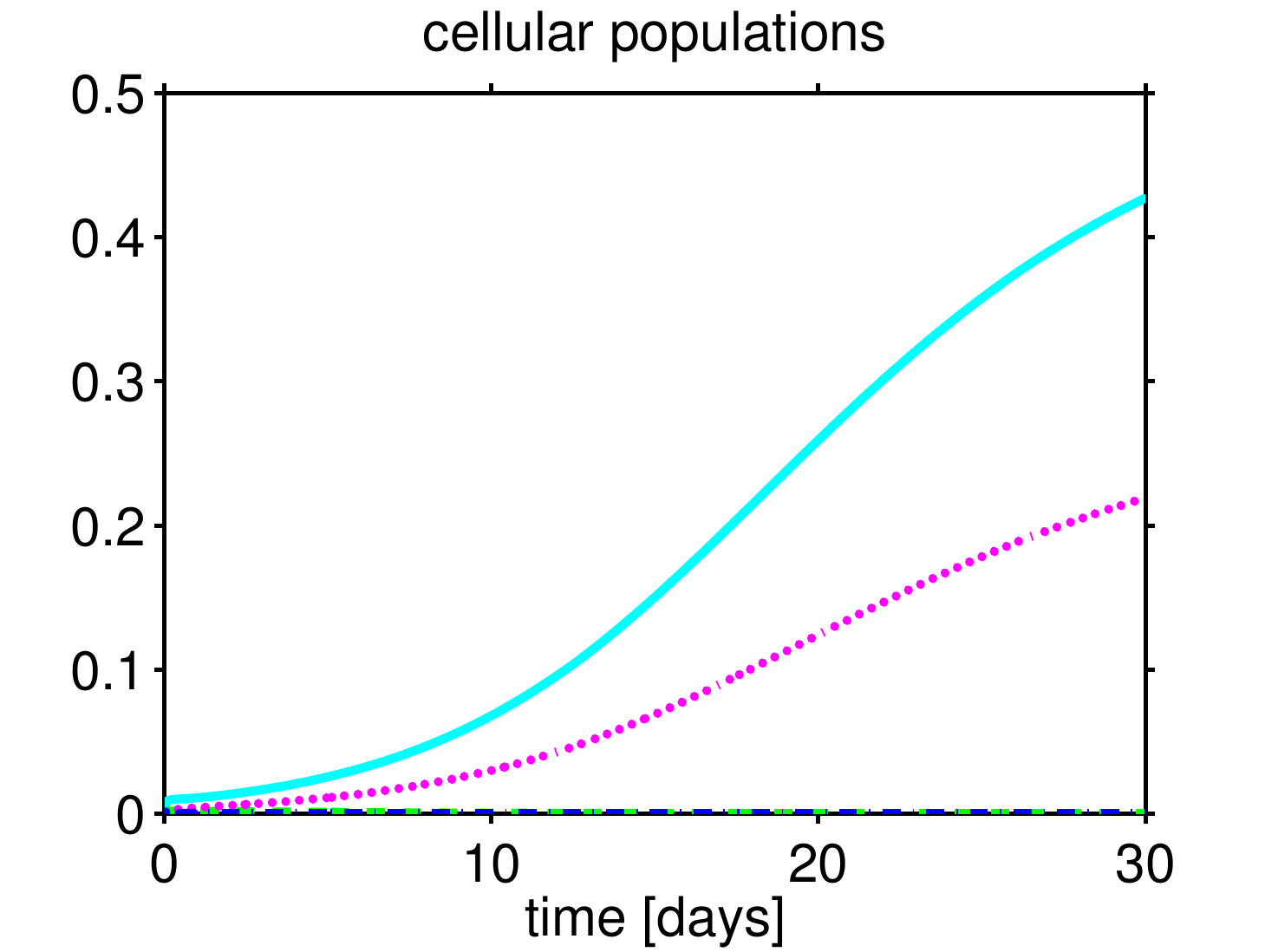}}
\end{minipage}
\caption{Temporal evolution of cellular populations and ECM in the perfused culture for $c_\text{ext}=c_\text{thr}$. Initial condition IC2. Left: $k_\text{g} = k_\text{g1}$. Right: $k_\text{g} = k_\text{g2}$.
\rs{Solid line: $\phi_\n$; dashed line: $\phi_\v$;
dotted line: $\phi_\q$; dash-dot line: $\phi_\ECM$.}
}
\label{fig:Perfused_phi_cthr_IC2}
\end{figure}

Fig.~\ref{fig:Perfused_phi_IC2} shows that at high seeding densities cell behavior does not undergo \textcolor{black}{a significant variation}, 
probably because in the perfused case the conditions for oxygen 
supply improve and oppose the negative consequences of a high 
crowding of cells. For this reason we present below 
only the computational results relative to the solid mixture components,
since the other variables evolve in a very similar way 
as in the initial condition IC1.
As shown in Fig.~\ref{fig:Perfused_phi_cthr_IC2}, the initial high seeding density guarantees the survival
of biomass cellular component even in correspondence of a reduced
level of oxygen concentration.

\section{Discussion of simulation results}
\label{sec:discussion}
In this paper we developed \textcolor{black}{the numerical approximation
of the mathematical mechanobiological model introduced in Part I of this
article.} 
% Our aim was to emphasize} the power of the model to describe 
% \textcolor{black}{in a quantitative manner}
% the biophysical phenomena that simultaneously occur at the
% cellular scale level and that influence the 
% macroscopic response and functioning of in vitro tissue growth.
% As a matter of fact, although 
\textcolor{black}{Despite the intrinsic simplicity of the one-dimensional
formulation, the computational} 
tests described in Sect.~\ref{sec:simulations} constitute 
a useful validation of the
\rs{biophysical} consistency of the mechano-physiological model and 
a \textcolor{black}{preliminary} attempt toward
a qualitative and quantitative prediction on the evolution of the construct components and, consequently, on the formation of the engineered tissue.
\textcolor{black}{Below we address the more significant 
outcomes of the conducted simulations.}

\begin{enumerate} 
\item The illustrated 
numerical results indicate that the in vitro cell cultivation process is strongly sensitive to
variations of \textcolor{black}{(i) the initial seeding density of cells, 
(ii) the value of the maximum growth rate and (iii) 
the mechanical boundary conditions.}
In particular, the amount of seeded cells \textcolor{black}{turns out to be} 
determinant for cell responsiveness: initial cell density
should be high enough to ensure optimal conditions for proliferation, but not so high that grow factors are
rapidly depleted from the medium and the contact inhibition 
phenomenon prevents the formation of new colonies (see Fig.~\ref{fig:Static_cox_IC2}).
Furthermore the value of parameter $k_\text{g}$ influences the long-term behavior of the biomass. As a matter of fact, \textcolor{black}{model simulations
indicate that if $k_\text{g}$ is smaller than} 
the reference value $k_\text{g,0}$, cell metabolic activity and ECM synthesis significantly decrease in the
cultured construct and are completely exhausted at about 10 days of culture (Figs.~\ref{fig:Static_phi} top, \ref{fig:Static_phi_cthr} left).
On the contrary, \textcolor{black}{model predictions show that} 
if the maximum growth rate exceeds the reference value, cell and ECM volumetric fractions increase until convergence to a finite value that represents a stable steady state of the mathematical system (Figs.~\ref{fig:Static_phi} bottom, \ref{fig:Static_phi_cthr} right).
This finding represents a favorable result from the 
\textcolor{black}{experimentalist} point of view, because it predicts the formation,
at the end of the cultivation and under specific conditions, of the bio-artificial texture to be used for replacing damaged tissues,
that constitutes the real aim of \rs{TE.}
A similar objective can be reached by conveniently assigning the mechanical boundary conditions at the interface
between the biomass construct and the interstitial fluid. 
As a matter of fact, \textcolor{black}{model results indicate that} 
when the biomass is stimulated by both external fluid
velocity and pressure, even if $k_\text{g}$ is 
tuned on a under-threshold value, the amount of cells and ECM in the construct
remains considerable until the end of the simulation (see Fig.~\ref{fig:Perfused_phi} top panel).
This outcome reinforces the notion that mechanical stimulation in perfused cultures may promote
chondrogenesis and ECM production~\cite{Raimondi2006a,Raimondi2008,Chung2007}.
Actually, in order to achieve this optimal result, 
nutrient concentration at the fluid-biomass interface should
not fall under a critical level
otherwise cell functionality could be rapidly reduced until
cell apoptosis (see Fig.~\ref{fig:Perfused_phi_cthr} left panel).
%The oxygen level contained in the fluid represents another crucial parameter
%that influences system response,
% limiting the analysis to the dynamical system that describes the conservation of mass of the solid mixture components.
\item The behavior of the solid mixture components is in 
excellent agreement with the experimental trends obtained by cultivation of engineered tissues in bioreactors. In particular, 
the temporal evolution of construct cellularity and ECM content,
especially for an under-threshold value of $k_\text{g}$, 
agree with experimental results shown in several papers~\cite{Novakovic,Obradovic, Chung2007,Davisson}.

\item The characterization of the (an)isotropicity 
of the biomass intrinsic stress state through the equivalent parameter $\xi$ introduced in~\eqref{eq:csi} demonstrated to be a successful strategy 
to model the mechanical regulation of culture progression and to link
the mechanisms occurring at the micro-scale level to the macroscopic functioning of the growing tissue \rs{(see~\cite{Nava2015}).}
As a matter of fact, \textcolor{black}{model predictions show 
that the parameter $\xi$ is an effective 
indicator of the propagation of the isotropic and anisotropic 
waves within the construct and allows 
an easy and immediate identification of the adhesion
mechanisms developing at the single cell-level, that, accordingly, drive the evolution of the volumetric fraction $\phi_\v$ and $\phi_\n$ respectively.}
\end{enumerate}

\section{\rs{Future perspectives}}
\label{sec:conclusions}
% 
% In this Part II of our research work on mechanobiological models
% for Tissue Engineering applications we have numerically solved the
% equation system proposed in Part I in a one-dimensional
% geometrical setting to allow an easy preliminary validation
% of the formulation.
% 
% \null
% 
% Extensive simulation tests have outlined a generally sound response 
% of the computational model with respect to biophysical conjectures.
% 
% \null
% 
% In particular, numerical results indicate that the 
% in vitro cell cultivation process is strongly sensitive to
% variations of (i) the initial seeding density of cells, 
% (ii) the value of the maximum growth rate and (iii) 
% the mechanical boundary conditions.
% 
% \null
Below, we mention several future steps that we intend to take 
in the prosecution of this promising research activity.
\begin{enumerate}
\item 
A sensitivity analysis to provide the bio-scientist an indication of the critical values of the 
parameters that significantly perturb the cell cultivation fate.
The analysis will be the object 
of a subsequent paper in which we will determine the 
homogeneous stable steady states of the model.
\item 
The introduction of a visco-elastic component in the constitutive law for the 
total stress with the purpose of damping out the propagation of 
the anisotropy/isotropy waves, as recently proposed in~\cite{bociu}.
\item 
The inclusion of other mixture constituents, such 
as proteoglycan and collagen as done in~\cite{Klisch}.
% 
% or model the temporal and spatial proteoglycan accumulation~\cite{KleinSah}, or
% focus on the development of cartilaginous tissue around an isolated chondrocyte~\cite{Please}.
\item 
The extension of the computational algorithms to treat
a fully three-dimensional representation of the scaffold pore 
to allow a deeper model validation 
against previous existing simulation results and experimental data 
(see, e.g.,~\cite{Raimondi2005,Cioffi2008,Raimondi2}).
\end{enumerate}

% \textcolor{black}{The good response of the proposed model in the 1D 
% geometrical setting is a promising motivation to extend the
% simulation analysis conducted in Sect.~\ref{sec:simulations} 
% to a fully three-dimensional representation of the scaffold pore.
% This would allow us to perform a deeper model validation 
% against previous existing simulation results and experimental data 
% (see, e.g.,~\cite{Raimondi2005,Cioffi2008,Raimondi2}).}
% 
% 

\begin{table}[h!]
\small{
\centering
\begin{tabular}{|l|l|l|l|l|}
\hline
Symbol&Definition&Value&Units&Reference\\
\hline
$c_0$&$O_2$ concentration for $t=0$&$5\times10^{-6}$&
$\unit{g \; cm^{-3}}$ &this work\\
\hline
$c_\text{sat}$&$O_2$ saturation concentration&$6.4\times10^{-6}$& $\unit{g\;cm^{-3}}$&\cite{Causin}\\
\hline
$c_\text{thr}$&$O_2$ threshold concentration&$1.6\times10^{-6}$&
$\unit{g\; cm^{-3}}$&\cite{MaraNava}\\
\hline
$c_\text{apo}$&$O_2$ apoptosis concentration&$3.2\times10^{-7}$&
$\unit{g\; cm^{-3}}$&\cite{MaraNava}\\
\hline
%$D_\text{c}$&$O_2$ diffusion coefficient&$3.2\times10^{-5}$&$\unit{cm^2s^{-1}}$&\cite{Sacco}\\
$K_\mathrm{eq}$&$O_2$ local mass equilibrium coefficient&$0.1$&-&\cite{CAIM_Causin}\\
\hline
$D_{\mathrm{c},\s}$&$O_2$ diffusivity in the solid phase&$0.75\times10^{-5}$&$\unit{cm^2\; s^{-1}}$&\cite{CAIM_Causin}\\
\hline
$D_{\mathrm{c},\fl}$&$O_2$ diffusivity in the fluid phase&$1\times10^{-5}$&$\unit{cm^2\; s^{-1}}$&\cite{CAIM_Causin}\\
% \hline
% $\mathbf{v}_\text{in}(t=0)$&inlet velocity of perfusion fluid&$50\times10^{-4}$&$cm\; s^{-1}$&\cite{Causin}\\
\hline
$V_b$ &inlet velocity of perfusion fluid&$50\times10^{-4}$&
$\unit{cm\; s^{-1}}$&\cite{Causin}\\
\hline
$T_b$ &stress due to perfusion fluid&$100$&$\unit{mPa}$&\cite{Causin}\\
\hline
% \rs{$D_\fl$????}&fluid diffusion coefficient&?&$cm^2\; s^{-1}$&\cite{}\\
% \hline
$\mu_\fl$&fluid dynamic viscosity at $20^\circ$ C&
$1.002 \cdot 10^{-2}$ & $\unit{g \; cm^{-1}\;s^{-1}}$&\cite{incropera1990fundamentals}\\
\hline
$R_\n=R_\v$&$O_2$ consumption rate for n/v-cells&$3.9\times10^{-8}$&
$\unit{g\; (cm^3\; s)^{-1}}$&\cite{Sacco}\\
\hline
$R_\q$&$O_2$ consumption rate for q-cells &$10^{-8}$&
$\unit{g\; (cm^3\; s)^{-1}}$&this work\\
\hline
$K_{1/2}$&$O_2$ half saturation constant&$3.2\times10^{-6}$&
$\unit{g\; cm^{-3}}$&\cite{Sacco}\\
\hline
$\beta_{{\rm A} \rightarrow {\rm B}}$ & transition rate from state A
to state B &$10^{-5}$&$\unit{s^{-1}}$&this work\\
\hline
%
% $\beta_{\n\rightarrow\q}$&$\n\rightarrow\q$ transition rate&$10^{-5}$&$s^{-1}$&this work\\
% \hline
% $\beta_{\n\rightarrow\v}$&$\n\rightarrow\v$ transition rate&$10^{-5}$&$s^{-1}$&this work\\
% \hline
% $\beta_{\q\rightarrow\n}$&$\q\rightarrow\n$ transition rate&$10^{-8}$&$s^{-1}$&this work\\
% \hline
% $\beta_{\v\rightarrow\q}$&$\v\rightarrow\q$ transition rate&$10^{-6}$&$s^{-1}$&this work\\
% \hline
% $\beta_{\q\rightarrow\v}$&$\q\rightarrow\v$ transition rate&$10^{-8}$&$s^{-1}$&this work\\
% \hline
% $\beta_{\v\rightarrow\n}$&$\v\rightarrow\n$ transition rate&$10^{-5}$&$s^{-1}$&this work\\
% \hline
$k_\text{apo}$&apoptosis transition rate&$3.858\times10^{-7}$&$\unit{s^{-1}}$&\cite{Sacco}\\
\hline
$k_\text{qui}$&quiescence transition rate&$3.858\times10^{-7}$&
$\unit{s^{-1}}$&this work\\
\hline
$k_\text{deg}$&ECM degradation rate&$7.7\times10^{-7}$&
$\unit{s^{-1}}$
&\cite{Please}\\
\hline
$k_\text{g0}$ & maximum specific cell growth rate&$5.8 \times10^{-6}$&$\unit{s^{-1}}$&\cite{Sacco}\\
\hline
$k_\text{g1}$ & "low" specific cell growth rate&
$1 \times 10^{-7}$&$\unit{s^{-1}}$&this work\\
\hline
$k_\text{g2}$ & "high" specific cell growth rate&
$1 \times 10^{-5}$&$\unit{s^{-1}}$&this work\\
\hline
$E$&expansion coefficient&$20$&&\cite{Causin}\\
\hline
$k_\text{GAG}$&GAG synthesis rate&$8.61\times10^{-11}$&$\unit{cm^6\;(cell\;s\;g)^{-1}}$&\cite{Causin}\\
\hline
$K_\text{sat}$&Monod saturation constant&$1.927\times10^{-6}$&$\unit{g\;cm^{-3}}$&\cite{Cheng2}\\
\hline
$D_\eta$&cells and ECM diffusion coefficient&$1\times10^{-9}$&
$\unit{cm^2\; s^{-1}}$&this work\\
\hline
$\lambda_\eta$&cells and ECM Lam\'{e}'s parameter&$5.1937\times10^{3}$&$\unit{dyne\,cm^{-3}}$&\cite{DeSantis}\\
\hline
$\mu_\eta$&cells and ECM Lam\'{e}'s parameter&$1.8248\times10^{3}$&$\unit{dyne\,cm^{-3}}$&\cite{DeSantis}\\
\hline
$\phi_{\ECM,\mathrm{max}}$&maximum ECM volume fraction&$0.1$&&this work\\
\hline
$R_\text{cell}$&cell radius&$5 \times 10^{-4}$&$\unit{cm}$&this work\\
\hline
$V_\text{cell}$&cell volume&$5.236\times10^{-10}$&
$\unit{cm^3}$&this work\\
\hline
$\tau_\text{m}$&mitotic characteristic time&$172800$&$\unit{s}$
&\cite{Sengers}\\
\hline
$K_\text{ref}$&reference permeability&$1.67 \cdot 10^{-5}$
&$\unit{cm^{3} s\; g^{-1}}$& this work\\
\hline
\end{tabular}
}
\caption{Numerical values of model parameters used in the
simulation tests.}
\label{parameters}
\end{table}

\section*{Acknowledgements}
\rs{Chiara Lelli was partially supported by
Grant 5 per Mille Junior 2009
"Computational Models for Heterogeneous Media. Application
to Micro Scale Analysis of Tissue-Engineered Constructs"
CUPD41J10000490001 Politecnico di Milano.
Manuela T. Raimondi has received
funding from the European Research Council (ERC) under the European
Union's Horizon 2020 research and innovation program (Grant
Agreement No. 646990-NICHOID). These results reflect only the
author's view and the agency is not responsible for any use that may be
made of the information contained.}

\bibliographystyle{plain}      % mathematics and physical sciences
\bibliography{PAPERMB}   % name your BibTeX data base

\begin{thebibliography}{10}

\bibitem{bociu}
L.~Bociu, G.~Guidoboni, R.~Sacco, and J.~T. Webster.
\newblock Analysis of nonlinear poro-elastic and poro-visco-elastic models.
\newblock Under Review in Archive for Rational Mechanics and Analysis, 2015.

\bibitem{BrezziFortin1991}
F.~Brezzi and M.~Fortin.
\newblock {\em Mixed and Hybrid Finite Element Methods}.
\newblock Springer Verlag, New York, 1991.

\bibitem{Brooks1982}
A.~N. Brooks and T.~J.R. Hughes.
\newblock Streamline {U}pwind/{P}etrov-{G}alerkin formulations for convection
  dominated flows with particular emphasis on the incompressible
  {N}avier-{S}tokes equations.
\newblock {\em Computer Methods in Applied Mechanics and Engineering},
  32(1\u20133):199 -- 259, 1982.

\bibitem{CAIM_Causin}
P.~Causin and R.~Sacco.
\newblock A computational model for biomass growth simulation in tissue
  engineering.
\newblock {\em Communications in Applied and Industrial Mathematics},
  2(1):1--20, 2011.

\bibitem{Causin}
P.~Causin, R.~Sacco, and M.~Verri.
\newblock A multiscale approach in the computational modeling of the
  biophysical environment in artificial cartilage tissue regeneration.
\newblock {\em Biomechanics and Modeling in Mechanobiology}, 12(4):763--780,
  2013.

\bibitem{Cheng2}
G.~Cheng, P.~Markenscoff, and K.~Zygourakis.
\newblock A {3D} hybrid model for tissue growth: The interplay between cell
  population and mass transport dynamics.
\newblock {\em Biophys J.}, 97:401--414, 2009.

\bibitem{Chung2007}
C.~A. Chung, C.~W. Chen, C.~P. Chen, and C.~S. Tseng.
\newblock Enhancement of cell growth in tissue-engineering constructs under
  direct perfusion: Modeling and simulation.
\newblock {\em Biotechnology and Bioengineering}, 97.

\bibitem{Cioffi2008}
M.~Cioffi, J.~Kueffer, S.~Stroebel, G.~Dubini, I.~Martin, and D.~Wendt.
\newblock Computational evaluation of oxygen and shear stress distributions in
  {3D} perfusion culture systems: Macro-scale and micro-structured models.
\newblock {\em Journal of Biomechanics}, 41(14):2918--2925, 2008.

\bibitem{Davisson}
T.~Davisson, R.~L. Sah, and A.~Patcliffe.
\newblock Perfusion increases cell content and matrix synthesis in chondrocyte
  three-dimensional cultures.
\newblock {\em Tissue Eng}, 8:807--816, 2002.

\bibitem{incropera1990fundamentals}
F.P. Incropera and D.P. DeWitt.
\newblock {\em Fundamentals of heat and mass transfer}.
\newblock Wiley, 1990.

\bibitem{Issa}
R.~I. Issa, B.~Engebretson, L.~Rustom, P.~S. McFetridge, and V.~I. Sikavitsas.
\newblock The effect of cell seeding density on the cellular and mechanical
  properties of a mechanostimulated tissue-engineered tendon.
\newblock {\em Tissue Engineering: Part A}, 17:1479--1487, 2011.

\bibitem{Klisch}
S.~M. Klisch, R.~L. Sah, and A.~Hoger.
\newblock A cartilage growth mixture model for infinitesimal strains: solutions
  of boundary-value problems related to in vitro growth experiments.
\newblock {\em Biomech. Model. Mechanobiol.}, 3:209--223, 2005.

\bibitem{Raimondi2}
M.~Lagan\`{a} and M.~T. Raimondi.
\newblock A miniaturized, optically accessible bioreactor for systematic {3D}
  tissue engineering research.
\newblock {\em Biomed. Microdevices}, 14:225--234, 2012.

\bibitem{MaraNava}
A.~Mara and M.~Nava.
\newblock Modellizzazione multifisica del processo di rigenerazione tessutale
  all'interno di un bioreattore perfuso, 2011.
\newblock Master Thesis, Politecnico di Milano.

\bibitem{Nava2015}
M.~M. Nava, R.~Fedele, and M.~T. Raimondi.
\newblock Computational prediction of strain-dependent diffusion of
  transcription factors through the cell nucleus.
\newblock {\em Biomechanics and Modeling in Mechanobiology}, pages 1--11, 2015.

\bibitem{Nava2012}
M.M. Nava, M.T. Raimondi, and R.~Pietrabissa.
\newblock Controlling self-renewal and differentiation of stem cells via
  mechanical cues.
\newblock {\em Journal of Biomedicine and Biotechnology}, 2012:12, 2012.

\bibitem{Nikolaev}
N.I. Nikolaev, B.~Obradovic, H.K. Versteeg, G.~Lemon, and D.J. Williams.
\newblock A validated model of gag deposition, cell distribution, and growth of
  tissue engineered cartilage cultured in a rotating bioreactor.
\newblock {\em Biotechnol Bioeng.}, 105(4):842--853, 2010.

\bibitem{Obradovic}
B.~Obradovic, J.~H. Meldon, L.~E. Freed, and G.~Vunjak-Novakovic.
\newblock Glycosaminoglycan deposition in engineered cartilage: Experiments and
  mathematical model.
\newblock {\em AIChE J.}, 46 (9):1860--1871, 2000.

\bibitem{QuarteroniValli94}
A.~Quarteroni and A.~Valli.
\newblock {\em {N}umerical {A}pproximation of {P}artial {D}ifferential
  {E}quations}.
\newblock Springer-Verlag, New York, Berlin, 1994.

\bibitem{Raimondi2006}
M.~T. Raimondi.
\newblock Engineered tissue as a model to study cell and tissue function from a
  biophysical perspective.
\newblock {\em Current Drug Discovery Technologies}, 3(4):245--268, 2006.

\bibitem{Raimondi2004}
M.~T. Raimondi, F.~Boschetti, L.~Falcone, F.~Migliavacca, A.~Remuzzi, and
  G.~Dubini.
\newblock The effect of media perfusion on three-dimensional cultures of human
  chondrocytes: integration of experimental and computational approaches.
\newblock {\em Biorheology}, 41:401--10, 2004.

\bibitem{Raimondi2005}
M.~T. Raimondi, F.~Boschetti, F.~Migliavacca, M.~Cioffi, and G.~Dubini.
\newblock Micro fluid dynamics in three-dimensional engineered cell systems in
  bioreactors.
\newblock In N.~Ashammakhi and R.~L. Reis, editors, {\em Topics in Tissue
  Engineering}, volume~2, chapter~9. 2005.

\bibitem{Raimondi2008}
M.~T. Raimondi, G.~Candiani, M.~Cabras, M.~Cioffi, M.~Lagan\`{a}, M.~Moretti,
  and R.~Pietrabissa.
\newblock Engineered cartilage constructs subject to very low regimens of
  interstitial perfusion.
\newblock {\em Biorheology}, 45:471--8, 2008.

\bibitem{Raimondi2006a}
M.T. Raimondi, M.~Moretti, M.~Cioffi, C.~Giordano, F.~Boschetti, K.~Lagan\`a,
  and R.~Pietrabissa.
\newblock The effect of hydrodynamic shear on {3D} engineered chondrocyte
  systems subject to direct perfusion.
\newblock {\em Biorheology}, 43(3-4):215--222, 2006.

\bibitem{RoosStynesTobiska1996}
H.~G. Roos, M.~Stynes, and L.~Tobiska.
\newblock {\em Numerical methods for singularly perturbed differential
  equations}.
\newblock Springer-Verlag, Berlin Heidelberg, 1996.

\bibitem{SaccoABB2014}
R.~Sacco, L.~Carichino, C.~de~Falco, M.~Verri, F.~Agostini, and T.~Gradinger.
\newblock A multiscale thermo-fluid computational model for a two-phase cooling
  system.
\newblock {\em Computer Methods in Applied Mechanics and Engineering},
  282(0):239 -- 268, 2014.

\bibitem{Sacco}
R.~Sacco, P.~Causin, P.~Zunino, and M.~T. Raimondi.
\newblock A multiphysics/multiscale {2D} numerical simulation of scaffold-based
  cartilage regeneration under interstitial perfusion in a bioreactor.
\newblock {\em Biomech Model Mechanobiol}, 10(4):577--589, 2011.

\bibitem{DeSantis}
G.~De Santis, A.~B. Lennon, F.~Boschetti, B.~Verhegghe, P.~Verdonck, and P.~J.
  Prendergast.
\newblock How can cells sense the elasticity of a substrate? an analysis using
  a cell tendegrity model.
\newblock {\em European cells and Materials}, 22:202--213, 2011.

\bibitem{Sengers}
B.~G. Sengers, C.~W.~J. Oomens, and F.~P.~T Baaijens.
\newblock An integrated finite-element approach to mechanics, transport and
  biosyntheis in tissue engineering.
\newblock {\em Journal of Biomechanical Engineering}, 126(1):82--91, 2004.

\bibitem{ateshian}
M.~A. Soltz and G.~A. Ateshian.
\newblock Experimental verification and theoretical prediction of cartilage
  interstitial fluid pressurization at an impermeable contact interface in
  confined compression.
\newblock {\em J. Biomech.}, 31(10):927--934, 1998.

\bibitem{Please}
A.~J. Trewenack, C.~P. Please, and K.~A. Landman.
\newblock A continuum model for the development of tissue-engineered cartilage
  around a chondrocyte.
\newblock {\em Mathematical Medicine and Biology}, 26:241--262, 2009.

\bibitem{Novakovic}
G.~Vunjak-Novakovic, B.~Obradovic, I.~Martin, P.~M. Bursac, R.~Langer, and
  L.~E. Freed.
\newblock Dynamic cell seeding of polymer scaffolds for cartilage tissue
  engineering.
\newblock {\em Biotechnology Progress}, 14(2):193--202, 1998.

\end{thebibliography}

\end{document}